\input amstex
\input xy
\xyoption{all}
\documentstyle{amsppt}
\document

\magnification 1000

\def\hbar{{t}}

\def\gen{\frak{g}}

\def\ten{\frak{t}}

\def\ben{\frak{b}}

\def\ien{\frak{i}}

\def\ken{\frak{k}}

\def\Fen{\frak{F}}

\def\len{\frak{l}}

\def\nen{\frak{n}}

\def\Nen{\frak{N}}

\def\Men{\frak{M}}

\def\Den{\frak{D}}

\def\Fen{\frak{F}}

\def\Xen{\frak{X}}

\def\Xen{\frak{X}}

\def\a{{\alpha}}
\def\g{{\gamma}}
\def\o{{\omega}}
\def\l{{\lambda}}
\def\b{{\beta}}

\def\qcb{{{qc}}}
\def\Lotimes{{\buildrel L\over\otimes}}

\def\[b{{\pmb [}}
\def\]b{{\pmb ]}}
\def\Ocb{{\pmb{\Oc}}}

\def\Dcb{{\pmb{\Dc}}}
\def\DG{{\pmb\Dc}^G}
\def\CG{{\pmb\Cc}^G}

\def\gb{{\bold g}}

\def\rb{{\bold r}}

\def\xb{{\bold x}}

\def\Ab{{\bold A}}
\def\Xb{{\bold X}}
\def\Bb{{\bold B}}

\def\Dcb{{\pmb{\Dc}}}

\def\It{{\pmb I}}
\def\Tt{{\pmb T}}
\def\Gt{{\pmb G}}

\def\Hb{{\bold H}}

\def\Rb{{\bold R}}

\def\Tb{{\bold T}}

\def\Xb{{\bold X}}
\def\Yb{{\bold Y}}

\def\c{{\roman c}}

\def\K{{\pmb {\roman K}}}

\def\alphav{{\check\alpha}}

\def\Spec{{\roman{Spec}}}

\def\Ext{\roman{Ext}}
\def\Tor{\text{Tor}}
\def\Hom{{\roman{Hom}}}
\def\Coh{{\roman{Coh}}}

\def\dim{{\roman{dim}}}

\def\Ind{{\roman{Ind}}}

\def\End{{\roman{End}}}
\def\on{\ {\roman{on}}\ }

\def\Id{{\roman{id}}}
\def\Irr{{\roman{Irr}}}

\def\op{{{\circ}}}

\def\AA{{\Bbb A}}

\def\CC{{\Bbb C}}

\def\NN{{\Bbb N}}

\def\ZZ{{\Bbb Z}}

\def\Ac{{\Cal A}}

\def\Cc{{\Cal C}}
\def\Dc{{\Cal D}}
\def\Ec{{\Cal E}}
\def\Fc{{\Cal F}}

\def\Gc{{\Cal G}}
\def\Hc{{\Cal H}}

\def\Ic{{\Cal I}}

\def\Lc{{\Cal L}}

\def\Oc{{\Cal O}}
\def\Pc{{\Cal P}}

\def\Sc{{\Cal S}}
\def\Tc{{\Cal T}}

\def\Sch{{\pmb{\Cal{S}ch}}}
\def\Isch{{\pmb{\Cal{I}sch}}}
\def\Pro{{\pmb{\Cal{P}ro}}}
\def\Ind{{\pmb{\Cal{I}nd}}}
\def\Coh{{\pmb{\Cal{C}oh}}}

\def\Qcoh{{\pmb{\Cal{Q}coh}}}
\def\Pcoh{{\pmb{\Cal{P}coh}}}
\def\Alg{{\pmb{\Cal{A}lg}}}
\def\Set{{\pmb{\Cal{S}et}}}
\def\Group{{\pmb{\Cal{G}rp}}}

\def\Spa{{\pmb{\Cal{S}pace}}}
\def\Acb{{\pmb{\Ac}}}
\def\Ccb{{\pmb{\Cc}}}

\def\Xc{{\Cal X}}

\def\rot{{\roman{rot}}}

\def\ft{{\roman{ft}}}

\def\and{{\text{and}}}
\def\mod{{\roman{mod}}}
\def\cent{{\roman{cen}}}

\def\qua{{\roman{qua}}}
\def\nil{{{\roman{nil}}}}

\def\ds{\displaystyle}

\def\bbox{{\sqcap \hskip-6.5pt \sqcup}}                 
\def\qed{\hfill{$\bbox$}}                              
\overfullrule=0pt                                    

\def\7dag{{{\!\!\!\!\!\!\!\dag}}}
\def\6dag{{{\!\!\!\!\!\!\dag}}}
\def\5dag{{{\!\!\!\!\!\dag}}}
\def\4dag{{{\!\!\!\!\dag}}}
\def\3dag{{{\!\!\!\dag}}}
\def\2dag{{{\!\!\dag}}}
\def\1dag{{{\!\dag}}}

\def\pro{{\roman{lim}}}
\def\ind{{\roman{colim}}}

\def\la{{\langle}}
\def\ra{{\rangle}}

%
\def\Vsquare#1{\alphaox{\Square{$#1$}}\kern-\Thickness}

\nologo

\topmatter
\title Double affine Hecke algebras and affine flag manifolds, I
\endtitle
\author M. Varagnolo, E. Vasserot\endauthor
\address D\'epartement de Math\'ematiques,
Universit\'e de Cergy-Pontoise, 2 av. A. Chauvin, BP 222, 95302
Cergy-Pontoise Cedex, France, Fax : 01 34 25 66 45\endaddress \email
michela.varagnolo\@math.u-cergy.fr\endemail
\address D\'epartement de Math\'ematiques,
Universit\'e Paris 7, 175 rue du Chevaleret, 75013 Paris, France,
Fax : 01 44 27 78 18
\endaddress
\email vasserot\@math.jussieu.fr\endemail
\endtopmatter
\document

\head \endhead

\head 0.~Introduction\endhead

This paper is the first of a series of papers reviewing the
geometric construction of the double affine Hecke algebra via affine
flag manifolds. The aim of this work is
to explain the main results in \cite{V}, \cite{VV}, but also to give
a simpler approach to some of them, and to give the proof of some
`folklore' related statements whose proofs are not available in the
published literature. This work should be therefore be viewed as a
companion to loc.~cit., and is by no means a logically independent
treatment of the theory from the very begining. In order that the
length of each paper remains reasonable, we have split the whole
exposition into several parts. This one concerns the most basic
facts of the theory : the geometric construction of the double
affine Hecke algebra via the equivariant, algebraic
K-theory and the classification of the
simple modules of the category $\Ocb$ of the double affine Hecke algebra.
It is our hope that by providing a detailed
explanation of some of the difficult aspects of the foundations,
this theory will be better understood by a wider audience.

This paper contains three chapters. The first one is a reminder on
$\Oc$-modules over non Noetherian schemes and over ind-schemes. The
second one deals with affine flag manifolds. The last chapter
concerns the classification of simple modules in the category $\Ocb$
of the double affine Hecke algebra. Let us review these parts in
more details.

In the second chapter of the paper we use two different versions of
the affine flag manifold. The first one is an ind-scheme of
ind-finite type, while the second one is a pro-smooth, coherent,
separated, non Noetherian, and non quasi-compact scheme. Thus, in
the first chapter we recall some basic fact on $\Oc$-modules over
coherent schemes, pro-schemes, and ind-schemes. The first section is
a reminder on pro-objects and ind-objects in an arbitrary category.
We give the definition of direct and inverse 2-limits of categories.
Next we recall the definition of the K-homology of a scheme.
We'll use non quasi-compact non Noetherian scheme.
Also, it is convenient to consider a quite general setting involving
unbounded derived categories, pseudo-coherent complexes and perfect
complexes. Fortunately, since all the schemes we'll consider are
coherent the definition of the K-theory remains quite close to the
usual one. To simplify the exposition it is convenient to introduce
the derived direct image of a morphism of non Noetherian schemes,
its derived inverse image, and the derived tensor product in the
unbounded derived categories of $\Oc$-modules. Finally we consider
the special case of pro-schemes (compact schemes, pro-smooth
schemes, etc) and of ind-schemes. They are important tools in this
work. This section finishes with equivariant $\Oc$-modules and some
basic tools in equivariant K-theory (induction, reduction of the
group action, the Thom isomorphism, and the Thomason concentration
theorem).

The second chapter begins with the definition of the affine flag
manifold, which is an ind-scheme of ind-finite type, and with the
definition of the Kashiwara affine flag manifold, which is a non
quasi-compact coherent scheme. This leads us in Section 2.3.6 to the
definition of an associative multiplication on a group of
equivariant K-theory $\K^\It(\Nen)$. Here $\Nen$ is an ind-scheme
which can be regarded as the affine analogue of the Steinberg
variety for reductive groups. Then, in section 2.4.1, we define an
affine analogue of the concentration map for convolution rings in
K-theory used in \cite{CG}. It is a ring homomorphism relates
$\K^\It(\Nen)$ to the K-theory of the fixed points subset for a
torus action. This concentration map is new, and it simplifies the
proofs in \cite{V}. The double affine Hecke algebra is introduced in
section 2.5.1 and its geometric realization is proved in Theorem
2.5.6. We use here an approach similar to the one in \cite{BFM},
where a degenerate version of
the double affine Hecke algebra is constructed
geometrically.
Compare also \cite{GG}, where the regular representation of the
double affine Hecke algebra is constructed geometrically.
The proof we give uses a reduction to the
fixed points of a torus acting on the affine analogue of the
Steinberg variety, and the concentration map in K-theory.

The third chapter is a review of the classification of the simple modules in the
category $\Ocb$ of the double affine Hecke algebra. The main theorem was proved
in \cite{V}. The proof we give here is simpler than in loc.~cit.~because
it uses the concentration map. The first section contains generalities
on convolution algebras in the cohomology of schemes of infinite type which are
locally of finite type. The proof of the classification is given in the
second section.

\vskip3mm

\head Contents\endhead

\item{1.}Schemes and ind-schemes
\itemitem{1.1.}Categories and Grothendieck groups
\itemitem{1.2.}K-theory of schemes
\itemitem{1.3.}K-theory of ind-coherent ind-schemes
\itemitem{1.4.}Group actions on ind-schemes
\itemitem{1.5.}Equivariant K-theory of ind-schemes
\item{2.}Affine flag manifolds
\itemitem{2.1.}Notation relative to the loop group
\itemitem{2.2.}Reminder on the affine flag manifold
\itemitem{2.3.}K-theory and the affine flag manifold
\itemitem{2.4.}Complements on the concentration in K-theory
\itemitem{2.5.}Double affine Hecke algebras
\item{3.}Classification of the simple admissible modules of the
Double affine Hecke algebra
\itemitem{3.1.}Constructible sheaves and convolution algebras
\itemitem{3.2.}Simple modules in the category $\Ocb$

\vskip3mm

\head 1.~Schemes and ind-schemes\endhead

\subhead 1.1.~Categories and Grothendieck groups\endsubhead

\subhead 1.1.1.~Ind-objects and pro-objects in a category\endsubhead
A standard reference for the material in this section is \cite{SGA4,
sec.~8}, \cite{KS1}, \cite{KS2}.

Let $\Set$ be the category of sets. Given a category $\Ccb$ let
$\Ccb^\op$ be the opposite category. The category $\Ccb^\wedge$ of
{\it presheaves over $\Ccb$} is the category  of functors
$\Ccb^\op\to\Set$. We'll abbreviate $\Ccb^\vee$ for the category
$((\Ccb^\op)^\wedge)^\op$. Yoneda's lemma yields fully faithful
functors
$$\gathered
\Ccb\to\Ccb^\wedge,\ X\mapsto\Hom_\Ccb(\cdot,X),\quad
\Ccb\to\Ccb^\vee,\ X\mapsto\Hom_\Ccb(X,\cdot).\cr
\endgathered$$

Let $\Acb$ be a category and $\a\mapsto X_\a$ be a functor
$\Acb\to\Ccb$ or $\Acb^\op\to\Ccb$ (also called a {\it system} in
$\Ccb$ indexed by $\Acb$ or $\Acb^\op$). Let $\ind_\a X_\a$ and
$\pro_\a X_\a$ denote the colimit or the limit of this system
whenever it is well-defined. If the category $\Acb$ is small or
filtrant the colimit and the limit are said to be {\it small} or
{\it filtrant}. A poset $\Acb=(A,\leqslant)$ may be viewed as a
category, with $A$ as the set of objects and with a morphism
$\a\to\b$ whenever $\a\leqslant\b$. A direct set is a poset $\Acb$
which is filtrant as a category. A {\it direct system} over $\Ccb$
is a functor $\Acb\to\Ccb$ and an {\it inverse system} over $\Ccb$
is a functor $\Acb^\op\to\Ccb$, where $\Acb$ is a direct set. A {\it
direct colimit} (also called {\it inductive limit}) is the colimit
of a direct system. An {\it inverse limit} (also called {\it
projective limit}) is the limit of an inverse system. Both are small
and filtrant.

A {\it complete} or {\it cocomplete} category is one that has all
small limits or all small colimits. A {\it Grothendieck category} is
a cocomplete Abelian category with a generator such that the small
filtrant colimits are exact.

Given a direct system or an inverse system in  $\Ccb$ we define the
following functors
$$\gathered
``\ind_\a" X_\a :\ \Ccb^\op\to\Set,\quad
Y\mapsto\ind_\a\Hom_\Ccb(Y,X_\a),\cr ``\pro_\a" X_\a :\
\Ccb\to\Set,\quad Y\mapsto\ind_\a\Hom_\Ccb(X_\a,Y).
\endgathered$$
The categories of {\it ind-objects} of $\Ccb$ and the category of
{\it pro-objects} of $\Ccb$ are the full subcategory $\Ind(\Ccb)$ of
$\Ccb^\wedge$ and the full subcategory $\Pro(\Ccb)$ of $\Ccb^\vee$
consisting of objects isomorphic to some $``\ind_\a"X_\a$ and
$``\pro_\a" X_\a$ respectively. Note that we have
$\Pro(\Ccb)=\Ind(\Ccb^\op)^\op$. By the Yoneda functor we may look
upon $\Ccb$ as a full subcategory of $\Ind(\Ccb)$ or $\Pro(\Ccb)$.
We'll say that an ind-object or a pro-object is {\it representable}
if it is isomorphic to an object in $\Ccb$. Note that, for each
object $Y$ of $\Ccb$ we have the following formulas, see, e.~g.,
\cite{KS1, sec.~1.11}, \cite{KS2, sec.~2.6}
$$\gathered
\Hom_{\Ind(\Ccb)}(Y,``\ind_\a" X_\a)=
\Hom_{\Ccb^\wedge}(Y,``\ind_\a" X_\a)= \ind_\a\Hom_{\Ccb}(Y,X_\a),
\cr \Hom_{\Pro(\Ccb)}(``\pro_\a" X_\a, Y)=
\Hom_{\Ccb^\vee}(``\pro_\a" X_\a, Y)= \ind_\a\Hom_{\Ccb}(X_\a, Y).
\endgathered$$

\subhead 1.1.2.~Direct and inverse 2-limits\endsubhead Let
$\Acb=(A,\leqslant)$ be a directed set. Given a direct system of
categories $(\Ccb_\a,i_{\a\b}:\Ccb_\a\to\Ccb_\b)$ the {\it 2-limit}
(also called the {\it 2-colimit}) of this system is the category
$\Ccb=2\ind_\a\Ccb_\a$
whose objects are the pairs $(\a,X_\a)$ with $X_\a$ an object of
$\Ccb_\a$. The morphisms are given by
$$\Hom_\Ccb((\a,X_\a),(\b,X_\b))=\ind_{\g\geqslant\a,\b}
\Hom_{\Ccb_\g}(i_{\a\g}(X_\a),i_{\b\g}(X_\b)).$$ Given an inverse
system of categories $(\Ccb_\a,i_{\a\b}:\Ccb_\b\to\Ccb_\a)$ the
{\it 2-limit} of this system is the category
$\Ccb=2\pro_\a\Ccb_\a$
whose objects are the families of objects $X_\a$ of $\Ccb_\a$ and of
isomorphisms $i_{\a\b}(X_\b)\simeq X_\a$ satisfying the obvious
composition rules. The morphisms are defined in the obvious way. See
\cite{W, app.~A} for more details on 2-colimits and 2-limits.

\subhead 1.1.3.~Grothendieck groups and derived
categories\endsubhead Given an Abelian category $\Acb$ let
$\Ccb(\Acb)$ be the category of complexes of objects of $\Acb$ with
differential of degree $+1$ and chain maps as morphisms, let
$\Dcb(\Acb)$ be the corresponding (unbounded) derived category, let
$\Dcb(\Acb)^-$ be the full subcategory of complexes bounded above,
let $\Dcb(\Acb)^b$ be the full subcategory of bounded complexes.
Finally let $\[b\Acb\]b$ the Grothendieck group of $\Acb$.

The Grothendieck group $\[b\pmb\Tc\]b$ of a triangulated
category $\pmb\Tc$ is the quotient of the free Abelian group with
one generator for each object $X$ of $\pmb\Tc$ modulo the
relations $X=X'+X''$ for each distinguished triangle
$$X'\to X\to X''\to X'[1].$$
Here the symbol [1] stands for the shift functor in the triangulated category
$\pmb\Tc$.
Throughout, we'll use the same symbol for an object of
$\pmb\Tc$ and is class in $\[b\pmb\Tc\]b$.

Recall that the Grothendieck group of $\Dcb(\Acb)^b$ is canonically
isomorphic to $\[b\Acb\]b$, and that two quasi-isomorphic complexes
of $\Ccb(\Acb)$ have the same class in $\[b\Acb\]b$.


\vskip1mm

\proclaim{1.1.4.~Proposition} Let $(\Ccb_\a)$ be a direct system of
Abelian categories (resp.~of triangulated categories) and exact
functors. Then the direct 2-limit $\Ccb$ of $(\Ccb_\a)$ is also an
Abelian category (resp.~a triangulated category) and we have a
canonical group isomorphism $\[b\Ccb\]b=\ind_\a
\[b\Ccb_\a\]b.$
\endproclaim

\vskip3mm

\subhead 1.2.~ K-theory of schemes\endsubhead

This section is a recollection of standard results from \cite{SGA6},
\cite{TT} on the K-theory of schemes, possibly of infinite type.

\vskip3mm

\subhead 1.2.1.~Background\endsubhead  For any Abelian category
$\Acb$ a complex in $\Ccb(\Acb)$ is cohomologically bounded if the
cohomology sheaves vanish except for a finite number of them. The
canonical functor yields an equivalence from $\Dcb(\Acb)^b$ to the
full subcategory of $\Dcb(\Acb)$ consisting of cohomologically
bounded complexes \cite{KS1, p.~45}.

A quasi-compact scheme is a scheme that has a finite covering by
affine open subschemes (e.g., a Noetherian scheme or a scheme of
finite type is quasi-compact) and a quasi-separated scheme is a
scheme such that the intersection of any two affine open subschemes
is quasi-compact (e.g., a separated scheme is quasi-separated). More
generally, a scheme homomorphism $f : X\to Y$ is said to be
quasi-compact, resp. quasi-separated, if for every affine open
$U\subset Y$ the inverse image of $U$ is quasi-compact, resp.
quasi-separated. Elementary properties of quasi-compact and
quasi-separated morphisms can be found in \cite{GD, chap.~I,
sec.~6.1}. For instance quasi-compact and quasi-separated morphisms
are stable under composition and pullback, and if $f:X\to Y$ is a
scheme homomorphism with $Y$ quasi-compact and quasi-separated then
$X$ is quasi-compact and separated iff $f$ is quasi-compact and
quasi-separated. {\it Throughout, by the word scheme we'll always
mean a separated $\CC$-scheme and by the word scheme homomorphism
we'll always mean a morphism of separated $\CC$-schemes}. In
particular a scheme homomorphism will always be separated (hence
quasi-separated) \cite{GD, chap.~I, sec.~5.3}.

Given a scheme $X$, the word $\Oc_X$-module will mean a sheaf on the
scheme $X$ which is a sheaf of modules over the sheaf of rings
$\Oc_X$. Unless otherwise stated, modules are left modules. This
applies also to $\Oc_X$-modules. Since $\Oc_X$ is commutative this
specification is indeed irrelevant. Let $\Ocb(X)$ be the Abelian
category of all $\Oc_X$-modules. Given a closed subscheme $Y\subset
X$ let $\Ocb(X\on Y)$ be the full subcategory of $\Oc_X$-modules
supported on $Y$.

Let $\Coh(X)$, $\Qcoh(X)$ be the categories of coherent and
quasi-coherent $\Oc_X$-modules. They are Abelian subcategories of
$\Ocb(X)$ which are stable under extensions.
Quasi-coherent sheaves are preserved by tensor products, by
arbitrary colimits, and by inverse images \cite{GD, chap.~I,
sec.~2.2}. They are well-behaved on quasi-compact (quasi-separated)
schemes : under this assumption quasi-coherent $\Oc_X$-modules are
preserved by direct images and any quasi-coherent $\Oc_X$-module is the
limit of a direct system of finitely presented $\Oc_X$-modules.
Further, if $X$ is quasi-compact (quasi-separated) the category
$\Qcoh(X)$ is a Grothendieck category. In particular for any such
$X$ there are enough injective objects in $\Qcoh(X)$ \cite{GD,
chap.~I, sec.~6.7,~6.9}, \cite{TT, sec.~B.3}. Given a closed
subscheme $Y\subset X$ let $\Coh(X\on Y)$, $\Qcoh(X\on Y)$ be the
full subcategories of sheaves supported on $Y$.

We'll abbreviate $\Ccb(X)=\Ccb(\Ocb(X))$ and
$\Dcb(X)=\Dcb(\Ocb(X))$. Let $\Dcb(X)_\qcb$ be the full subcategory
of $\Dcb(X)$ of complexes of $\Oc_X$-modules with quasi-coherent
cohomology.

\subhead 1.2.2.~Remark\endsubhead
B\"oksted and Neeman proved that if $X$ is quasi-compact
(separated) then the canonical functor is an equivalence
$$\Dcb(\Qcoh(X))\to \Dcb(X)_\qcb\leqno(1.2.1)$$
\cite{BN, cor.~5.5}, \cite{Li,
prop.~3.9.6}. Further, the standard derived functors in 1.2.10-12 below,
evaluated on quasi-coherent sheaves, are the same taken in $\Ocb(X)$
or in $\Qcoh(X)$, see e.g., \cite{TT, cor.~B.9}.
So from now on we'll identify the categories
$\Dcb(\Qcoh(X))$ and $\Dcb(X)_\qcb$.

\vskip3mm

A commutative ring $R$ is {\it coherent} iff it is coherent as a
$R$-module, or, equivalently, if every finitely generated ideal of
$R$ is finitely presented. For instance a Noetherian ring is
coherent, the quotient of a coherent ring by a finitely generated
ideal is a
coherent ring and the localization of a coherent ring is again coherent.

\subhead 1.2.3.~Definitions\endsubhead Let $X$ be any scheme. We say
that

\vskip1mm

$(a)$ $X$ is {\it coherent} if its structure ring $\Oc_X$ is
coherent,

\vskip1mm

$(b)$ $X$ is {\it locally of countable type} if the $\CC$-algebra
$\Oc_X(U)$ is generated by a countable number of elements for any
affine open subset $U\subset X$,

\vskip1mm

$(c)$ a closed subscheme $Y\subset X$ is {\it good} if the ideal of
$Y$ in $\Oc_X(U)$ is finitely generated for any affine open subset
$U\subset X$.

\vskip2mm

If the scheme $X$ is coherent then an $\Oc_X$-module is coherent iff
it is finitely presented, and we have $f^*(\Coh(Y))\subset\Coh(X)$
for any morphism $f:X\to Y$. If $X$ is quasi-compact and coherent
then any quasi-coherent $\Oc_X$-module is the direct colimit of a
system of coherent $\Oc_X$-modules. Finally a good subscheme $Y$ of
a coherent scheme $X$ is again coherent and the direct image of
$\Oc_Y$ is a coherent $\Oc_X$-module. See \cite{EGAI, chap.~0,
sec.~5.3} for details.

\vskip3mm

\subhead 1.2.4.~K-theory of a quasi-compact coherent
scheme\endsubhead For an arbitrary scheme $X$ the K-homology group
(=K-theory) may differ from $[\Coh(X)]$, one reason being that
$\Oc_X$ may not be an object of $\Coh(X)$. Let us recall briefly
some relevant definitions and results concerning pseudo-coherence.
Details can be found in \cite{SGA6, chap.~I}, \cite{TT} and
\cite{Li, sec.~4.3}. We'll assume that $X$ is quasi-compact and
coherent.

\subhead{1.2.5.~Definition-Lemma}\endsubhead $(a)$ A complex of
quasi-coherent $\Oc_X$-modules $\Ec$ is {\it pseudo-coherent} if it
is locally quasi-isomorphic to a bounded above complex of vector
bundles. Since $X$ is coherent, this simply means that $\Ec$ has
coherent cohomology sheaves vanishing in all sufficiently large
degrees \cite{SGA6, cor.~I.3.5(iii)}. In particular any coherent
$\Oc_X$-module is a pseudo-coherent complex.

\vskip1mm

$(b)$ Let $\Pcoh(X)$ be the full subcategory of $\Dcb(X)_\qcb$
consisting of the cohomologically bounded pseudo-coherent complexes.
Given a closed subscheme $Y\subset X$ the full subcategory of
complexes which are acyclic on $X-Y$ is $\Pcoh(X\ \roman{on}\ Y)$.
It is a triangulated category.

\vskip2mm

Note that for a general scheme the equivalence of categories (1.2.1)
does not hold and a pseudo-coherent complex may consist of non
quasi-coherent $\Oc_X$-modules. The K-homology group of the pair
$(X,Y)$ is \cite{SGA6, def.~IV.2.2}
$$\K(X\on Y)=\[b\Pcoh(X\on Y)\]b,\quad\K(X)=\K(X\on X).$$
By 1.2.5 the K-homology groups are well-behaved on quasi-compact coherent
schemes. More precisely we have the following.

\proclaim{1.2.6.~Proposition} Assume that $X$ is quasi-compact and
coherent. We have $\K(X\on Y)=\[b\Coh(X\on Y)\]b$ for any closed
subscheme $Y\subset X$. If $Y\subset X$ is good there is a canonical
isomorphism $\K(Y)\to\K(X\on Y)$.
\endproclaim

If $X$ is coherent but not quasi-compact we define the group $\K(X)$
as follows. Fix a covering $X=\bigcup_wX^w$ by quasi-compact open
subsets. The restrictions yield a inverse system of categories with
a functor
$$\Coh(X)\to 2\pro_w\Coh(X^w).$$
By functoriality of the K-theory we have also an inverse system of
Abelian groups. We define
$$\K(X)=\pro_w\K(X^w)=\pro_w\[b\Coh(X^w)\]b.$$  The group
$\K(X)$ does not depend on the choice of the open covering. It may
be regarded as a completion of the K-homology group of $X$, as
defined in \cite{SGA6}.

\subhead 1.2.7.~Remark\endsubhead Let $X$ be a quasi-compact scheme.
A {\it perfect complex} over $X$ is a complex of quasi-coherent
$\Oc_X$-modules which is locally quasi-isomorphic to a bounded
complex of vector bundles.
The K-cohomology groups of $X$ is
the Grothendieck group of the the full subcategory of $\Dcb(X)_\qcb$
of perfect complexes. We'll not use it.

\subhead 1.2.8.~Basic properties of the K-theory of a coherent
quasi-compact scheme\endsubhead Recall that for any $\Oc_X$-modules
$\Ec$, $\Fc$ the {\it local hypertor} is the $\Oc_X$-module $\Tc
or_i^{\Oc_X}(\Ec,\Fc)$ whose stalk at a point $x$ is
$\Tor_i^{\Oc_{X,x}}(\Ec_x,\Fc_x)$.

\vskip2mm

\subhead 1.2.9.~Definitions\endsubhead $(a)$ An $\Oc_X$-module $\Ec$
has a {\it finite tor-dimension} if there is an integer $n$ such
that $\Tc or_i^{\Oc_X}(\Ec,\Fc)=0$ for each $i>n$ and each
$\Fc\in\Qcoh(X)$.

\vskip1mm

$(b)$ A scheme homomorphism $f:X\to Y$ has {\it finite
tor-dimension} if there is an integer $n$ such that $\Tc
or_i^{\Oc_Y}(\Oc_X,\Ec)=0$ for each $i>n$ and each $\Ec\in\Qcoh(Y)$.
Equivalently, $f$ has finite tor-dimension if there is an integer
$n$ such that for each $x\in X$ there is an exact sequence of
$\Oc_{Y,f(x)}$-modules
$$0\to P_n\to P_{n-1}\to\dots\to P_0\to\Oc_{X,x}\to 0$$
with $P_i$ flat over $\Oc_{Y,f(x)}$.

\vskip1mm

$(c)$ A scheme $X$ {\it satisfies the Poincar\'e duality} if any
quasi-coherent $\Oc_X$-module has a finite tor-dimension.

\vskip3mm

Poincar\'e duality is a local property. Note that, since taking a
local hypertor commutes with direct colimits \cite{EGAIII,
prop.~6.5.6}, a coherent scheme satisfies the Poincar\'e duality iff
any coherent $\Oc_X$-module has a finite tor-dimension. Note that if
$X$ satisfies the Poincar\'e duality then any cohomologically
bounded pseudo-coherent complex is perfect \cite{TT, thm.~3.21}.

Now, let us recall a few basic properties of direct/inverse image of
complexes of $\Oc$-modules. We'll use derived functors in the
unbounded derived category of $\Oc$-modules. This simplifies the
exposition. Their definition and properties can be found in
\cite{Li}. To simplify, in Sections 1.2.10 to 1.2.15 we'll also
assume that all schemes are quasi-compact and coherent. Therefore,
all morphisms will also be quasi-compact.

\subhead 1.2.10.~Derived inverse image\endsubhead For any morphism
$f:Z\to X$ the inverse image functor $Lf^*$ maps $\Dcb(X)_\qcb$,
$\Dcb(X)^-$ into $\Dcb(Z)_\qcb$, $\Dcb(Z)^-$ respectively. It
preserves pseudo-coherent and perfect complexes \cite{Li,
prop.~3.9.1}, \cite{TT, sec.~2.5.1}. Further if $\Ec$ is a
pseudo-coherent complex then the complex $Lf^*(\Ec)$ is
cohomologically bounded if the map $f$ has finite tor-dimension.
Under this assumption, for each closed subscheme $Y\subset X$ the
functor $Lf^*$ yields a group homomorphism
$$Lf^*:\K(X\on Y)\to \K(Z\on f^{-1}(Y)).$$
If the schemes $X$, $Z$ satisfy the Poincar\'e duality then $f$ has a
finite tor-dimension.

\subhead 1.2.11.~Derived tensor product\endsubhead Let $\otimes_X$
denote the tensor product of $\Oc$-modules on any scheme $X$. The
standard theory of the derived tensor product of $\Oc$-modules
applies to complexes in $\Dcb(X)^-$, see e.g., \cite{Hr, p.93}.
Following Spaltenstein \cite{Sp} we can  extend the theory to
arbitrary complexes in $\Dcb(X)$, see also \cite{Li, sec.~2.5}. This
yields a functor
$${\buildrel L\over\otimes}_X:\Dcb(X)\times\Dcb(X)\to
\Dcb(X)$$ which maps $\Dcb(X)_\qcb\times\Dcb(X)_\qcb$,
$\Dcb(X)^-\times\Dcb(X)^-$ to $\Dcb(X)_\qcb$, $\Dcb(X)^-$
respectively. It preserves pseudo-coherent complexes \cite{TT,
sec.~2.5.1}. If $\Ec$, $\Fc$ are pseudo-coherent complexes their
derived tensor product is cohomologically bounded if either $\Ec$ is
perfect or $\Fc$ is perfect \cite{TT, sec.~3.15}. Recall that if $X$
satisfies the Poincar\'e duality then any cohomologically bounded
pseudo-coherent complex is perfect. Under this assumption, for each
closed subschemes $Y,Z\subset X$ there is a (derived) tensor product
$$\buildrel L\over\otimes_X : \K(X\ \roman{on}\ Y)\times \K(X\
\roman{on}\ Z)\to \K(X\ \roman{on}\ Y\cap Z).$$ Given a map $f$ as
in 1.2.10 there is a functorial isomorphism in $\Dcb(X)$ [Li,
prop.~3.2.4]
$$Lf^*(\Ec\Lotimes_X\Fc)=Lf^*(\Ec)\Lotimes_ZLf^*(\Fc)$$
We'll refer to this relation by saying that the derived tensor
product commutes with $Lf^*$.

\subhead 1.2.12.~Derived direct image\endsubhead For any $f:X\to Z$
the direct image functor $Rf_*$ is right adjoint to $Lf^*$ and it
maps $\Dcb(X)_\qcb$, $\Dcb(X)^b$ into $\Dcb(Z)_\qcb$, $\Dcb(Z)^b$
respectively \cite{Li, prop.~3.9.2}. We say that the map $f$ is {\it
pseudo-coherent} if it factors, locally on $X$, as $f=p\circ i$
where $i$ is a closed embedding with $i_*\Oc_X$ coherent and $p$ is
smooth. Kiehl's finiteness theorem insures that if $f$ is proper and
pseudo-coherent then $Rf_*$ preserves pseudo-coherent complexes
\cite{Li, cor.~4.3.3.2}. Therefore if $f$ is proper and
pseudo-coherent, for any closed subscheme $Y\subset X$, the functor
$Rf_*$ yields a group homomorphism
$$Rf_*:\K(X\on Y)\to \K(Z\on f(Y)).$$

\subhead 1.2.13.~Example\endsubhead A good embedding is proper and
pseudo-coherent. In this case we have indeed an exact functor
$f_*:\Coh(X)\to\Coh(Z)$. It yields the isomorphism $\K(X)\to\K(Z\on
X)$ in 1.2.5. Note that a closed embedding $X\subset Z$ with
$Z=\AA^\NN$ and $X$ of finite type is not pseudo-coherent. Here
$\AA^\NN=\Spec(\CC[x_i;i\in\NN])$ is a coherent scheme.

\subhead 1.2.14.~Projection formula\endsubhead For any $f:X\to Z$
there is a canonical isomorphism called the {\it projection formula}
\cite{Li, prop.~3.9.4}
$$Rf_*(\Ec\Lotimes_{X}Lf^*(\Fc))=Rf_*(\Ec)\Lotimes_Z\Fc,\quad
\forall \Ec\in\Dcb(X)_\qcb,\,\Fc\in\Dcb(Z)_\qcb.$$

\subhead 1.2.15.~Base change\endsubhead Consider the following
Cartesian square
$${\xymatrix{ X'\ar[d]_g&Y'\ar[l]_-{f'}  \ar[d]_-{g'} \\ X&\ar[l]_{f}  Y.}}
$$ Assume that it is {\it tor-independent}, i.e., assume that we have
$$\Tor_i^{\Oc_{X,x}}(\Oc_{X',x'},\Oc_{Y,y})=0,
\quad\forall i>0,\,\forall x\in X,\forall x'\in X', \forall y\in
Y,\,x=g(x')=f(y).$$ Then we have a functorial base-change
isomorphism \cite{Li, thm.~3.10.3}
$$Lg^*Rf_*(\Ec)\simeq Rf'_*L(g')^*(\Ec),\quad\forall\Ec\in\Dcb(Y)_\qcb.$$

\vskip3mm

\subhead 1.2.16.~Compact schemes\endsubhead A simple way to produce
quasi-compact schemes of infinite type is to use pro-schemes. Let us
explain this.

\subhead 1.2.17.~Lemma-Definition\endsubhead A scheme is {\it
compact} if it is the limit of an inverse system of finite type
schemes with affine morphisms. A scheme is compact iff it is
quasi-compact \cite{TT, thm.~C.9}.

\subhead 1.2.18.~Remarks\endsubhead Let $X$ be a compact scheme and
$(X^\a)$ be an inverse system of schemes
as above. Then the canonical maps $p_\a:X\to
X^\a$ are affine. Further the following hold.

\vskip1mm

$(a)$ If $\Fc$ is a coherent $\Oc_{X}$-module there is an $\a$
and a coherent $\Oc_{X^\a}$-module $\Fc^\a$ such that
$\Fc=(p_\a)^*(\Fc^\a)$. Given two coherent $\Oc_X$-modules $\Fc$,
$\Gc$ and two coherent $\Oc_{X^\a}$-modules $\Fc^\a$, $\Gc^\a$ as
above we set $\Fc^\b=(p_{\a\b})^*(\Fc^\a)$ and
$\Gc^\b=(p_{\a\b})^*(\Gc^\a)$ for each $\b\geqslant\a$. Then we have
\cite{TT, sec.~C4}, \cite{EGAIV, sec.~8.5}
$$\Hom_{\Oc_X}(\Fc,\Gc)=\ind_{\b\geqslant\a}\Hom_{\Oc_{X^\b}}(\Fc^\b,\Gc^\b).$$

$(b)$ If $f:Y\to X$ is a scheme finitely presented over $X$ then
there is an $\a\in A$ and a finitely presented $f_\a:Y^\a\to X^\a$
such that \cite{TT, sec.~C.3}
$$\gathered
f=f_\a\times\Id,\quad Y=Y^\a\times_{X^\a}X=\pro_{\b}Y^\b,\quad
Y^\b=Y^\a\times_{X^\a}X^\b,\quad\b\geqslant\a.
\endgathered$$

\vskip3mm

\subhead 1.2.19.~Definition\endsubhead A compact scheme
$X=``\pro_\a"X^\a$ {\it satisfies the property $(S)$} if $A=\NN$ and
$(X^\a)_{\a\in A}$ is an inverse system of smooth schemes of finite type
with smooth affine morphisms. A scheme is {\it pro-smooth} if it is
covered by a finite number of open subsets satisfying $(S)$.

\proclaim{1.2.20.~Proposition} A pro-smooth scheme is quasi-compact,
coherent, and it satisfies the Poincar\'e duality.
\endproclaim

\noindent{\sl Proof :} A pro-smooth scheme $X$ is coherent by
\cite{K1, prop.~1.1.6}. Let us prove that $X$ satisfies the
Poincar\'e duality. Let $\Fc\in\Coh(X)$. We must prove that
$\Tor^{\Oc_{X,x}}_i(\Fc_x,\Ec_x)=0$ for each $i\gg 0$, each $x\in
X$, and each $\Ec\in\Qcoh(X)$.
Since the question is local around $x$ we can assume that $X$ is a
compact scheme satisfying the property $(S)$. By 1.2.18$(a)$ there
is an $\a\in A$ and $\Fc^\a\in\Coh(X^\a)$ such that
$\Fc=(p_\a)^*(\Fc^\a)$. Write again $x=p_\a(x)$. Since $X^\a$ is
smooth of finite type the $\Oc_{X^\a,x}$-module $\Fc^\a_x$ has
finite tor-dimension. Since the map $p_\a$ is affine we have
$$\Tor_i^{\Oc_{X,x}}(\Fc_x,\Ec_x)=
\Tor_i^{\Oc_{X^\a,x}}(\Fc^\a_x,(p_\a)_*(\Ec_x)).$$
Since $(p_\a)_*(\Ec)$ is quasi-coherent and the scheme $X^\a$ is
smooth of finite type, and since taking the Tor's commutes with
direct colimits, the rhs vanishes for large $i$.

\qed

\vskip3mm

\subhead 1.2.21.~Remarks\endsubhead $(a)$ Let $\Sch$ be the category
of schemes and $\Sch^\ft$ be the full subcategory of schemes of
finite type. The category of compact schemes can be identified with
a full subcategory in $\Pro(\Sch^\ft)$ via the assignment $\pro_\a
X^\a\mapsto``\pro_\a"X^\a$. From now on we'll omit the quotation
marks for compact schemes.

\vskip1mm

$(b)$ Let $X$ be a quasi-compact coherent scheme and $Y\subset X$ be
a good subscheme. Since $X$ is a compact scheme we can fix an
inverse system of finite type schemes $(X^\a)$ with affine morphisms
$p_{\a\b}:X^\b\to X^\a$ such that $X=\pro_\a X^\a$. Since the scheme
$X$ is coherent and since $Y$ is a good subscheme, the inclusion
$Y\subset X$ is finitely presented. Thus, by 1.2.18$(b)$ there is an
$\a\in A$ and a closed subscheme $Y^\a\subset X^\a$ such that
$Y=p_\a^{-1}(Y^\a)$. Setting $Y^\b=p_{\a\b}^{-1}(Y^\a)$ for each
$\b\geqslant\a$ we get a direct system of categories $\Coh(X^\a\on
Y^\a)$ with functors $(p_{\a\b})^*$.
Now, assume that the pro-object $X=\pro_\a X^\a$ satisfies the property $(S)$.
The pull-back by the canonical
map $p_\a:X\to X^\a$ yields an equivalence of categories
\cite{EGAIV, thm.~8.5.2}
$$2\ind_\a\Coh(X^\a\on Y^\a)\to\Coh(X\on Y),$$
and we have a group isomorphism
$$\ind_\a \K(X^\a\on Y^\a)=\K(X\on Y).$$
See also \cite{SGA6, sec.~IV.3.2.2}, \cite{TT, prop.~3.20}.

\subhead 1.2.22.~Pro-finite-dimensional vector bundles\endsubhead An
important particular case of compact schemes is given by
pro-finite-dimensional vector bundles.

\subhead 1.2.23.~Definition\endsubhead
A {\it pro-finite-dimensional vector bundle}
$\pi : X\to Y$ is a scheme homomorphism
which is represented as the inverse limit of a system
of vector bundle $\pi_n:X^n\to Y$ (of finite rank) with $n$ an integer
$\geqslant 0$,
such that the morphism $X^m\to X^n$ is a vector bundle homomorphism
for each $m\geqslant n$.

\proclaim{1.2.24.~Proposition} A pro-finite-dimensional vector
bundle $\pi:X\to Y$ is flat. If $Y$ is compact then $X$ is compact.
If $Y$ is pro-smooth then $X$ is pro-smooth. If $Y$ is coherent then
$X$ is coherent.
\endproclaim

\noindent{\sl Proof :} The first claim is obvious.  By 1.2.18$(b)$
any vector bundle over $Y$ is, locally over $Y$, pulled-back from a
vector bundle over some $Y^\a$ where $(Y^\a)$ is a inverse system as
in 1.2.19. This implies the second and the third claim. The last one
follows from \cite{K1, prop.~1.1.6}.

\qed

%

\vskip3mm

\subhead 1.3.~ K-theory of ind-coherent ind-schemes\endsubhead

\subhead 1.3.1.~Spaces and ind-schemes\endsubhead
Let $\Alg$ be the category of
associative, commutative $\CC$-algebras with 1.
The category of {\it spaces} is the category
$\Spa$ of functors $\Alg\to\Set$.
By Yoneda's lemma $\Sch$ can be considered
as a full subcategory in the category $\Sch^\wedge$ of presheaves on
$\Sch$. It can be as well realized as a full subcategory in $\Spa$
via the functor
$$\Sch\to\Spa,\quad X\mapsto \Hom_\Sch(\Spec(\cdot),X).$$
By a {\it subspace} we mean a subfunctor. A subspace $Y\subset X$ is
said to be {\it closed, open} if for every scheme $Z$ and every
$Z\to X$ the subspace $Z\times_XY\subset Z$ is a closed, open
subscheme.

\subhead 1.3.2.~Definitions\endsubhead  $(a)$ An {\it ind-scheme} is
an ind-object $X$ of $\Sch$ represented as $X=``\ind_\a" X_\a$ where
$A=\NN$ and $(X_\a)_{\a\in A}$ is a direct system of quasi-compact
schemes with closed embeddings $i_{\a\b}:X_\a\to X_\b$ for each
$\a\leqslant\b$.

\vskip1mm

$(b)$ A {\it closed ind-subscheme} $Y$ of the ind-scheme $X$ is a closed
subspace of $X$.
An {\it open ind-subscheme} $Y$ of the ind-scheme $X$ is an ind-scheme which is
an open subspace of $X$.

\vskip3mm

Since direct colimits exist in the category $\Spa$ we may regard
$\Isch$ as a full subcategory of $\Spa$. Hence, to unburden the
notation we'll omit the quotation marks for ind-schemes.

\vskip3mm


\subhead{1.3.3.~Remarks}\endsubhead
$(a)$ A closed subscheme of an ind-scheme is always quasi-compact.
\vskip1mm

$(b)$ We may consider ind-objects of $\Sch$ which are representented
by a direct system of non quasi-compact schemes $X_\a$ with closed
embeddings. To avoid any confusion we'll call them ind$'$-schemes.
\vskip1mm

$(c)$ Given a closed ind-subscheme $Y\subset X$, for each $\a\in A$
the closed immersion $X_\a\subset X$ yields a closed subscheme
$Y_\a=X_\a\times_XY\subset X_\a$.
Further the closed immersion $X_\a\subset X_\b$, $\a\leqslant\b$,
factors to a closed immersion $Y_\a\subset Y_\b$.
The ind-scheme $Y$ is represented as
$Y=\ind_\a Y_\a.$
\vskip1mm

$(d)$ For each ind-scheme $X$ and each quasi-compact scheme $Y$ we
have
$$\Hom_{\Isch}(Y,X)=
\ind_\a\Hom_{\Sch}(Y,X_\a). $$

\vskip3mm

\subhead 1.3.4.~Definitions\endsubhead $(a)$ An ind-scheme $X$ is
{\it ind-proper} or of {\it ind-finite type} if it can be
represented as the direct colimit of system of proper schemes or of
finite type schemes respectively with closed embeddings.

\vskip1mm

$(b)$ An ind-scheme $X$ is {\it ind-coherent} if it can be
represented as the direct colimit of a system of coherent
quasi-compact schemes with good embeddings.

\vskip3mm

\subhead 1.3.5.~Coherent and quasi-coherent $\Oc$-modules over
ind-coherent ind-schemes\endsubhead Let $X$ be an ind-coherent
ind-scheme and $Y\subset X$ be a closed ind-subscheme. Given $X_\a$,
$Y_\a$ as in 1.3.3$(c)$ we have a direct system of Abelian
categories $\Coh(X_\a\on Y_\a)$ with exact functors $(i_{\a\b})_*$.
We define the following Abelian categories
$$
\Coh(X\on Y)=2\ind_\a\Coh(X_\a\on Y_\a),\quad\Coh(X)=\Coh(X\on X).$$
These categories do not depend on the direct system $(X_\a)$ up to
canonical equivalences. An object of $\Coh(X)$ is called a coherent
$\Oc_X$-module.

We define also quasi-coherent $\Oc_X$-modules in the following way
\cite{BD, sec.~7.11.3}, \cite{D, sec.~6.3.2}. We have a inverse
system of categories $\Qcoh(X_\a)$ with functors $(i_{\a\b})^*$. We
set
$$\Qcoh(X)=2\pro_\a\Qcoh(X_\a).$$
The category $\Qcoh(X)$ is a tensor category, but it need not be
Abelian. It is independent on the choice of the system $(X_\a)$ up
to canonical equivalences of categories. A quasi-coherent
$\Oc_X$-module can be regarded as a rule that assigns to each scheme
$Z$ with a morphism $Z\to X$ a quasi-coherent $\Oc_Z$-module
$\Ec_Z$, and to each scheme homomorphism $f:W\to Z$ an isomorphism
$f^*\Ec_Z\simeq\Ec_W$ satisfying the obvious composition rules.

Finally we define the Grothendieck group of the pair $(X,Y)$ by
$$\K(X\on
Y)=\[b\Coh(X\on Y)\]b.$$ Note that we have $\K(X\on Y)=\ind_\a
\K(X_\a\on Y_\a)$ where $\K(X_\a\on Y_\a)=\[b\Coh(X_\a\on Y_\a)\]b$
for each $\a$.

\subhead 1.3.6.~Remarks\endsubhead
$(a)$
There is another notion of quasi-coherent $\Oc_X$-modules on an ind-scheme,
called $\Oc^!_X$-modules in \cite{BD, sec.~7.11.4}.
They form an Abelian category. We'll not need this.

\vskip1mm

$(b)$ Any morphism of ind-coherent ind-schemes $f:X\to Y$ yields a functor
$f^*:\Qcoh(Y)\to\Qcoh(X)$. If $f$ is an open embedding the base change
yields an exact functor
$f^*:\Coh(Y)\to\Coh(X)$ and a group homomorphism
$f^*:\K(Y)\to\K(X).$

\vskip3mm

\subhead 1.3.7.~Definition\endsubhead
Let $X$ be an ind-scheme. A closed ind-subscheme $Y\subset X$ is
{\it good} if for every scheme $Z\to X$ the closed subscheme
$Z\times_XY\subset Z$ is good.

\vskip3mm

Note that if $X$ is ind-coherent and $Y\subset X$ is a good
ind-subscheme then $Y$ is again an ind-coherent ind-scheme. If
$f:Y\to X$ is an ind-proper homomorphism of ind-schemes of
ind-finite type, or a good ind-subscheme of an ind-coherent
ind-scheme, then there is a functor $f_*:\Coh(Y)\to \Coh(X)$ and a
group homomorphism $\quad Rf_*:\K(Y)\to \K(X).$

\vskip3mm

\subhead 1.4.~Group actions on ind-schemes\endsubhead

\subhead 1.4.1.~Ind-groups and group-schemes\endsubhead Let $\Group$
be the category of groups. A group-scheme is a scheme representing a functor
$\Alg\to\Group$. An
ind-group is an ind-scheme representing a functor $\Alg\to\Group$.

\subhead 1.4.2.~Definition\endsubhead We abbreviate {\it linear
group} for linear algebraic group.
A {\it pro-linear group} $G$ is a
compact, affine, group-scheme which is represented as the inverse
limit of a system of linear groups $G=\pro_{n} G^{n}$ with $n$ any integer
$\geqslant 0$,
such that the morphism $G^m\to G^n$ is a group-scheme homomorphism
for each $m\geqslant n$.

\subhead 1.4.3.~Examples\endsubhead Let $G$ be a linear group. For
each $\CC$-algebra $R$ the set of $R$-points of $G$ is
$G(R)=\Hom_\Sch(\Spec(R),G).$

\vskip1mm

$(a)$ The algebraic group $G(\CC[\varpi]/(\varpi^n))$ represents the
functor $R\mapsto G(R[\varpi]/(\varpi^{n}))$. The functor $R\mapsto
G(R[[\varpi]])$ is represented by a group-scheme, denoted by
$K=G(\CC[[\varpi]])$. The group-scheme $K$ is a pro-linear group,
since it is the limit of the inverse system of linear groups
$G(\CC[\varpi]/(\varpi^n))$ with $n\geqslant 0$.

\vskip1mm

$(b)$ The functor  $R\mapsto G(R[\varpi^{-1}])$ is represented by an
ind-group, denoted by $G(\CC[\varpi^{-1}])$.

\vskip1mm

$(c)$ The functor  $R\mapsto G(R((\varpi)))$ is represented by an
ind-group, denoted by $G(\CC((\varpi)))$.

\vskip3mm

Throughout we'll use the same symbol for an ind-scheme $X$ and the
set of $\CC$-points $X(\CC)$. For instance the symbol $K$ will
denote both the functor above and the group of
$\CC[[\varpi]]$-points of the linear group $G$.

\vskip3mm

\subhead 1.4.4.~Group actions on an ind-scheme\endsubhead Let $G$ be
an ind-group and $X$ be an ind-scheme. We'll say that $G$ acts on
$X$ if there is a morphism of functors $G\times X\to X$ satisfying
the obvious composition rules. A {\it $G$-equivariant ind-scheme} is
an ind-scheme with a (given) $G$-action. We'll abbreviate
$G$-ind-scheme for $G$-equivariant ind-scheme. We'll also call
ind-$G$-scheme a $G$-ind-scheme which is represented as the direct
colimit of a system of quasi-compact $G$-schemes $(X_\a)$ as in
1.3.2.

\subhead 1.4.5.~Definition\endsubhead Let $G=\pro_n G^n$ be a
pro-linear group.

$(a)$ A (compact) $G$-scheme $X$ is {\it admissible} if it is represented
as the inverse limit of a system of $G$-schemes of finite type with affine
morphisms $(X^\a)$ such that, for each $\a$, the $G$-action on
$X^\a$ factors through a $G^{n}$-action if $n\geqslant n_\a$
for some integer $n_\a$.

\vskip1mm

$(b)$ A morphism of admissible $G$-schemes $f:X\to Y$ is {\it
admissible} if there are inverse systems of $(X^\a)$, $(Y^\a)$ as
above such that $f$ is the limit of a morphism of systems of
$G$-schemes $(f^\a):(X^\a)\to (Y^\a)$ and the following square
is Cartesian for each $\a\leqslant \b$
$$\xymatrix{X^\b\ar[r]^{f_\b}\ar[d]&Y^\b\ar[d]\cr
X^\a\ar[r]^{f_\a}&Y^\a.
}$$

\vskip1mm

$(c)$ An ind-$G$-scheme $X$ is {\it admissible} if it is the direct
colimit of a system of compact admissible $G$-schemes with admissible closed
embeddings.

\subhead 1.4.6.~Remarks\endsubhead
$(a)$ Let $X$ be a $G$-torsor, with $G=\pro_nG^n$ a pro-linear group. For each
$n$ let $G_n$ be the kernel of the canonical morphism $G\to G^n$. If
the quotient scheme $X/G$ is of finite type then the $G$-scheme $X$
is admissible. Indeed $X$ is the inverse limit of the system of
$G$-schemes $(X/G_n)$, and the $G$-action on $X/G_n$ factors through
a $G^n$-action.

\vskip1mm

$(b)$ If $f:Y\to X$ is a finitely presented morphism of $G$-schemes
with $X$ admissible then $Y$ and $f$ are also admissible \cite{TT,
sec.~C.3}.

\vskip3mm

\subhead 1.5.~Equivariant K-theory of ind-schemes\endsubhead

To simplify, in this section we'll assume that all schemes are
quasi-compact.

\subhead 1.5.1.~Equivariant quasi-coherent $\Oc$-modules over a
scheme\endsubhead
Let $G$ be a group-scheme and $X$ be a
$G$-scheme. Let $a,p:G\times X\to X$ be the action and the obvious
projection.

\subhead 1.5.2.~Definition\endsubhead A {\it $G$-equivariant
quasi-coherent $\Oc_X$-module} is a quasi-coherent $\Oc_X$-module
$\Ec$ with an isomorphism $\theta:a^*(\Ec)\to p^*(\Ec)$. The obvious
cocycle condition is to hold. Let $\Qcoh^{G}(X)$ be the category of
$G$-equivariant quasi-coherent $\Oc_X$-modules. Given a closed
subset $Y\subset X$ we define the category $\Qcoh^G(X\on Y)$ of
$G$-equivariant quasi-coherent $\Oc_X$-modules supported on $Y$ in
the obvious way.

\vskip2mm

%

The category $\Qcoh^{G}(X)$ is Abelian. The forgetful functor
$$for:\Qcoh^{G}(X)\to\Qcoh(X)$$ is exact and it {\it reflects
exactness}, i.e., whenever a sequence in $\Qcoh^G(X)$ is exact in
$\Qcoh(X)$ it is also exact in $\Qcoh^G(X)$. A $G$-equivariant
quasi-coherent $\Oc_X$-module is said to be {\it coherent}, {\it of
finite type} or {\it finitely presented} if it is coherent, of
finite type or finitely presented as an $\Oc_X$-module. We define
the categories $\Coh^G(X\on Y)$ and $\Coh^G(X)$ in the obvious way.

Let $\CG(X)_\qcb$ be the category of complexes of $G$-equivariant
quasi-coherent $\Oc_X$-modules, and let $\DG(X)_\qcb$ be the derived
category of $G$-equivariant quasi-coherent $\Oc_X$-modules. Note
that this notation may be confusing. We do not claim that
$\DG(X)_\qcb$ is the same as the derived category of $G$-equivariant
$\Oc_X$-modules with quasi-coherent cohomology. For a coherent
quasi-compact $G$-scheme $X$ we set
$$\K^G(X\on Y)=\[b\Coh^G(X\on Y)\]b,\quad\K^G(X)=\K^G(X\on X).$$
The representation ring of $G$ is defined by $\Rb^G=\K^G(pt)$. It
acts on the group $\K^G(X\on Y)$ by tensor product.

To define the standard derived functors for equivariant sheaves we
need more material. There are a number of foundational issues to be
addressed in translating the theory of derived functors of
quasi-coherent sheaves from the non equivariant setting to the
equivariant one. Here we briefly consider the issues that are
relevant to the present paper.

\subhead{1.5.3.~Definitions}\endsubhead $(a)$ An {\it ample family
of line bundles} on $X$ is a family of line bundles $\{\Lc_i\}$ such
that for every quasi-coherent $\Oc_X$-module $\Ec$ the evaluation
map yields an epimorphism
$$\bigoplus_i\bigoplus_{n> 0}\Gamma(X,\Ec\otimes_X\Lc_i^{\otimes n})
\otimes\Lc_i^{-\otimes n}\to\Ec.$$ We'll say that $X$ {\it
satisfies the property $(A_G)$} if it has an ample family of
$G$-equivariant line bundles.

\vskip1mm

$(b)$ We say that $X$ {\it satisfies the (resolution) property
$(R_G)$} if for every $G$-equivariant quasi-coherent $\Oc_X$-module
$\Ec$ there is a $G$-equivariant, flat, quasi-coherent
$\Oc_X$-module $\Pc$ and a surjection of $G$-equivariant
$\Oc_X$-modules $f:\Pc\to\Ec$. We'll also demand that we can choose
$\Pc$ and $f$ in a functorial way with respect to $\Ec$.

\vskip1mm

$(c)$ We say that a $G$-equivariant complex of quasi-coherent
$\Oc_X$-modules $\Ec$ admits a {\it $K$-flat resolution} if there is
a $G$-equivariant quasi-isomorphism $\Pc\to\Ec$ with $\Pc$ a
$G$-equivariant complex of quasi-coherent $\Oc_X$-modules such that
$\Pc\otimes_X\Fc$ is acyclic for every acyclic complex $\Fc$ in
$\CG(X)_\qcb$, see \cite{Sp, def.~5.1}.

\vskip1mm

$(d)$ We say that a $G$-equivariant complex of quasi-coherent
$\Oc_X$-modules $\Ec$ admits a {\it $K$-injective resolution} if
there is a $G$-equivariant quasi-isomorphism $\Ec\to\Ic$ with $\Ic$
a $G$-equivariant complex of quasi-coherent $\Oc_X$-modules such
that the complex of chain homomorphisms $\Fc\to\Ic$ in
$\CG(X)_\qcb$ is acyclic for every acyclic complex $\Fc$ in
$\CG(X)_\qcb$, see \cite{Sp, def.~1.1}.

\vskip3mm

If $G$ is the trivial group we'll abbreviate $(A)=(A_G)$ and $(R)=(R_G)$.

\vskip3mm

\subhead 1.5.4.~Remarks\endsubhead $(a)$ The property $(A_G)$
implies the property $(R_G)$. It implies also that any
$G$-equivariant quasi-coherent $\Oc_X$-module of finite type is the
quotient of a $G$-equivariant vector bundle, because $X$
is quasi-compact \cite{GD, chap.~0, (5.2.3)}.

\vskip1mm

$(b)$ If $G$ is linear and $X$ is Noetherian, normal, and satisfies
the property $(A)$, then $X$ satisfies also the property
$(A_G)$ \cite{T3, lem.~2.10 and sec.~2.2}. Since any
quasi-projective scheme satisfies $(A)$, we recover the
well-known fact that $X$ satisfies the property $(R_G)$ if it is
quasi-projective and normal and if $G$ is linear.

\vskip1mm

$(c)$ If $G$ is linear and $X$ is  Noetherian and regular, then $X$
satisfies $(A_G)$ by part $(b)$, because it satisfies $(A)$
\cite{SGA6, II.2.2.7.1}.

\vskip1mm

$(d)$ Let $X$ be an admissible $G$-scheme represented as
the inverse limit of a system of $G$-schemes $(X^\a)$
as in 1.4.5$(a)$. If $X^\a$ satisfies $(A_G)$ for
some $\a$ then $X$ satisfies also $(A_G)$, as well as $X^\b$ for
each $\b\geqslant\a$ \cite{TT, ex.~2.1.2$(g)$}.
Thus if $X$ is an admissible $G$-scheme which satisfies
the property $(S)$ in 1.2.19 then it satisfies also the
property $(A_G)$ (as well as $X^\a$ for each $\a$) by part $(c)$
above.

\vskip3mm

The $G$-equivariant quasi-coherent sheaves are well-behaved on
quasi-compact schemes satisfying the property $(A_G)$. {\it In the
rest of Section 1.5 all $G$-schemes are assumed to be quasi-compact
and to satisfy $(A_G)$}.

\proclaim{1.5.5.~Lemma} Assume that $X$ is coherent. Then any
$G$-equivariant quasi-coherent $\Oc_X$-module is the direct colimit
of a system of $G$-equivariant coherent $\Oc_X$-modules.
\endproclaim

\noindent{\sl Proof :} For any $G$-equivariant quasi-coherent
$\Oc_X$-module $\Ec$ the property $(A_G)$ yields an epimorphism in $\Qcoh^G(X)$
$$\Fc=\bigoplus_i\bigoplus_{n> 0}\Gamma(X,\Ec\otimes_X\Lc_i^{\otimes n})
\otimes\Lc_i^{-\otimes n}\to\Ec.$$  Any (rational) $G$-module is
locally finite, see e.g., \cite{J, sec.~I.2.13}. Choose a finite
number of $i$'s and $n$'s and a finite dimensional $G$-submodule of
$\Gamma(X,\Ec\otimes_X\Lc_i^{\otimes n})$ for each $i$ and each $n$
in these finite sets. Then $\Fc$ is represented as the union of a system of
$G$-equivariant locally free $\Oc_X$-submodules of finite type
$\Fc_\a\subset\Fc$. Taking the image under the epimorphism above we
can represent $\Ec$ as the direct colimit of a system of
$G$-equivariant $\Oc_X$-submodules of finite type $\Ec_\a$ with
surjective maps $\phi_\a:\Fc_\a\to\Ec_\a$. The kernel of $\phi_\a$
is again a $G$-equivariant quasi-coherent $\Oc_X$-module.
Considering its finitely generated $G$-equivariant quasi-coherent
$\Oc_X$-submodules, we prove as in \cite{GD, chap.~I, cor.~6.9.12}
that $\Ec$ is the direct colimit of a system of $G$-equivariant
finitely presented $\Oc_X$-modules. Since $X$ is a coherent scheme,
any finitely presented $\Oc_X$-module is coherent.

\qed

\vskip3mm

\proclaim{1.5.6.~Proposition} (a) Any complex in $\CG(X)_\qcb$
admits a $K$-flat resolution.

(b) There is a left derived tensor product
$\DG(X)_\qcb\times\DG(X)_\qcb\to\DG(X)_\qcb.$

(c) If $f:X\to Y$ is a morphism of $G$-schemes there is a left
derived functor $Lf^*:\DG(Y)_\qcb\to\DG(X)_\qcb$.

\endproclaim

\noindent{\sl Proof :} The non equivariant case is treated in
\cite{Sp}. The equivariant case is very similar and is left to the
reader. For instance, part $(a)$ is proved as in \cite{Li,
sec.~2.5}, while parts $(b)$, $(c)$ follow from $(a)$ and the
general theory of derived functors \cite{Li, sec.~2.5, ~3.1},
\cite{KS1}, \cite{KS2}.

\qed

\vskip3mm

\proclaim{1.5.7.~Proposition}  (a) The category $\Qcoh^G(X)$ is a
Grothendieck category. It has enough injective objects. Any complex
of $\CG(X)_\qcb$ has a $K$-injective resolution.

(b) If $f:X\to Y$ is a morphism of $G$-schemes there is a right
derived functor $Rf_*:\DG(X)_\qcb\to\DG(Y)_\qcb$.
\endproclaim

\noindent{\sl Proof :} Part $(b)$ follows from $(a)$ and the general
theory of derived functors \cite{Li}, \cite{KS1}, \cite{KS2}. Let us
concentrate on part $(a)$. The second claim is a well-known
consequence of the first one. The third claim follows also from the
first one by \cite{S}. See also \cite{AJS, thm.~5.4}, \cite{KS2,
sec.~14}. So we must check that $\Qcoh^G(X)$ is a Grothendieck
category. To do so we must prove that it has a generator, that it is
cocomplete, and that direct colimits are exact.
Fix a small category $\Acb$ and a functor $\Acb\to\Qcoh^G(X)$,
$\a\mapsto\Ec_\a$. Composing it with the forgetful functor we get a
functor $\Acb\to\Qcoh(X)$ with a colimit
$$\Ec=\ind_\a for(\Ec_\a),$$
because the category $\Qcoh(X)$ is cocomplete. For the same reason
we have also the following colimits
$$\ind_\a a^*(for(\Ec_\a)),\quad \ind_\a p^*(for(\Ec_\a)).$$
Since the functors $a^*$, $p^*$ have right adjoints,
a general result yields
$$a^*(\Ec)=\ind_\a a^*(for(\Ec_\a)),\quad
p^*(\Ec)=\ind_\a p^*(for(\Ec_\a)).$$ Next, since $(\Ec_\a)$ is a
system of $G$-equivariant quasi-coherent sheaves we have an
isomorphism of systems $a^*(for(\Ec_\a))\to p^*(for(\Ec_\a))$.
Taking the colimit we get an isomorphism of quasi-coherent sheaves
$a^*(\Ec)\to p^*(\Ec).$ This isomorphism yields a $G$-equivariant
structure on $\Ec$. The resulting $G$-equivariant sheaf is a colimit
in $\Qcoh^G(X)$. Thus the category $\Qcoh^G(X)$ is cocomplete and
the functor $for$ preserves colimits. Since $for$ reflects exactness
and $\Qcoh(X)$ is a Grothendieck category we obtain that the direct
colimit is an exact functor in $\Qcoh^G(X)$.
Finally we must prove that the Abelian category $\Qcoh^G(X)$ has a
generator.
This is obvious, because the proof of 1.5.5 implies that the tensor
powers of the $\Lc_i$'s generate the category $\Qcoh^G(X)$.

\qed

\vskip3mm

\subhead 1.5.8.~Compatibility of the derived functors\endsubhead It
seems to be unknown whether for any quasi-compact $G$-scheme
satisfying the property $(A_G)$ the derived tensor product, the
derived pull-back and the derived direct image satisfy the
equivariant analogue of the properties in 1.2.10-1.2.13. Here we
briefly discuss a weaker version of those which is enough for the
present paper.

First $Rf_*$ is right adjoint to $Lf^*$, next $Rf_*$ preserves the
cohomologically bounded complexes, and finally $Lf^*$ commutes with
the derived tensor product. These three properties are proved as in
the non equivariant case, see e.g., \cite{Li, prop.~3.2,~3.9},
\cite{KS2, sec.~14,18}. For instance the second one follows from the
spectral sequence $R^pf_*\circ H^q\Rightarrow R^{p+q}f_*$, where
$R^pf_*=H^p\circ Rf_*$, and the third one from the fact that both
derived functors can be computed via K-flat resolutions.

Next $Lf^*$ and the derived tensor product both commute with the
forgetful functor $for$ because they can be computed via K-flat
resolutions in both the equivariant and the non equivariant cases,
and because $for$ takes flat $\Oc_X$-modules in $\Qcoh^G(X)$ to flat
$\Oc_X$-modules in $\Qcoh(X)$.

The remaining properties require some work and more hypothesis.
We'll say that a quasi-coherent $\Oc_X$-module is {\it
$f_*$-acyclic} if it is annihilated by $R^pf_*$ for each $p>0$. By
the general theory of derived functors, for any $G$-equivariant
quasi-coherent sheaf $\Ec$ the complex $Rf_*(\Ec)$ can be computed
using a $G$-equivariant resolution of $\Ec$ by $f_*$-acyclic
sheaves. Assume that $G$ is a pro-linear group. The following lemma
is standard.

\proclaim{1.5.9.~Lemma} Let $f:X\to Y$ be a morphism of $G$-schemes.
If $X$ is normal and quasi-projective then any $G$-equivariant
quasi-coherent sheaf has a $G$-equivariant right resolution by
$f_*$-acyclic quasi-coherent sheaves.
\endproclaim

\noindent{\sl Proof :} By Sumihiro's theorem there is a
$G$-equivariant ample line bundle $\Lc$ on $X$, see e.g., \cite{CG,
sec.~5.1}. For a large enough integer $n>0$ the sheaf $\Lc^n$ is
generated by its global sections, and we have a $G$-equivariant
inclusion $\Oc_X\subset\Gc_n=\Lc^n\otimes V_n$, where $V_n$ is a
finite dimensional rational $G$-module such that the cokernel is a
locally free $\Oc_X$-module. For any $G$-equivariant coherent sheaf
$\Ec$ we have an inclusion $\Ec\subset\Ec\otimes_X\Gc_n$ such that
$\Ec\otimes_X\Gc_n$ is $f_*$-acyclic (because $\Lc$ is ample). If
$\Ec$ is a $G$-equivariant quasi-coherent sheaf we can represent it
as the direct limit $\Ec=\ind_\a\Ec_\a$ of a system of
$G$-equivariant coherent sheaves.  Choose integers $n_\a$ such that
$\Ec_\a\otimes_X\Gc_{n_\a}$ is $f_*$-acyclic for each $\a$ and
$\Gc_{n_\a}\subset\Gc_{n_\b}$ for $\a\leqslant\b$. We have an
inclusion of $G$-equivariant quasi-coherent sheaves
$\Ec\subset\Ec\otimes_X\Gc$ where $\Gc=\ind_\a\Gc_{n_\a}$ and
$\Ec\otimes_X\Gc$ is $f_*$-acyclic. The cokernel
$(\Ec\otimes_X\Gc)/\Ec$ is again a $G$-equivariant quasi-coherent
sheaf. By induction we get a resolution of $\Ec$.

\qed

\vskip2mm

We'll say that a coherent $G$-scheme $X$ is {\it
almost-quasi-projective} if it is represented as the inverse limit
of a system of normal quasi-projective $G$-schemes $X^\a$ with
affine morphisms such that the $G$-action on $X^\a$ factors through
$G^{n_\a}$ for some integer $n_\a$, compare 1.4.5$(a)$.

Let $f:X\to Y$ be a morphism of almost-quasi-projective $G$-schemes.
Then any $G$-equivariant quasi-coherent sheaf $\Ec$ over $X$ has a
$G$-equivariant right resolution by $f_*$-acyclic quasi-coherent
sheaves by 1.2.18. Therefore the general theory of derived functors
implies that $Rf_*(for(\Ec))=for(Rf_*(\Ec))$ in $\Dcb(Y)_\qcb$.
Further, for any $G$-equivariant quasi-coherent sheaves $\Ec$, $\Fc$
over $X$, $Y$ respectively the projection formula holds, i.e., we
have
$$Rf_*(\Ec)\Lotimes_Y\Fc=Rf_*(\Ec\Lotimes_XLf^*(\Fc)).$$ Indeed,
by adjunction we have a natural projection map from the lhs to the
rhs. To prove that this map is invertible it is enough to observe
that it is invertible in the non equivariant case, because the
forgetful functor commutes with $Rf_*$, $Lf^*$ and the derived
tensor product. The details are left to the reader. Note that, since
the projection formula and the base change hold for coherent
sheaves, they hold a fortiori in $G$-equivariant K-theory. This is
enough for this paper.

Recall that we have assumed that all $G$-schemes are quasi-compact
and satisfy the property $(A_G)$. This insures that the standard
derived functors are well-defined. {\it In the rest of Section 1.5
we'll also assume that the standard derived functors satisfy the properties
in 1.5.8}.

\subhead 1.5.10.~Equivariant coherent sheaves over an ind-coherent
ind-scheme\endsubhead Let $G$ be a group-scheme, $X$ be an
ind-coherent ind-$G$-scheme, and $Y\subset X$ be a closed
ind-subscheme which is preserved by the $G$-action. Let $(Y_\a)$,
$(X_\a)$ be systems of quasi-compact $G$-schemes as in 1.4.4,
such that the inclusion $Y\subset X$ is represented by a system of inclusions
$Y_\a\subset X_\a$. Since the maps $i_{\a\b}:X_\a\to X_\b$ are
good $G$-subschemes we have a direct system of Abelian categories
$\Coh^G(X_\a\on Y_\a)$ and exact functors $(i_{\a\b})_*$. We set
$$\gathered\Coh^G(X\on Y)=2\ind_\a\Coh^G(X_\a\on Y_\a),\cr
\K^G(X\on Y)=\[b\Coh^G(X\on Y)\]b=\ind_\a \K^G(X_\a\on
Y_\a).\endgathered\leqno(1.5.1)$$

\proclaim{1.5.11.~Proposition} The category $\Coh^G(X\on Y)$ is
independent of the choice of the system $(X_\a)$, up to canonical
equivalence. The group $\K^G(X\on Y)$ is independent of the choice
of the system $(X_\a)$, up to canonical isomorphism.
\endproclaim

\noindent{\sl Proof :} Let $(\tilde X_{\tilde\a})$ be another direct
system of closed subschemes of $X$ representing $X$. So we have
$X=\ind_\a X_\a$ and $X=\ind_{\tilde\a}\tilde X_{\tilde\a}.$
The second equality means that each $X_\a$ is included into some
$\tilde X_{\tilde\a}$ as a closed subset and vice-versa. Therefore
the 2-limits of both systems are identified. \qed

\vskip3mm

Once again we write $\Coh^G(X)=\Coh^G(X\on X)$ and
$\K^G(X)=\K^G(X\on X)$. Note that the tensor product $\otimes_X$ yields an
action of the ring $\Rb^G$ on the category $\Coh^G(X\on Y)$ and on the Abelian group $\K^G(X\on Y)$.

\vskip3mm

\subhead 1.5.12.~Admissible ind-coherent ind-schemes and reduction
of the group action\endsubhead Let $G$ be a pro-linear group. Fix a
system $(G^n)$ as in 1.4.2.
For each integer $n\geqslant 0$
let $G_{n}$ be the kernel
of the canonical map $G\to G^{n}$.
Let $X$ be an admissible coherent
$G$-scheme. Let $X^\a$, $n_\a$ be as in 1.4.5$(a)$.
We have a direct system of categories $\Coh^{G^{n}}(X^\a)$,
$n\geqslant n_\a$, with exact functors.
The pull-back by the canonical map $p_\a:X\to X^\a$
yields a functor
$$2\ind_\a 2\ind_{n\geqslant n_\a}\Coh^{G^{n}}(X^\a)\to\Coh^G(X).\leqno(1.5.2)$$

\proclaim{1.5.13.~Proposition} $(a)$
Assume that the pro-object $X=\pro_\a X^\a$ satisfies the property $(S)$.
Then the functor (1.5.2) is an equivalence of Abelian categories,
and it yields a group isomorphism
$$\ind_\a \ind_{n\geqslant n_\a}\K^{G^{n}}(X^\a)\to\K^G(X).$$
$(b)$ If $G_{0}$ is pro-unipotent, the canonical map
$\K^{G^0}(X)\to\K^G(X)$ is invertible.
\endproclaim

\noindent{\sl Proof :} The proof of $(b)$ is standard, see e.g.,
\cite{CG}. Let us concentrate on part $(a)$. The functor (1.5.2) is
fully faithful by 1.2.18$(a)$. Let us check that it is essentially
surjective. To do so, fix a $G$-equivariant coherent $\Oc_X$-module
$\Ec$. By 1.2.18$(a)$ there is an $\a$ and a coherent
$\Oc_{X^\a}$-module $\Ec^\a$ such that $\Ec=p_\a^*(\Ec^\a)$. We must
check that we can choose $\Ec^\a$ such that the $G$-action on $\Ec$
factors to a $G^n$-action on $\Ec^\a$ for some $n\geqslant n_\a$.
The unit of the adjoint pair of functors $(p_\a^*,(p_\a)_*)$ yields
an inclusion of quasi-coherent $\Oc_{X^\a}$-modules
$$\Ec^\a\subset(p_\a)_*p_\a^*(\Ec^\a)=(p_\a)_*(\Ec).$$
Since $X^\a$ is a Noetherian $G$-scheme and $(p_\a)_*(\Ec)$ is a
quasi-coherent $G$-equivariant $\Oc_{X^\a}$-module, we know that
$(p_\a)_*(\Ec)$ is the union of all its $G$-equivariant coherent
subsheaves. Fix a $G$-equivariant coherent $\Oc_{X^\a}$-module
$\Fc^\a$ containing $\Ec^\a$. The $G$-action on $\Fc^\a$ factors to
an action of the linear group $G^n$ for some $n\geqslant n_\a$. Let
$\Gc^\a\subset\Fc^\a$ be the $G^n$-equivariant quasi-coherent
subsheaf generated by $\Ec^\a$. It is again a coherent
$\Oc_{X^\a}$-module, because $\Fc^\a$ is coherent and $X^\a$ is
Noetherian. Since $\Ec$ is already $G$-equivariant the inclusion
$$\Ec=p_\a^*(\Ec^\a)\subset p_\a^*(\Gc^\a)$$ is indeed an
equality of $\Oc_X$-modules $\Ec\simeq p_\a^*(\Gc^\a)$.

\qed

\vskip3mm

Now, let $X$ be an admissible ind-coherent ind-$G$-scheme
represented as the direct colimit of a system of admissible
$G$-schemes $(X_\a)$ as in 1.4.5$(c)$. By (1.5.1) we have
$$
\Coh^G(X)=2\ind_\a \Coh^{G}(X_\a),\quad \K^G(X)=\ind_\a\K^{G}(X_\a).
$$
If $G_{0}$ is pro-unipotent then 1.5.13 yields isomorphisms
$$\K^{G^{0}}(X)=\K^{G}(X),\quad\Rb^{G^{0}}=\Rb^G.$$
This is called the {\it reduction of the group action}.

\subhead 1.5.14.~Thom isomorphism and pro-finite-dimensional
vector bundles over ind-schemes\endsubhead
A {\it vector bundle} over the ind-scheme $Y$
is an ind-scheme homomorphism $X\to Y$ which is represented as the
direct colimit of a system of vector bundles $X_\a\to Y_\a$. More precisely,
we require that
$X=\ind_\a X_\a$, $Y=\ind_\a Y_\a$,
and for $\a\leqslant \b$ we have a Cartesian square
$$\xymatrix{ X_\a\ar[d]\ar[r]&X_\b\ar[d]\\ Y_\a\ar[r]  &Y_\b}$$
such that the vertical maps are vector-bundles
and the upper horizontal map is a morphism of vector bundles.
To any vector bundle over an ind-scheme $X$ we can associate its
sheaf of sections which is a quasi-coherent $\Oc_X$-module, see 1.3.5.

A {\it pro-finite-dimensional vector bundle} over the ind-scheme $Y$
is defined in the same way by replacing everywhere vector bundles by
pro-finite-dimensional vector bundles, see 1.2.23.
In other words, it is an ind-scheme homomorphism which is represented
as the ``double limit" of a system of
vector-bundles $X^{n}_\a\to Y_\a$ with $n$ an integher $\geqslant 0$.
Further, for each $\a$ and each $m\geqslant n$
we have a vector-bundle homomorphism
$X^{m}_\a\to X^{n}_\a$ over $Y_\a$, and for each $n$ and each $\a\leqslant \b$
we have an isomorphism of vector-bundles $X^{n}_\a\to
X^{n}_\b\times_{Y_\b}Y_\a$. We require that these data satisfy the
obvious composition rules. In particular for each $m\geqslant n$ and
each $\b\geqslant\a$ the following square is Cartesian
$$\xymatrix
{X^{m}_\a\ar[d]\ar[r]&X^{m}_\b\ar[d]\cr X^{n}_\a\ar[r]
&X^{n}_\b.}
$$

Note that a pro-finite-dimensional vector bundle over an ind-coherent
ind-scheme is again an ind-coherent ind-scheme.

Let $\pi:X\to Y$ be an admissible $G$-equivariant
pro-finite-dimensional vector bundle over an admissible ind-coherent
ind-$G$-scheme $Y$. From 1.5.10 and base change we get an exact
functor $\pi^*:\Coh^G(Y)\to\Coh^G(X)$. It factors to a group
homomorphism $\pi^*:\K^G(Y)\to\K^G(X).$ The {\it Thom isomorphism}
implies that $\pi^*$ is invertible.

\subhead 1.5.15.~Descent and torsors over ind-schemes\endsubhead
Fix a pro-linear group $G=\pro_n G^n$.
For each integer $n\geqslant 0$
let $G_{n}$ be the kernel
of the canonical map $G\to G^{n}$.

Let $Y$ be a scheme.
A {\it $G$-torsor} over $Y$
is a scheme homomorphism
$P\to Y$
which is represented as the inverse limit of a system
consisting of a $G^n$-torsor $P^n\to Y$ for each integer
$n\geqslant 0$
such that the morphism of $Y$-schemes $P^m\to P^n$, $m\geqslant n$,
intertwines the $G^{m}$-action on the
lhs and the $G^{n}$-action on the rhs, via the canonical
group-scheme homomorphism $G^{m}\to G^{n}$.

Now, assume that $Y=\ind_\a Y_\a$ is an ind-scheme.
A {\it $G$-torsor} over $Y$
is an ind-scheme homomorphism $P\to Y$ which is represented as the
direct colimit of a system of $G$-torsors $P_\a\to Y_\a$. More precisely,
we require that
$P=\ind_\a P_\a$,
that for each $\a$ we have
a system of $G^n$-torsors $P^n_\a\to Y_\a$ representing the
$G$-torsor $P_\a\to Y_\a$, and that
for each $n$ and each $\b\geqslant \a$ we have an isomorphism of $G^{n}$-torsors
$P^{n}_\a\to P^{n}_\b\times_{Y_\b}Y_\a.$
In particular,
for each $m\geqslant n$ and each $\b\geqslant \a$ we have a Cartesian square
$$\xymatrix
{P^{m}_\a\ar[r]\ar[d]&P^{m}_\b\ar[d]\cr P^{n}_\a\ar[r]
&P^{n}_\b.}
$$

Note that a $G$-torsor over a pro-smooth scheme is again a pro-smooth scheme,
and that a $G$-torsor over an ind-coherent ind-scheme is again an
ind-coherent ind-scheme.

Now, assume that $Y$ is a scheme and let $P\to Y$ be  a $G$-torsor.
Note that for each integer $n\geqslant 0$ we have a $G^n$-torsor $P/G_n\to Y$.
In the rest of this subsection we consider the {\it induction functors}.

Let $X$ be an admissible $G$-scheme and let $X^\a$, $n_\a$ be as in 1.4.5$(a)$.
For each $\a$ and each integer $n\geqslant n_\a$
the quotient space $(X^\a)_Y=P\times_GX^\a$ is equal to
the $Y$-scheme $(P/G_n)\times_{G^n}X^\a$.
Further, if $\b\geqslant\a$ the canonical map $X^\b\to X^\a$ yields a
$Y$-scheme homomorphism $(X^\b)_Y\to (X^\a)_Y$.
Thus the quotient space $X_Y=P\times_G X$ is a $Y$-scheme
which is represented as the inverse limit
$X_Y=\lim_\a(X^\a)_Y$.
By (1.5.2) we have an equivalence of categories
$$2\ind_\a 2\ind_{n\geqslant n_\a}\Coh^{G^{n}}(X^\a)\to\Coh^G(X).$$
By faithfully flat descent we have a functor
$$\Coh^{G^{n}}(X^\a)\to
\Coh((P/G_n)\times_{G^n}X^\a)=
\Coh((X^\a)_Y).$$
This yields a functor
(called {\it induction} functor)
$$\Coh^{G}(X)\to \Coh(X_Y).$$

Next, let $Z$ be an admissible ind-$G$-scheme which is represented as the direct
colimit of a system of admissible $G$-schemes $Z_\a$ as in 1.4.5$(c)$.
Then the quotient space $Z_Y=P\times_G Z$ is an ind-scheme over $Y$
which is represented as the direct colimit
$Z_Y=\lim_\a(Z_\a)_Y$, and $(Z_\a)_Y$ is defined as above for each $\a$.
If $Z$ is in-coherent and $Y$ is coherent then the ind-scheme $Z_Y$ is again
ind-coherent and (1.5.1) yields a functor
(called {\it induction} functor)
$$\Coh^{G}(Z)\to \Coh(Z_Y).$$

The induction functor is defined in a similar way if $Y$ is an ind-scheme.
The details of the construction are left to the reader.

\subhead 1.5.16.~Remark\endsubhead We define in a similar way
induction functors for quasi-coherent sheaves. Next, let $H$ be a
group-scheme acting on the $G$-torsor $P\to Y$, i.e., the group $H$
acts on $P$, $Y$ and the action commutes with the $G$-action and
with the projection $P\to Y$. Then there is a $H$-action on $Z_Y$
and the induction yields a functor $\Coh^{G}(Z)\to \Coh^H(Z_Y).$

\subhead 1.5.17.~Complements on the concentration map\endsubhead
Now, assume that $G=S$ is a diagonalizable linear group. Let $\Xb^S$
be the group of characters of $S$. Each $\l\in \Xb^S$ defines a
one-dimensional representation $\theta_\l$ of $S$. Let $\theta_\l$
denote also its class in $\Rb^S$. The ring $\Rb^S$ is spanned by the
elements $\theta_\l$ with $\l\in\Xb^S$. So any element of $\Rb^S$
may be viewed as a function on $S$. For any $\Sigma\subset S$ let
$\Rb^S_\Sigma$ be the ring of quotients of $\Rb^S$ with respect to
the multiplicative set of the functions in $\Rb^S$ which do not
vanish identically on $\Sigma$.

Now, let $X$ be an ind-coherent ind-$S$-scheme. We'll say that
$\Sigma$ is {\it $X$-regular} if the fixed points subsets
$X^\Sigma$, $X^S$ are equal. Write
$$\K^S(X)_\Sigma=\Rb^S_\Sigma\otimes_{\Rb^S}\K^S(X).$$
The {\it Thomason concentration theorem} says that the map
$$\K^S(X^S)_\Sigma\to\K^S(X)_\Sigma$$
given by the direct image by the canonical inclusion $X^S\subset X$
is invertible if $X$ is a scheme of finite type and $\Sigma$ is
$X$-regular \cite{T1}, \cite{T2}. We'll use some form of the
concentration theorem in some more general situation, which we
consider below. In each case, the proof of the concentration theorem
can be reduced to the original statement of Thomason using the
discussion above. It is left to the reader.

\subhead 1.5.18\endsubhead Let $X$ be a pro-smooth admissible $S$-scheme.
It is easy to check that the fixed-points subset $X^S\subset X$ is a
closed subscheme which is again pro-smooth. Thus the obvious
inclusion $j:X^S\to X$ has a finite tor-dimension by 1.2.10, 1.2.20.
Hence it yields a $\Rb^S$-linear map $Lj^*:\K^S(X)\to\K^S(X^S)$.
This map can be viewed as follows. Any coherent $\Oc_X$-module $\Ec$
has locally a finite resolution by locally free modules of finite
rank. Hence the $p$-th left derived functor
$L_pj^*\Ec=H^{-p}(Lj^*\Ec)$ vanishes for $p\gg 0$. We have
$Lj^*(\Ec)=\sum_{p\geqslant 0}(-1)^pL_pj^*(\Ec)$.
If $\Sigma$ is $X$-regular we get a group isomorphism
$$Lj^*:\K^S(X)_\Sigma\to\K^S(X^S)_\Sigma.$$

\subhead 1.5.19\endsubhead Let $X$ be an admissible ind-$S$-scheme of
ind-finite type. The inclusion of the fixed points subset $i:X^S\to
X$ is a good ind-subscheme. Thus the direct image yields a
$\Rb^S$-linear map $i_*:\K^S(X^S)\to\K^S(X).$
If $\Sigma$ is $X$-regular we get a group isomorphism
$$i_*:\K^S(X^S)_\Sigma\to\K^S(X)_\Sigma.$$

\subhead 1.5.20\endsubhead Let $X$ be a pro-smooth admissible $S$-scheme and
$f:Y\to X$ be an admissible ind-$S$-scheme over $X$. We'll assume
that the map $f$ is locally trivial in the following sense :
there is an admissible ind-$S$-scheme $F$ of ind-finite type
and a $S$-equivariant finite affine open cover $X=\bigcup_wX^w$ such
that over each $X^w$ the map $f$ is isomorphic to the obvious
projection $X^w\times F\to X^w$, where the group $S$ acts diagonally
on the lhs. The ind-scheme $Y$ is ind-coherent by 1.2.20.
Further, the fixed points subset $X^S$ is
again pro-smooth. Setting $Y'=X^S\times_XY$ we get the following diagram
$$\xymatrix{
Y^S\ar[r]^{i}\ar[d]&Y'\ar[r]^j\ar[d]&Y\ar[d]^f\cr
X^S\ar[r]^{=}&X^S\ar[r]&X.}
$$
Over the open set $X^w$ the map $j$ is isomorphic to the obvious inclusion
$$(X^w)^S\times F\subset X^w\times F.$$
The inclusion $(X^w)^S\subset X^w$
has a finite tor-dimension by 1.2.10, 1.2.20.
Thus the map $j$ has also a finite tor-dimension.
By base change we have a $\Rb^S$-linear map
$Lj^*:\K^S(Y)\to\K^S(Y').$ Since $i$ is the inclusion of a good ind-subscheme the direct
image gives a map $i_*:\K^S(Y^S)\to\K^S(Y')$. If $\Sigma$
is $Y$-regular we
get a group isomorphism
$$\gathered
(i_*)^{-1}\circ Lj^*:\K^S(Y)_\Sigma\to \K^S(Y^S)_\Sigma.
\endgathered
$$

\vskip3mm

\head 2.~Affine flag manifolds\endhead

\subhead 2.1.~Notation relative to the loop group\endsubhead

\subhead 2.1.1\endsubhead Let $G$ be a simple, connected and simply
connected linear group over $\CC$ with the Lie algebra $\gen$. Let $T\subset G$
be a Cartan subgroup and $W$ be the Weyl group of the pair $(G,T)$.
Recall that $\Xb^T$ is the Abelian group of characters of $T$ and that
$\Yb^T$ is the Abelian group of cocharacters of $T$.
Let $\ten$ be the Lie algebra of $T$ and $\ten^*$ be the set of linear forms
on $\ten$. We'll view $\Xb^T$, $\Yb^T$ as lattices in $\ten^*$, $\ten$
in the usual way. Note that, since $G$
is simply connected, the lattice $\Xb^T$ is spanned by the
fundamental weight and the lattice $\Yb^T$ is spanned by the simple
coroots. Let $\Xb_+^T\subset\Xb^T$ and $\Yb_+^T\subset\Yb^T$ denote
the monoids of dominant characters and cocharacters.

Fix a Borel subgroup $B\subset G$. Let $\Delta$ be the set of roots
of $(G,T)$ and $\Pi\subset\Delta$ be the subset of simple roots
associated to $B$. Let $\Delta^\vee$ be the set of coroots. Let
$\theta$ be the highest root and $\check\theta$ be the corresponding
coroot. Let
$$\la\ ,\ \ra:\Xb^T\times\Yb^T\to\ZZ,
\quad(\ ,\ ):\ten^*\times\ten^*\to\CC$$
be the canonical perfect
pairing and
the nondegenerate $W$-invariant bilinear form
normalized by $(\theta,\theta)=2$.
We'll denote by $\kappa$ the corresponding homomorphism
$\ten\to\ten^*$ and we'll abbreviate
$(\check\l,\check\mu)=(\kappa(\check\l),\kappa(\check\mu))$
for each $\check\l,\check\mu\in\ten$.

Let $\tilde\Delta$ be the set of affine roots, $\tilde\Delta_e$ be
the subset of positive affine roots and $\tilde\Pi$ be the subset of
simple affine roots. Let $\a_0\in\tilde\Pi$ be the unique simple
affine root which does not belong to $\Pi$. We have
$\tilde\Pi=\{\a_0,\a_1,\dots,\a_n\}$ where $n$ is the rank of $G$.
Let $\tilde\Delta^\vee$ be the set of affine coroots.
We have $\tilde\Pi^\vee=\{\check
\a_0,\check\a_1,\dots,\check\a_n\}$ where
$\check\a_i$ is the affine coroot associated to the simple affine root
$\a_i$ for each $i$.

Let $\tilde W=W\ltimes\Yb^T$ be the affine Weyl group of $G$.
For any affine real root $\a$ let $s_\a\in\tilde W$ be the
corresponding affine reflection. We'll abbreviate $s_i=s_{\a_i}$ for
each $i$. Since $G$ is simply connected the group $\tilde W$ is a
Coxeter group with simple reflections the $s_i$'s.

We'll
abbreviate $w=(w,0)$ and $\xi_{\check\lambda}=(e,\check\lambda)$ for
each $(w,\check\lambda)\in\tilde W$.
In particular we'll regard to $W$ as a subgroup of
$\tilde W$ in the obvious way.
Here $e$ denotes the unit, both in $W$ and in $\tilde W$.

\subhead 2.1.2\endsubhead We'll fix a decreasing sequence of subsets
$\tilde\Delta_l\subset\tilde\Delta_e$, with $l\in\NN$, such
that
$$(\tilde\Delta_l+\tilde\Delta_e)\cap\tilde\Delta_e\subset
\tilde\Delta_l,\quad
\sharp(\tilde\Delta_e\setminus\tilde\Delta_l)<\infty,\quad
\bigcap_l\Delta_l=\emptyset.$$ For instance we may set
$\tilde\Delta_l=l\delta+\tilde\Delta_e$ where $\delta$ is the
smallest positive imaginary root. Put also
$\tilde\Delta_{l}^\op=-\tilde\Delta_l$.

\subhead 2.1.3\endsubhead We'll abbreviate
$G((\varpi))=G(\CC((\varpi)))$,
$\gen((\varpi))=\gen(\CC((\varpi)))$, etc.
Recall that $K=G(\CC[[\varpi]])$. Let $I\subset K$ be the
standard Iwahori subgroup and
$I^\op\subset K^\op=G(\CC[[\varpi^{-1}]])$
be the opposite Iwahori subgroup. Let $N$, $N^\op$ be the
pro-unipotent radicals of $I$, $I^\op$ respectively. The groups $I$,
$I^\op$, $N$, $N^\op$ are compact.

\subhead 2.1.4\endsubhead Let $\nen$, $\nen^\op$, $\ien$, $\ken$ be
the Lie algebras of $N$, $N^\op$, $I$, $K$. For any integer $l\geqslant 0$ let
$\nen_l\subset\nen$ and
$\nen_l^\op\subset\nen^\op$
be the product of all weight subspaces
associated to the roots in $\tilde\Delta_l$,
$\tilde\Delta_l^\op$ respectively. Put also
$$\nen^l=\nen/\nen_l,\quad
\nen^{\op,l}=\nen^\op/\nen_l^\op,\quad
\nen_w=w(\nen)\cap\nen,\quad\nen^\op_w=w(\nen^\op)\cap\nen,\quad
w\in\tilde W.$$ Let $N_l$, $N_l^\op$, etc, be the groups associated
with the Lie algebras $\nen_l$, $\nen_l^\op$, etc. We'll write
$\tilde\Delta_w$, $\tilde\Delta^\op_w$ for the set of roots of
$\nen_w$, $\nen_w^\op$. Note that $\nen$, $\nen_l$, $\nen_w$ have a
natural structure of admissible $I$-equivariant compact coherent
schemes and that $\nen_l$, $\nen_w$ are good subschemes of $\nen$.
Further the quotient $\nen^l$ is finite dimensional and $\nen^\op_w$
has a natural structure of $I$-scheme of finite type. The $I$-action
is the adjoint one. We'll call an Iwahori Lie subalgebra of
$\gen((\varpi))$ any Lie subalgebra which is $G((\varpi))$-conjugate
to $\ien$.

\subhead 2.1.5\endsubhead The group $\CC^\times$ acts on
$\CC((\varpi))$ by loop rotations, i.e., a complex number $z\in\CC^\times$
takes a formal series $f(\varpi)$ to $f(z\varpi)$. This yields
$\CC^\times$-actions on $G((\varpi))$, $I$ and $\gen((\varpi))$.
Consider the semi-direct products
$$\hat G=\CC^\times\ltimes G((\varpi)),\quad
\hat I=\CC^\times\ltimes I,\quad\hat I^\op=\CC^\times\ltimes
I^\op,\quad \hat T=\CC^\times\times T.$$ The group $\hat G$ acts
again on $\gen((\varpi))$.

\subhead 2.1.6\endsubhead Let $\tilde G$ be the maximal, ``simply
connected", Kac-Moody group associated to $G$ defined by Garland
\cite{G}. It is a group ind-scheme which is a central extension
$$1\to\CC^\times\to \tilde G\to\hat G\to 1.$$
See \cite{K, sec.~13.2} for details. Let $\tilde I$, $\tilde T$,
$\tilde K$ be the corresponding Iwahori, Cartan and maximal compact
subgroup. Note that
$\tilde K=\hat K\times\CC^\times$,
$\tilde I=\hat I\times\CC^\times$ and $\tilde
T=\hat T\times\CC^\times$, i.e., the central extension splits.
We define also the opposite Iwahori group
$\tilde I^\op=\hat I^\op\times\CC^\times$. Let $\tilde\gen$,
$\tilde\ien$, $\tilde\ken$ be the Lie algebras of $\tilde G$,
$\tilde I$, $\tilde K$. The group $\tilde G$ acts on
$\gen((\varpi))$ and $\tilde\gen$ by the adjoint action. By an Iwahori
Lie subalgebra of $\tilde\gen$ we simply mean a Lie subalgebra which
is $\tilde G$-conjugate to $\tilde\ien$.

\subhead 2.1.7\endsubhead We'll also use the groups
$$\Gt=\tilde G\times\CC^\times,\quad \It=\tilde I\times\CC^\times,\quad
\Tt=\tilde T\times\CC^\times.$$ The group $\Gt$ acts also on
$\tilde\gen$. We simply require that an element $z\in\CC^\times$
acts by multiplication by $z$ (=by dilatations). Note that
$\Tt=T\times(\CC^\times)^3$. To distinguish the different copies of
$\CC^\times$ we may use the following notation : $\CC^\times_\rot$
corresponds to loop rotation, $\CC^\times_\cent$ to the central
extension, and $\CC_\qua^\times$ to dilatations. Thus we have
$$\hat T=T\times\CC^\times_\rot,\quad\tilde T=\hat
T\times\CC^\times_\cent,\quad \Tt=\tilde T\times\CC^\times_\qua.$$
We'll also write $T_\cent=T\times\CC^\times_\cent$.

\subhead 2.1.8\endsubhead  We'll abbreviate $\tilde\Xb=\Xb^{\tilde
T}$, $\Xb=\Xb^{\Tt}$, $\tilde\Yb=\Yb^{\tilde T}$ and
$\Yb=\Yb^{\Tt}$.
The pairing $\la\ ,\ \ra$ extends to the canonical pairing
$$\la\ ,\ \ra:\tilde\Xb\times\tilde\Yb\to\ZZ.$$
Let $d$, $c$ be the canonical generators of $\Yb^{\CC^\times_\rot}$,
$\Yb^{\CC^\times_\cent}$. We have
$$\tilde\Yb=\Yb^T\oplus\ZZ d\oplus\ZZ c=\ZZ d\oplus\bigoplus_{i=0}^n\ZZ\check\a_i,\quad
\check\a_0=c-\check\theta.$$ The affine fundamental weights are the
unique elements $\o_i\in\tilde\Xb$, $i=0,1,\dots, n$, such that
$\la\o_i,\check\a_j\ra=\delta_{i,j}$ for each $i,j$. We have
$$\tilde\Xb=\Xb^T\oplus\ZZ\delta\oplus\ZZ{\omega_0}=
\ZZ\delta\oplus\bigoplus_{i=0}^n\ZZ\o_i,\quad \Xb=\tilde\Xb\oplus\ZZ
t,$$ where $\delta$ is the smallest positive imaginary root. Recall
that $\a_0=\delta-\theta$.
Then $\delta$, $\omega_0$, $t$ are the canonical generators of
$\Xb^{\CC^\times_\rot}$, $\Xb^{\CC^\times_\cent}$, and
$\Xb^{\CC^\times_\qua}$ respectively.

\subhead 2.1.9\endsubhead  There is a $\tilde W$-action on
$\tilde\Xb$, $\tilde\Yb$ such that the natural pairing is $\tilde
W$-invariant. It is given by :

\vskip1mm

\item{$\bullet$} $W$ fixes the elements $\delta$, $\omega$, $d$, $c$
and it acts in the usual way on $\Xb^T$, $\Yb^T$,

\vskip1mm

\item{$\bullet$} the element $\xi_{\check\l}$, $\check\l\in\Yb^T$,
maps $\mu\in\tilde\Xb$ to
$$\mu+\la\mu,c\ra\,\kappa(\check\l)-
\bigl(\la\mu,\check\l\ra+(\check\l,\check\l)\la\mu,c\ra/2\bigr)\delta,$$

\vskip1mm

\item{$\bullet$} the element $\xi_{\check\l}$, $\check\l\in\Yb^T$,
maps $\check\mu\in\tilde\Yb$ to
$$\check\mu+\la\delta,\check\mu\ra\check\l-
\bigl(\la\kappa(\check\l),\check\mu\ra+(\check\l,\check\l)\la\delta,\check\mu\ra/2\bigr)c.$$

\noindent This action is denoted by
${}^w\mu$, ${}^w\check\mu$ for each $w\in\tilde W$, $\mu\in\tilde\Xb$, and
$\check\mu\in\tilde\Yb$.

\subhead 2.1.10\endsubhead  \noindent There is also a $\tilde
W$-action on $\tilde T$. It is given by :

\vskip1mm

\item{$\bullet$} $W$ fixes $\CC^\times_\rot$, $\CC^\times_\cent$ and it acts on the usual way on $T$,

\vskip1mm

\item{$\bullet$} the element $\xi_{\check\l}$, $\check\l\in\Yb^T$,
maps the pair $(s,\tau)$ with $s\in T_\cent$ and
$\tau\in\CC^\times_\rot$ to the pair
$$\bigl(s\check\l(\tau)c(\kappa(\check\l)(sh))^{-1},\tau\bigr)\ \roman{with}\
h^2=\check\l(\tau).$$ Here we regard $\check\l,c$ as group
homomorphisms  $\CC^\times\to T_\cent$ and $\kappa(\check\l)$ as a
group homomorphism $T_\cent\to\CC^\times$.

\subhead 2.1.11\endsubhead  Since $I$ is a group-scheme the ring
$\Rb^I$ is well-defined. By devissage we have $\Rb^{I}=\Rb^T$.
Recall that $\Rb^T=\sum_{\l\in\Xb^T}\ZZ\theta_\l$ is the group
algebra of $\Xb^T$, see 1.5.17. We'll abbreviate $q=\theta_\delta$,
$t=\theta_t$ in $\Rb^\Tt$. We may also use the following
$\ZZ_t$-algebras
$$\ZZ_t\tilde\Xb=\sum_{\l\in\tilde\Xb}\ZZ_t\theta_\l,\quad
\ZZ_t\tilde\Yb=\sum_{\check\l\in\tilde\Yb}\ZZ_t\theta_{\check\l}.$$
Note that $\ZZ_t\tilde\Xb=\Rb^\Tt$.

\vskip3mm

\subhead 2.2.~Reminder on the affine flag manifold\endsubhead

\subhead 2.2.1.~The affine flag manifold\endsubhead
Let $\Fen=\Fen_G$ be the affine flag
manifold of $G$. It is an ind-proper ind-scheme of ind-finite type
whose set of $\CC$-points is
$$\Fen=G((\varpi))/I=\{\roman{Iwahori\ Lie\ subalgebra\ of}\ \gen((\varpi))\}.$$
The space $\Fen$ can be viewed as the sheaf
for the fppf topology over the flat affine site over $\CC$,
associated with the quotient pre-sheaf $\tilde G/\tilde I$.
In particular there is a canonical ind-scheme homomorphism
$\tilde G\to\Fen$ which is a $\tilde I$-torsor as in 1.5.15.
The set
of $\CC$-points of $\Fen$ is simply the quotient set $\tilde G/\tilde I$.
It will be convenient to regard an element of this set as an
Iwahori Lie subalgebra of $\tilde\gen$.
The ind-group $\tilde G$ acts on itself by left multiplication.
This action yields a $\tilde G$-action on $\Fen$.
The group-scheme $\tilde I$ acts also on $\Fen$, and the later
has the structure
of an admissible ind-$\tilde I$-scheme.
The $\tilde I$-orbits are numbered by the elements of $\tilde W$
$$\Fen=\bigsqcup_{w\in\tilde W}{\buildrel\circ\over\Fen}_w.$$
Let $\leqslant$ be the Bruhat order on $\tilde W$.
We have
$$\Fen=\ind_w\Fen_{w},\quad
\Fen_{w}=\bigsqcup_{v\leqslant w}{\buildrel\circ\over\Fen}_v.$$
Further $\Fen_w$ is a projective, normal, $\tilde I$-scheme for every $w$.
We have
$$\K(\Fen)=\ind_w\K(\Fen_{w}),\quad\K(\Fen_w)=\[b\Coh(\Fen_{w})\]b.$$
For a future use, we'll abbreviate $\Den=\Fen\times\Fen$. For each
$v$, $w$  let $\Den_{v,w}=\Fen_v\times\Fen_w$.

\subhead 2.2.2.~The Kashiwara affine flag manifold\endsubhead We'll
also use  the Kashiwara flag manifold $\Xen$. See \cite{K},
\cite{KT1}, \cite{KT2} for details. It is a coherent, pro-smooth
(non quasi-compact) scheme
locally of countable type with a left $\tilde I^\op$-action.
Recall that a
$G$-scheme $X$ is {\it locally free} if any point of $X$ has a
$G$-stable open neighborhood which is isomorphic, as a $G$-scheme,
to $G\times Y$ for
some scheme $Y$. In this case the quotient $X/G$ is representable by
a scheme.
The Kashiwara flag manifold
is constructed as a quotient $\Xen=\tilde G_\infty/\tilde I$, where
$\tilde G_\infty$ is a coherent scheme with a locally free left
action of $\tilde I^\op$ and a locally free right action of $\tilde I$.
In particular there is a canonical scheme homomorphism
$\tilde G_\infty\to\Xen$ which is a $\tilde I$-torsor as in 1.5.15.
There is a $\tilde I^\op$-orbit decomposition
$$\Xen=\bigsqcup_{w\in \tilde W}{\buildrel\circ\over\Xen}{}^w$$
where ${\buildrel\circ\over\Xen}{}^w$ is a locally closed
subscheme of codimension $l(w)$ (=the length of $w$ in $\tilde W$)
which is isomorphic to the
infinite-dimensional affine space
$\AA^\NN.$ The Zariski closure of ${\buildrel\circ\over\Xen}{}^w$ is
$\bigsqcup_{v\geqslant w}{\buildrel\circ\over\Xen}{}^v.$
The scheme $\Xen$ is covered by the following open subsets
$$\Xen^{w}=\bigsqcup_{v\leqslant w}{\buildrel\circ\over\Xen}{}^v.$$
Note that $\Xen^w$ is a $\tilde I^\op$-stable finite union of
translations of the big cell $\Xen^e$ and that
$\Xen^e\simeq\AA^\NN$. Thus $\Xen^w$ is quasi-compact and pro-smooth.
Since $\Xen$ is
not quasi-compact, we have
$$\K(\Xen)=\pro_w\K(\Xen^{w}),\quad\K(\Xen^w)=\[b\Coh(\Xen^{w})\]b.$$
For each $w$ there is a closed immersion $\Fen_{w}\subset\Xen^{w}$,
see \cite{KT2, prop.~1.3.2}.
Therefore the restriction of $\Oc$-modules yields a functor
$$\Qcoh(\Xen)\to\Qcoh(\Fen).$$ The tensor product of
quasi-coherent $\Oc_\Xen$-modules yields a functor
$$\otimes_\Xen:\ \Coh(\Fen)\times\Coh(\Xen)\to\Coh(\Fen).
$$ Since $\Xen^{w}$ is pro-smooth we have also a group homomorphism
$$\otimes^L_\Xen:\ \K(\Fen)\otimes\K(\Xen)\to\K(\Fen).$$
Finally, we have the following important property.

\proclaim{2.2.3.~Proposition} The $\tilde I^\op$-scheme $\Xen^w$ is
admissible and it satisfies $(A_{\tilde I^\op})$.
\endproclaim

\noindent{\sl Proof :} The admissibility follows from 1.4.6. Given
an integer $l\geqslant 0$ we consider the quotients
$$\Xen^{w,l}=N^\op_l\setminus\Xen^w,\quad
\tilde I^{\op,l}=\tilde I^\op/N^\op_l.$$ Note that $\tilde
I^{\op,l}$ is a linear group, that $\Xen^{w,l}$ is a smooth $\tilde
I^{\op,l}$-scheme and that the canonical map $\Xen^w\to\Xen^{w,l}$
is a $\tilde I^\op$-equivariant $N^\op_l$-torsor \cite{KT2,
lem.~2.2.1}. A priori $\Xen^{w,l}$ could be not separated. See the
remark after \cite{KT2, lem.~2.2.1}. The separatedness is proved in
\cite{VV, sec.~A.6}. See also 2.2.4 below. Since the $\tilde
I^{\op,l}$-scheme $\Xen^{w,l}$ is Noetherian and regular it
satisfies the property $(A_{\tilde I^{\op,l}})$. Then $\Xen^w$
satisfies also the property $(A_{\tilde I^\op})$ by 1.5.4.

\qed

\vskip3mm

\subhead 2.2.4.~Remarks\endsubhead
$(a)$ The scheme $\Xen^{w,l}$ above is
separated if $l$ is large enough, even in the more general case of
Kac-Moody groups considered in \cite{KT2}. This follows from
\cite{TT, prop.~C.7} and the fact that $\Xen^w$ is a separated
scheme.

$(b)$ The $\tilde T$-fixed points subsets
in ${\buildrel\circ\over\Fen}_w$ and
${\buildrel\circ\over\Xen}{}^w$ are reduced to the same
single point. We'll denote it by $\ben_w$. Note that $\ben_e$ is
identified with the Iwahori Lie algebra $\ien$ (or $\tilde\ien$)
for $e$ the unit element of $\tilde W$.

\subhead 2.2.5.~Pro-finite-dimensional vector-bundles over
$\Fen$\endsubhead Consider the ind-coherent ind-scheme of
ind-infinite type
$$\tilde\gen\times\Fen=\ind_{w,l}(\tilde\gen_{l}\times\Fen_{w}),$$
where $l\geqslant 0$ and $\tilde\gen_{l}\subset\tilde\gen$ is the sum of all
weight subspaces which do not belong to $\tilde\Delta_{l}^\circ$. Given a
Lie subalgebra $\ben\subset\tilde\gen$ let $\ben_\nil$ denote its
pro-nilpotent radical. Set
$$\dot\nen=\{(x,\ben)\in\tilde\gen\times\Fen;x\in\ben_\nil\}.$$
It is a pro-finite-dimensional vector bundle over $\Fen$. Thus it is
an ind-coherent ind-schemes such that
$$\dot\nen=\ind_w\dot\nen_{w},
\quad\dot\nen_{w}=\dot\nen\cap(\tilde\gen\times\Fen_w),$$ where
$\dot\nen_{w}$ is a compact coherent scheme for each $w$. Define
also
$$\Nen=\dot\nen\cap(\nen\times\Fen).$$
It is an ind-coherent admissible ind-$\tilde I$-scheme such that
$$\Nen=\ind_w\Nen_w,\quad\Nen_w=\dot\nen_w\cap(\nen\times\Fen_w).$$
Note that the $\tilde I$-scheme $\dot\nen_w$ satisfies the property
$(A_{\tilde I})$ because the canonical map $\dot\nen_w\to\Fen_w$ is
$\tilde I$-equivariant, affine, and $\Fen_w$ is normal and
projective, see 1.5.4 and 2.2.1.

\subhead 2.2.6.~Group actions on flag varieties and related
objects\endsubhead
Recall that the ind-group $\tilde G$ acts on the ind-scheme
$\Fen$ by left multiplication.
It
acts also on $\Den=\Fen\times\Fen$ diagonally,
on $\tilde\gen$ by conjugation, and
on $\tilde\gen\times\Fen$ diagonally. For each $w\in\tilde W$ let
$\Den_w\subset\Den$ be the smallest $\tilde G$-stable subset
containing the pair $(\ben_e,\ben_v)$ for each $v\leqslant w$.

Similarly, the group $\Gt$ acts also on $\Fen$, $\Den$, $\tilde\gen$, and
$\tilde\gen\times\Fen$. We simply require that an element
$z\in\CC^\times_\qua$ acts trivially on $\Fen$ and that $z$ acts by
multiplication by $z$ on $\tilde\gen$. This action preserves
$\nen\times\Fen$ and $\dot\nen$, and it restricts to an admissible
$\It$-action on both of them. Note that $\Fen$,
$\tilde\gen\times\Fen$, $\nen\times\Fen$ and $\dot\nen$ are
admissible ind-coherent ind-$\It$-schemes. We also equip $\Xen$ with
the canonical $\It^\op$-action such that $\CC^\times_\qua$ acts
trivially.

For a future use let us introduce the following notation.
Given $\l\in\Xb$ we can view $\theta_\l$ as a one-dimensional
representation of $\It$. Then for each $\It$-scheme $X$ and
each $\It$-equivariant $\Oc_X$-module $\Ec$ we'll write
$\Ec\la\l\ra$ for the $\It$-equivariant $\Oc_X$-module
$$\Ec\la\l\ra=\theta_\l\otimes\Ec.$$

\vskip3mm

\subhead 2.3.~K-theory and the affine flag manifold\endsubhead

\subhead 2.3.1.~Induction of ind-schemes\endsubhead
Recall that the Kashiwara flag manifold is equipped with a canonical
$\tilde I$-torsor $\tilde G_\infty\to\Xen$, where $\tilde G_\infty$ is a
coherent scheme with a $\tilde I^\op\times\tilde I$-action.
For any admissible ind-$\It$-scheme $Z$ we equip the quotient
$$Z_{\Xen}=\tilde G_\infty\times_{\tilde I} Z$$
with the $\It^\op$-action such that the subgroup $\tilde I^\op$ acts
by left multiplication on $\tilde G_\infty$ and $\CC^\times_\qua$
through its action on $Z$. We can regard $Z_\Xen$ as a bundle over
$\Xen$. For any subspace $X\subset\Xen$ let $Z_X$ be the restriction of
$Z_{\Xen}$ to $X$. We'll abbreviate
$Z^{(w)}=Z_{\Xen^w}$ and $Z_{(w)}=Z_{\Fen_w}$ for each
$w\in\tilde W.$
The discussion in Section 1.5.15 yields the following.

\proclaim{2.3.2.~Proposition} Let $Z$ be an ind-coherent admissible
ind-$\It$-scheme. Then $Z_{(w)}$ and $Z_\Fen$ are ind-coherent
admissible ind-$\It$-schemes, and $Z^{(w)}$ is an ind-coherent
admissible ind-$\It^\op$-scheme.
\endproclaim

\noindent Note that $Z_\Xen$ is only a ind$'$-scheme, because $\Xen$
is not quasi-compact. Now, we discuss a few examples which are
important for us.

\subhead 2.3.3.~Examples\endsubhead $(a)$ Set $Z=\Fen$. Consider the
natural projection $p:\Fen_\Fen\to\Fen$ and the
action map $a:\Fen_\Fen\to\Fen$.
The pair $(p,a)$ gives an ind-$\It$-scheme
isomorphism \vskip-3mm
$$\Fen_\Fen\to\Den=\Fen\times\Fen,\quad
(g,\ben)\ \mod\ \tilde I\mapsto(g(\ien),g(\ben)),
$$ where $\tilde I$ acts diagonally on $\Den$. Under
this isomorphism the maps $p, a$ are identified with the
projections $\Den=\Fen\times\Fen\to\Fen$ to the first and the second factors
respectively. Further, the ind-subscheme $(\Fen_w)_\Fen$
is taken to the ind-subscheme $\Den_w\subset\Den$.

\vskip1mm

$(b)$ Taking $Z=\nen\times\dot\nen$ the induction yields an
ind-scheme which is canonically isomorphic to
$\dot\nen\times\dot\nen$.

\vskip1mm

$(c)$ Taking $Z=\Nen$ the induction yields the ind-scheme
$\Nen_\Fen$. We'll abbreviate $\Men=\Nen_\Fen$. By 2.3.3$(a)$ we can
view $\Men$ as the admissible $\It$-equivariant
pro-finite-dimensional vector bundle over $\Den$ whose total space
is
$$\{(x,\ben,\ben')\in
\tilde\gen\times\Den;x\in\ben_\nil\cap\ben'_\nil\}.$$ The
$\It$-action is the diagonal one. In 2.2.5 we have defined $\Nen$ as
an ind-subscheme of $\dot\nen$ and of $\nen\times\Fen$. We may also
regard it as an ind-subscheme of $\nen\times\dot\nen$ by taking a
pair $(x,\ben)\in\Nen$ to the pair
$(x,(x,\ben))\in\nen\times\dot\nen$. Hence we have an inclusion
$\Men\subset\dot\nen\times\dot\nen$ which takes a triple
$(x,\ben,\ben')$ to the pair $((x,\ben), (x,\ben'))$. Composing this
inclusion with the obvious projections
$$q:\dot\nen\times\dot\nen\to\dot\nen\times\Fen,\quad
p:\dot\nen\times\dot\nen\to\Fen\times\dot\nen$$ we can also view
$\Men$ as a good ind-subscheme either of $\dot\nen\times\Fen$ or of
$\Fen\times\dot\nen$. For each $v,w$ we'll write
$$\Men_v=\{(x,\ben,\ben')\in\Men;(\ben,\ben')\in\Den_v\},\quad
\Men_{w,u}=\{(x,\ben,\ben')\in\Men;\ben\in\Fen_w,\ben'\in\Fen_u\}.$$
Note that $(\Nen_v)_\Fen\simeq\Men_v$ and that
$\Men=\ind_{w,u}\Men_{w,u}.$

\vskip1mm

$(d)$ Taking $Z=\nen$ the induction yields a pro-finite-dimensional
vector bundle $\nen^{(w)}$ over $\Xen^w$, and a
pro-finite-dimensional vector bundle $\nen_{(w)}$ over $\Fen_w$ for
each $w$. Note that we have $\dot\nen_w=\nen_{(w)}$, see 2.2.5. For
any integer $l\geqslant 0$ we'll set $\dot\nen^l_w=(\nen^l)_{(w)}$.
The canonical projection $\nen\to\nen^l$ yields a smooth affine
morphism $\dot\nen_w\to\dot\nen^l_w$. Both maps are denoted by the
symbol $p$.

\vskip3mm

\subhead 2.3.4.~Induction of $\It$-equivariant sheaves\endsubhead
Fix an admissible ind-coherent ind-$\It$-scheme $Z$. Consider the
induced ind-scheme $Z^{(w)}$ over $\Xen^w$ for each $w\in\tilde W$.
For any elements $v,w\in\tilde W$ such that $v\leqslant w$ the open
embedding $\Xen^v\subset\Xen^w$ yields an open embedding of
ind-schemes $Z^{(v)}\subset Z^{(w)}$. Fix a closed subgroup
$S\subset\Tt$. We obtain an inverse system of categories
$(\Qcoh^{S}(Z^{(w)}))$, an inverse system of categories
$(\Coh^{S}(Z^{(w)}))$, and an inverse system of $\Rb^{S}$-modules
$\K^{S}(Z^{(w)})$, see 1.3.6. We define
$$\gathered
\Qcoh^{S}(Z_\Xen)=2\pro_w\Qcoh^{S}(Z^{(w)}),\cr
\Coh^{S}(Z_\Xen)=2\pro_w\Coh^{S}(Z^{(w)}),\cr
\K^{S}(Z_\Xen)=\pro_w\K^{S}(Z^{(w)}).\endgathered$$ The discussion
in Section 1.5.15 implies the following

\vskip1mm

\item{$\bullet$}
for each $w\in\tilde W$ the induction yields exact functors
$\Qcoh^{\It}(Z)\to\Qcoh^{S}(Z^{(w)})$ and
$\Coh^{\It}(Z)\to\Coh^{S}(Z^{(w)})$ which commute with tensor
products and a group homomorphism $\K^{\It}(Z)\to\K^{S}(Z^{w}),$

\vskip1mm

\item{$\bullet$}
taking the inverse limit over all $w$'s we get functors
$\Qcoh^{\It}(Z)\to\Qcoh^{S}(Z_\Xen)$,
$\Coh^{\It}(Z)\to\Coh^{S}(Z_\Xen)$ which commute with tensor
products and a group homomorphism $\K^{\It}(Z)\to\K^{S}(Z_\Xen).$
\vskip1mm
\item{$\bullet$}
For each $w\in\tilde W$ and each $\It$-equivariant quasi-coherent
$\Oc_Z$-module $\Ec$ the restriction of the induced
$\Oc_{Z^{(w)}}$-module to the ind-scheme $Z_{(w)}$ is naturally
$\It$-equivariant. Hence the induction yields also functors
$\Qcoh^{\It}(Z)\to\Qcoh^\It(Z_{(w)})$,
$\Coh^{\It}(Z)\to\Coh^\It(Z_{(w)})$, and a group homomorphism
$\K^{\It}(Z)\to\K^\It(Z_{(w)}).$

\vskip1mm

\noindent For each $\Ec\in\Qcoh^\It(Z)$ we'll write $\Ec_\Xen$,
$\Ec_{(w)}$, and $\Ec^{(w)}$ for the induced $\Oc$-modules over
$Z_\Xen$, $Z_{(w)}$ and $Z^{(w)}$ respectively.

\vskip3mm

\subhead 2.3.5.~Examples\endsubhead $(a)$ For each $\l\in\Xb$ let
$\Oc_{\Xen}(\l)$ be the line bundle over $\Xen$ induced from the
character $\theta_\l$. The local sections of $\Oc_\Xen(\l)$ are the
regular functions $f:\tilde G_\infty\to\CC$ such that
$f(xb)=\l(b)f(x)$ for each $x\in\tilde G_\infty$ and $b\in\tilde I$.
Note that $\Oc_\Xen(t)=\Oc_\Xen\la t\ra$, where $t$ is as in 2.1.8.

Restricting $\Oc_{\Xen}(\l)$ to $\Fen$ we get also a line bundle
$\Oc_\Fen(\l)$ over the ind-scheme $\Fen$. We'll write
$\Oc_X(\l)=f^*\Oc_\Xen(\l)$ or $f^*\Oc_\Fen(\l)$ for any map
$f:X\to\Xen$ or $f:X\to \Fen$. For instance we have the line bundles
$\Oc_{\nen_\Xen}(\l)$, $\Oc_{\dot\nen}(\l)$, and $\Oc_{\Nen}(\l)$.
For any $\Oc_X$-module $\Ec$ we'll abbreviate
$$\Ec(\l)=\Ec\otimes_X\Oc_X(\l).$$

\vskip1mm

$(b)$ Given $\l,\mu\in\Xb$ we can consider the $ \It$-equivariant
line bundle $\Oc_\Fen(\mu)\la\l\ra$ over $\Fen$. By induction and
2.3.3$(a)$ it yields an $\It$-equivariant line bundle over $\Den$.
Recall that the $\It$-action on $\Den$ is the diagonal one.
Restricting the induced bundle $(\Oc_\Fen(\mu)\la\l\ra)_\Xen$ over
$\Fen_\Xen$ to $\Fen_\Fen\simeq\Den$ we get the line bundle
$$\Oc_\Den(\l,\mu)=
\Oc_\Fen(\l)\boxtimes\Oc_\Fen(\mu).$$ For any map $Z\to\Den$ we'll
write $\Oc_Z(\l,\mu)=f^*\Oc_\Den(\l,\mu).$ We write also
$\Oc_{\Nen_\Xen}(\l,\mu)=(\Oc_{\Nen}(\mu)\la\l\ra)_\Xen. $

\subhead 2.3.6.~Convolution product on $\K^{\It}(\Nen)$
\endsubhead
The purpose of this section is to define an associative
multiplication
$$\star:\K^{\It}(\Nen)\otimes\K^{\It}(\Nen)\to\K^{\It}(\Nen).$$
Fix $\Ec,\Fc\in\Coh^\It(\Nen)$. Recall that
$$\Coh^\It(\Nen)=2\ind_w\Coh^\It(\Nen_w).$$ Choose $v,w\in\tilde W$
such that $\Ec\in\Coh^\It(\Nen_w)$ and $\Fc\in\Coh^\It(\Nen_v)$. We
can regard $\Ec$ as a coherent $\Oc_{\dot\nen_w}$-module and $\Fc$
as a quasi-coherent $\Oc_{\nen\times\dot\nen_v}$-module. Note that
the closed embedding $\Nen_v\subset\nen\times\dot\nen_v$ is not
good. Fix $u\in\tilde W$ such that the isomorphism
$(\nen\times\dot\nen)_\Fen=\dot\nen\times\dot\nen$ in 2.3.3$(b)$
factors to a good embedding
$$\nu:(\nen\times\dot\nen_v)_{(w)}\to\dot\nen_w\times\dot\nen_u.\leqno(2.3.2)$$
Consider the obvious projections
$$\xymatrix{
\dot\nen_w&\dot\nen_w\times\dot\nen_u\ar[l]_-{f_2}\ar[r]^-{f_1}&\dot\nen_u.}$$
Then $f_2^*(\Ec)$ and $\nu_*(\Fc_{(w)})$ are both $\It$-equivariant
quasi-coherent $\Oc$-modules over $\dot\nen_w\times\dot\nen_u$, and
we can define the following complex in
$\Dcb^\It(\dot\nen_w\times\dot\nen_u)_\qcb$
$$\Gc=f_2^*(\Ec)\Lotimes_{\dot\nen_w\times\dot\nen_u}\nu_*(\Fc_{(w)}).
\leqno(2.3.3)$$
We'll view it as a complex of $\It$-equivariant quasi-coherent
$\Oc$-modules over the ind-scheme $\dot\nen\times\dot\nen$ supported
on the subscheme $\dot\nen_w\times\dot\nen_u$. We want to consider
its direct image by the map $f_1$. Since the schemes $\dot\nen_w$,
$\dot\nen_u$ are not of finite type, this requires some work.

\proclaim{2.3.7.~Proposition} The complex of $\Oc$-modules $\Gc$
over $\dot\nen\times\dot\nen$ does not depend on the choices of
$u,v,w$ up to quasi-isomorphisms. It is cohomologically bounded. Its
direct image $R(f_1)_*(\Gc)$ is a cohomologically bounded
pseudo-coherent complex over $\dot\nen_u$ with cohomology sheaves
supported on $\Nen_u$. The assignment $\Ec\otimes\Fc\mapsto
R(f_1)_*(\Gc)$ yields a group homomorphism
$\star:\K^{\It}(\Nen)\otimes\K^{\It}(\Nen)\to\K^{\It}(\Nen).$
\endproclaim

\noindent{\sl Proof :} We'll abbreviate
$$T=(\nen\times\dot\nen_v)_{(w)},\quad
Y=\dot\nen_w\times\dot\nen_u,\quad\phi_2=f_2\circ\nu,\quad\phi_1=f_1\circ\nu.$$
Thus we have the following diagram
$$\xymatrix{
\dot\nen_w&Y\ar[l]_-{f_2}\ar[r]^-{f_1}&\dot\nen_u\cr
&T.\ar[u]_-\nu\ar[lu]^-{\phi_2}\ar[ru]_-{\phi_1}&}\leqno(2.3.4)$$
The map $\phi_2$ is flat.
We claim that the complex of $\It$-equivariant quasi-coherent
$\Oc_T$-modules
$$\Hc=\phi_2^*(\Ec)\Lotimes_T\Fc_{(w)}$$
is cohomologically bounded.  It is enough to prove that the complex
$for(\Hc)$ is cohomologically bounded.  Since the derived tensor
product commutes with the forgetful functor we may forget the
$\It$-action everywhere. Hence we can use base change and the
projection formula in full generality. To unburden the notation in
the rest of the proof we'll omit the functor $for$.

Now, for an integer $l\geqslant 0$ we have the maps in 2.3.3$(d)$
$$p:\nen\to\nen^l,\quad p:\dot\nen_w\to\dot\nen^l_w.\leqno(2.3.5)$$ Since $\Ec$
is an object of $\Coh(\dot\nen_w)$, by 1.5.13 there is an $l$ and an
object $\Ec^l$ of $\Coh(\dot\nen^l_w)$ such that $\Ec=p^*(\Ec^l)$.
Next, recall that $\Fc$ is an object of
$\Qcoh(\nen\times\dot\nen_v)$. In the commutative diagram
$$\xymatrix{\nen\times\dot\nen_v\ar[r]^{p\times 1}&\nen^l\times\dot\nen_v\ar[r]^{1\times p}&
\nen^l\times\dot\nen^l_v\cr&\Nen_v\ar[lu]\ar[u]&}$$ the right
vertical map is a good inclusion. Thus the $\Oc$-module $(p\times
1)_*(\Fc)$ over $\nen^l\times\dot\nen_v$ is coherent. So if $l$ is
large enough there is a coherent sheaf $\Fc^l$ over
$\nen^l\times\dot\nen^l_v$ such that
$$(p\times 1)_*(\Fc)=(1\times p)^*(\Fc^l).$$ We'll abbreviate
$T^l=(\nen^l\times\dot\nen^l_v)_{(w)}$. Set
$\phi_{2,l}=f_{2,l}\circ\nu_l$ and $\phi_{1,l}=f_{1,l}\circ\nu_l$,
where $\nu_l$, $f_{2,l}$ and $f_{1,l}$ are the obvious inclusion and
projections in the diagram
$$\xymatrix{
\dot\nen^l_w&\dot\nen^l_w\times\dot\nen^l_u\ar[l]_-{f_{2,l}}
\ar[r]^-{f_{1,l}}& \dot\nen^l_u\cr
&T^l.\ar[ul]^{\phi_{2,l}}\ar[ur]_{\phi_{1,l}}\ar[u]_-{\nu_l}&}$$ Let
us consider the complex
$\Hc^l=\phi_{2,l}^*(\Ec^l)\Lotimes_{T^l}\Fc^l_{(w)}$ over $T^l$. The
projection $p$ in (2.3.5) gives a chain of maps
$$\xymatrix{T\ar[r]^-q&(\nen^l\times\dot\nen_v)_{(w)}\ar[r]^-r&T^l}.$$
Note that $Rq_*=q_*$ and $Lr^*=r^*$ because $q$ is affine and $r$ is
flat. Thus a short computation using base change and the projection
formula implies that $$q_*(\Hc)=r^*(\Hc^l).$$ So to prove that $\Hc$
is cohomologically bounded it is enough to prove that $\Hc^l$ itself
is cohomologically bounded. This can be proved using the Kashiwara
affine flag manifold as follows. Write $X^l=\dot\nen^l_w$ and
$$T'=(\nen^l\times\dot\nen^l_v)^{(w)},\quad X'=(\nen^l)^{(w)}.
$$ Consider the Cartesian square
$$\xymatrix{
T'\ar[r]^{\phi_2'}&X'\cr
T^l\ar[u]^i\ar[r]^{\phi_{2,l}}&X^l,\ar[u]_i}$$ where the vertical
maps are the embeddings induced by the inclusion
$\Fen_w\subset\Xen^w$. Recall that $\Ec^l$ is a coherent
$\Oc_{X^l}$-module and that $\Fc^l_{(w)}$ is the restriction to
$T^l$ of the coherent $\Oc_{T'}$-module $\Fc'=(\Fc^l)^{(w)}$. Since
the scheme $\Xen^w$ satisfies the property $(S)$, by 1.2.18 there is
also a Cartesian square
$$\xymatrix{
T^\a\ar[r]^{\phi_{2,\a}}&X^\a\cr
T'\ar[u]^{p_\a}\ar[r]^{\phi'_2}&X'\ar[u]_{p_\a}}$$ where the vertical
maps are smooth and affine, $X^\a$ is smooth of finite type and the
composed maps $j=p_\a\circ i$ are closed embeddings. Further we can
assume that $\Fc'=(p_\a)^*(\Fc^\a)$ for some coherent
$\Oc_{T^\a}$-module $\Fc^\a$. We have
$$i_*(\Hc^l)=(\phi'_2)^*i_*(\Ec^l)\Lotimes_{T'}\Fc'.$$Thus we have also
$$j_*(\Hc^l)=(\phi_{2,\a})^*j_*(\Ec^l)\Lotimes_{T^\a}\Fc^\a.$$
Now $(\phi_{2,\a})^*j_*(\Ec^l)$ and $\Fc^\a$ are both coherent
$\Oc_{T^\a}$-modules and $j_*(\Ec^l)$ is perfect because $X^\a$ is
smooth, see 1.2.11. Hence the complex $j_*(\Hc^l)$ is
pseudo-coherent and cohomologically bounded. Thus the complex
$\Hc^l$ is also pseudo-coherent and cohomologically bounded, because
$j$ is a closed embedding. So $\Hc$ is also cohomologically bounded.

Now we can prove that $\Gc$ and $R(f_1)_*(\Gc)$ are cohomologically
bounded. Once again we can omit the $\It$-action. Since
$R\nu_*=\nu_*$, using the projection formula we get
$\Gc=\nu_*(\Hc)$. Thus the complex $\Gc$ is cohomologically bounded.
Hence $R(f_1)_*(\Gc)$ is also cohomologically bounded because the
derived direct image preserves cohomologically bounded complexes.

To prove that the complex $R(f_1)_*(\Gc)$ is pseudo-coherent it is
enough to observe that we have $R(f_1)_*(\Gc)=p^*R(\phi_{1,l})_*(\Hc^l)$ and
that $\Hc^l$ is pseudo-coherent.

The first claim of the proposition is obvious and is left to the
reader. For instance, since $\Gc=\nu_*(\Hc)$ the complex of
$\Oc$-modules $\Gc$ over the ind-scheme $\dot\nen\times\dot\nen$
does not depend on the choice of $u$. The independence on $v$, $w$
is proved in a similar way.

\qed

\vskip3mm

The following proposition will be proved in 2.4.9 below.

\proclaim{2.3.8.~Proposition} The map $\star$ equips
$\K^{\It}(\Nen)$ with a ring structure.
\endproclaim

\subhead 2.3.9.~Remarks\endsubhead $(a)$ The map $\star$ is an
affine analogue of the convolution product used in \cite{CG}. It is
$\Rb^{\It}$-linear in the first variable (see part $(c)$ below)
but not in the second one.
The definition of $\star$ we have
given here is inspired from \cite{BFM, sec.~7.2}. Observe, however,
that the complex $\Gc$ is not a complex of coherent sheaves over
$\dot\nen\times\dot\nen$, contrarily to what is claimed in loc.~cit.
(in a slightly different setting).

\vskip1mm

$(b)$ Since $\Nen_e\subset\Nen$ is a good subscheme, for each
$\l\in\Xb$ we have the $\It$-equivariant coherent $\Oc_\Nen$-module
$\Oc_{\Nen_e}(\l)=\Oc_{\Nen_e}\la\l\ra$.
Consider the diagram
$$\xymatrix{
\dot\nen_w&\dot\nen_w\times\dot\nen_w\ar[l]_{f_2}\ar[r]^{f_1}&\dot\nen_w\cr
&\dot\nen_w\ar[u]_\delta.&}$$
where $f_1$, $f_2$ are the obvious projections and $\delta$ is the
diagonal inclusion.
Given $\l\in\Xb$ and an object $\Ec$
in $\Coh^{\It}(\Nen_w)$ let $\Ec(\l)$ be the ``twisted" sheaf defined
in 2.3.5$(a)$. We have
$$\aligned
\Ec\star\Oc_{\Nen_e}(\l)
=R(f_1)_*\bigl(f_2^*(\Ec)\,\otimes_{\dot\nen_w\times\dot\nen_w}
\delta_*\Oc_{\dot\nen_w}(\l)\bigr) =\Ec(\l).
\endaligned$$
Thus the associativity of $\star$ yields
$$\Ec\star\Fc(\l)=(\Ec\star\Fc)(\l),\quad\forall\,\Ec,\Fc\in\Coh^\It(\Nen).$$

\vskip1mm

$(c)$ Consider the diagram
$$\xymatrix{
\nen&\nen\times\dot\nen_v\ar[l]_{f_2}\ar[r]^{f_1}&\dot\nen_v\cr
&\Nen_v\ar[u]_\delta&}$$
where $f_1$, $f_2$ are the obvious projections and $\delta$ is the
diagonal inclusion.
Given $\l\in\Xb$ and an object $\Fc$
in $\Coh^{\It}(\Nen_v)$ let $\Fc\la\l\ra$ be the ``twisted" sheaf defined
in 2.2.6.
We have
$$\aligned
\Oc_{\Nen_e}\la\l\ra\star\Fc
=R(f_1)_*\bigl(f_2^*(\Oc_\nen\la\l\ra)
\otimes_{\nen\times\dot\nen_v}\delta_*(\Fc)\bigr)
=\Fc\la\l\ra.
\endaligned$$
Thus the associativity of $\star$ yields
$$\Ec\la\l\ra\star\Fc=(\Ec\star\Fc)\la\l\ra,\quad\forall\,\Ec,\Fc\in\Coh^\It(\Nen).$$

\vskip3mm

\subhead 2.4.~Complements on the concentration in
K-theory\endsubhead

\subhead 2.4.1.~Definition of the concentration map
$\rb_\Sigma$\endsubhead  Let $S\subset \Tt$ be a closed subgroup.
We'll say that $S$ is {\it regular} if the schemes $\nen^S$,
$\Xen^S$ are both locally of finite type. Note that if $S$ is
regular then we have $\Xen^S=\Fen^S$ (as sets, because the lhs is a
scheme of infinite type and locally finite type while the rhs an
ind-scheme of ind-finite type). Next, we'll say that a subset
$\Sigma\subset S$ is {\it regular} if we have $\nen^S=\nen^\Sigma$
and $\Xen^S=\Xen^\Sigma$. In this subsection we'll assume that $S$
and $\Sigma$ are both regular.
Let $\Fen(\a)$, $\a\in A$, be the connected components of $\Fen^S$.
We have $$\dot\nen^S=\bigsqcup_{\a\in A}\dot\nen(\a),$$ where
$\dot\nen(\a)$ is a vector bundle over $\Fen(\a)$ for each $\a$.
Since $S$ is regular we have
$$\gathered
\Nen_\Xen^S=\Men^S=
\bigsqcup_{\a,\b}\Men(\a,\b),\quad\Men(\a,\b)=\Men\cap(\dot\nen(\a)\times\dot\nen(\b)),\endgathered$$
where $\Men$ is as in 2.3.3. Here we have abbreviated
$\Nen^S_\Xen=(\Nen_\Xen)^S$. We define
$$\K^S(\Nen_\Xen^S)=\pro_w\ind_v\K^S\bigl((\Nen_v)^{(w),S}\bigr)=
\prod_\a\bigoplus_\b\K^S(\Men(\a,\b)).$$ Observe that in 2.3.4 we
have defined the group $\K^S(\Nen_\Xen)$ in a similar way by setting
$$\K^S(\Nen_\Xen)=\pro_w\ind_v\K^S\bigl((\Nen_v)^{(w)}\bigr).$$
Now we can define the concentration map. Consider the closed
embeddings
$$\xymatrix{(\nen\times\Fen)^S_\Xen\ar[r]^-i&\nen_\Xen^S\times_\Xen\Fen_\Xen\ar[r]^-j&
(\nen\times\Fen)_\Xen=\nen_\Xen\times_\Xen\Fen_\Xen.}\leqno(2.4.1)$$
The scheme $\nen_\Xen$ is pro-smooth and the inclusion
$\Nen\subset\nen\times\Fen$ is good. Thus 1.5.20 yields a group
homomorphism
$$
\pmb\g_\Sigma=(i_*)^{-1}\circ
Lj^*:\K^S(\Nen_\Xen)\to\K^S(\Nen_\Xen^S)_\Sigma.$$ Composing it with
the induction  $\Gamma:\Ec\mapsto\Ec_\Xen$ yields a group
homomorphism
$$\rb_\Sigma:\xymatrix{\K^{\It}(\Nen)\ar[r]&\K^{S}(\Nen_\Xen)\ar[r]&\K^S(\Nen_\Xen^S)_\Sigma}.\leqno(2.4.2)$$
The map $\rb_\Sigma$ is
called the {\it concentration map}.

\subhead 2.4.2.~Remark\endsubhead The map $\rb_\Sigma$ is an affine
analogue of the concentration map defined in \cite{CG,
thm.~5.11.10}. It can also be described in the following way. Let
$\Ec$ be an $\It$-equivariant coherent $\Oc$-module over $\Nen$. Fix
$v\in\tilde W$ such that $\Ec\in\Coh^\It(\Nen_v)$. Given any
$w\in\tilde W$ we fix $\nu$, $u$ as in (2.3.2). Under the direct
image by $\nu$ we can view the $S$-equivariant coherent $\Oc$-module
$\Ec_{(w)}$ as a $S$-equivariant quasi-coherent $\Oc$-module over
$\dot\nen_w\times\dot\nen_u$ supported on $\Men_v\cap\Men_{w,u}$.
The obvious projection
$$q:\dot\nen_w\times\dot\nen_u\to\dot\nen_w\times\Fen_u$$
yields a good inclusion $\Men_{w,u}\subset\dot\nen_w\times\Fen_u.$
Thus, under the direct image by $q$ we can also regard $\Ec_{(w)}$
as a $S$-equivariant coherent $\Oc$-module over $\dot\nen_w\times\Fen_u$.
Then we
consider the following chain of inclusions
$$\xymatrix{\dot\nen_w^S\times\Fen_u^S\ar[r]^i&
\dot\nen_w^S\times\Fen_u\ar[r]^j&\dot\nen_w\times\Fen_u.}\leqno(2.4.3)$$
Since $\Ec_{(w)}$ is flat over $\Fen_w$ and since
$\dot\nen_w\to\Fen_w$ is a pro-finite-dimensional vector bundle, we
have a cohomologically bounded pseudo-coherent complex
$\Ec'=Lj^*(\Ec_{(w)})$ over $\dot\nen_w^S\times\Fen_u$. Next, the
Thomason theorem yields an invertible map
$$i_*:\K^S(\dot\nen^S_w\times\Fen^S_u)_\Sigma\to\K^S(\dot\nen^S_w\times\Fen_u)_\Sigma.$$
Thus we have a well-defined element $\Ec''=(i_*)^{-1}(\Ec')$. It can
be regarded as an element of $\K^S(\Men_{u,w}^S)_\Sigma$ for a
reason of supports. If $w,u$ are large enough then $\Men(\a,\b)$ is
a closed and open subset of $\Men^S_{u,w}$. The component of
$\rb_\Sigma(\Ec)$ in $\K^S(\Men(\a,\b))_\Sigma$ is the restriction
of $\Ec''$ to $\Men(\a,\b)$.

\vskip2mm

\proclaim{2.4.3.~Proposition} If $S=\Sigma=\Tt$ then $\rb_S$ is an
injective map.
\endproclaim

\noindent{\sl Proof :} We have $\Xen^e=N^\op=\tilde I^\op/\tilde T$
as a $\tilde I^\op$-scheme. Thus we have $\Nen^{(e)}=N^\op\times
\Nen$. The induction yields an inclusion
$$\K^{\It}(\Nen)\to\K^{\Tt}(\Nen^{(e)}),\ \Ec\mapsto\Ec^{(e)}.\leqno(2.4.4)$$
Therefore the induction map
$\K^{\It}(\Nen)\to\K^\Tt(\Nen_\Xen)$ is also injective, because
composing it with the canonical map
$$\K^{\Tt}(\Nen_{\Xen})=
\pro_w\K^{\Tt}(\Nen^{(w)}) \to \K^{\Tt}(\Nen^{(e)})$$ yields
(2.4.4). Thus, to prove that $\rb_\Tt$ is injective it is enough to
check that the canonical map
$\K^\Tt(\Nen^{(e)})\to\K^\Tt(\Nen^{(e)})_\Tt$ is injective. This is
obvious because the $\Rb^\Tt$-module $\K^\Tt(\Nen^{(e)})$ is
torsion-free (use an affine cell decomposition of $\Nen$).

\qed

\subhead 2.4.4.~ Concentration of $\Oc$-modules supported on
$\Nen_e$
\endsubhead
Let $\Ec$ be an $\It$-equivariant vector bundle over $\Nen_e$. Since
the inclusion $\Nen_e\subset\nen\times\Fen$ is good we may view
$\Ec$ as an object of $\Coh^\It(\nen\times\Fen)$. Now we consider
the diagrams (2.4.1) and (2.4.2). The induced coherent sheaf
$\Gamma(\Ec)=\Ec_\Xen$ is flat over ${\nen_\Xen}$. Thus we have
$Lj^*(\Ec_\Xen)=j^*(\Ec_\Xen)$. Thus we obtain
$$\rb_{\Sigma}(\Ec)=(i_*)^{-1}j^*(\Ec_\Xen).$$
Next we have $j^{-1}((\Nen_e)_\Xen)=i((\Nen_e)^S_\Xen)$. This
implies that
$$\rb_{\Sigma}(\Ec)=j^*(\Ec_\Xen).$$
Therefore we have proved the following.

\proclaim{2.4.5.~Proposition} If $\Ec$ is an $\It$-equivariant
vector bundle over $\Nen_e$ then $\rb_{\Sigma}(\Ec)$ is the
restriction of the coherent sheaf $\Ec_\Xen$ over $\Nen_\Xen$ to the
fixed-points subscheme $\Nen^S_\Xen$.
\endproclaim

\subhead 2.4.6.~ Concentration of $\Oc$-modules supported on
$\Nen'_{s_\a}$
\endsubhead
Fix a simple affine root $\a\in\tilde\Pi$. Recall that $s_\a$ is the
corresponding simple reflection and that $\nen^\op_{s_\a}$ is a
1-dimensional $\Tt$-module whose class in the ring $\Rb^{\Tt}$ is
$\theta_{t+\a}$. Recall also that $\Nen\subset\nen\times\Fen$ and
that $\nen_{s_\a}\subset\ben_\nil$ are good inclusions for each
$\ben\in\Fen_{s_\a}$. So we have a good $\It$-equivariant subscheme
$\Nen'_{s_\a}\subset\Nen$ given by
$$\Nen'_{s_\a}=\nen_{s_\a}\times\Fen_{s_\a}.$$
By a good subscheme we means that $\Nen'_{s_\a}$ is a good
ind-subscheme, as in 1.3.7, which is a scheme. Note that
$\Nen'_{s_\a}$ is pro-smooth, because it is a pro-finite-dimensional
vector bundle over the smooth scheme $\Fen_{s_\a}$. Let $\Ec$ be an
$\It$-equivariant vector bundle over $\Nen'_{s_\a}$. We'll view it
as a $\It$-equivariant coherent $\Oc$-module over $\Nen$ or
$\nen\times\Fen$. The purpose of this section is to compute the
element $\rb_{\Sigma}(\Ec)$.

First we assume that $S=\Tt$. Consider the diagrams (2.4.1) and
(2.4.2). The coherent sheaf $\Gamma(\Ec)=\Ec_\Xen$ is flat over
$(\Nen'_{s_\a})_\Xen$. So it is also flat over $(\nen_{s_\a})_\Xen$.
However it is not flat over $\nen_\Xen$. To compute $Lj^*(\Ec_\Xen)$
we need a resolution of $\Ec_\Xen$ by flat
$\Oc_{\nen_\Xen}$-modules. For this it is enough to construct a
resolution of $\Ec$ by flat $\Oc_\nen$-modules, and to apply
induction to it. We have a closed immersion
$$\Nen'_{s_\a}\subset\Nen''_{s_\a},\quad
\Nen''_{s_\a}=\nen\times\Fen_{s_\a}.$$ The Koszul resolution of
$\Oc_{\Nen'_{s_\a}}$ by locally-free $\Oc_{\Nen''_{s_\a}}$-modules
is the complex
$$\Lambda_{\Nen''_{s_\a}}(\a)=
\Bigl\{\Oc_{\Nen''_{s_\a}}\la t+\a\ra\to\Oc_{\Nen''_{s_\a}}\Bigr\}$$
situated in degrees $[-1,0]$. We may assume that
$$\Ec=\Oc_{\nen_{s_\a}}\boxtimes\Fc,$$ where $\Fc$ is a
$\It$-equivariant locally free $\Oc_{\Fen_{s_\a}}$-module. Set
$$\Ec'=\Oc_{\nen}\boxtimes\Fc.$$ It is a $\It$-equivariant locally free
$\Oc_{\Nen''_{s_\a}}$-module whose restriction to $\Nen'_{s_\a}$ is
equal to $\Ec$. We have
$$\aligned
\rb_{\Sigma}(\Ec)
&=(i_*)^{-1}Lj^*\Gamma(\Ec'\otimes_{\Oc_{\Nen''_{s_\a}}}\Lambda_{\Nen''_{s_\a}}(\a)),\cr
&=(i_*)^{-1}j^*\Gamma(\Ec')-(i_*)^{-1}j^*\Gamma(\Ec'\la t+\a\ra).
\endaligned$$
Since $S=\Tt$ we have
$j^{-1}((\Nen'_{s_\a})_\Xen)=j^{-1}((\Nen''_{s_\a})_{\Xen})$. Thus,
for each $\Oc_{\Nen''_{s_\a}}$-module $\Fc$ we have
$j^*\Gamma(\Fc)=j^*\Gamma(\Fc|_{\Nen'_{s_\a}})$. This implies that
$$\rb_{\Sigma}(\Ec)=
(i_*)^{-1}j^*\Gamma(\Ec)-(i_*)^{-1}j^*\Gamma(\Ec\la t+\a\ra).
$$
Next, observe that $\Nen_\Xen^S=\Den^S$. Thus the map
$$i:(\Nen'_{s_\a})_\Xen^S\to j^{-1}((\Nen'_{s_\a})_\Xen)$$ is equal to the obvious inclusion
$$\Den_{s_\a}^S\subset\Den'_{s_\a},\quad\Den'_{s_\a}=
\Den_{s_\a}\cap(\Fen^S\times\Fen).$$ Now, we have an exact sequence
of $\Oc_{\Den'_{s_\a}}$-modules
$$0\to\Oc_{\Den'_{s_\a}}(0,-\a)\to\Oc_{\Den'_{s_\a}}\to\Oc_{\Den^S_{s_\a}}\to 0.$$
Therefore we have
$$\rb_{\Sigma}(\Ec)=(1-\Oc_{\Nen_\Xen^S}(\a+t,0))\,
(1-\Oc_{\Nen_\Xen^S}(0,-\a))^{-1}\Ec_\Xen|_{\Nen_\Xen^S},$$ where
$\Ec_\Xen|_{\Nen_\Xen^S}$ is the restriction to $\Nen_\Xen^S$
of the induced sheaf $\Ec_\Xen$ over $\Nen_\Xen$.

For any $S\subset\Tt$ we obtain in the same way the following
formula, compare \cite{VV, (2.4.6)}.

\proclaim{2.4.7.~Proposition} We have
$$\rb_{\Sigma}(\Ec)=\cases
(1-\Oc_{\Nen_\Xen^S}(\a+t,0))\,(1-\Oc_{\Nen_\Xen^S}(0,-\a))^{-1}
\Ec_\Xen|_{\Nen_\Xen^S}&\roman{if}\ \theta_\a\neq 1,t,\cr
(1-\Oc_{\Nen_\Xen^S}(\a+t,0))\, \Ec_\Xen|_{\Nen_\Xen^S}&\roman{if}\
\theta_\a=1\neq t,\cr (1-\Oc_{\Nen_\Xen^S}(0,-\a))^{-1}
\Ec_\Xen|_{\Nen_\Xen^S}&\roman{if}\ \theta_\a=t\neq 1,\cr
\Ec_\Xen|_{\Nen_\Xen^S}&\roman{if}\ \theta_\a=t=1.
\endcases$$

\endproclaim

\subhead 2.4.8.~Multiplicativity of $\rb_\Sigma$\endsubhead Let
$S\subset\Tt$ be a regular closed subgroup. We have
$$\gathered
\Men^S=\ind_{w,u}\Men_{w,u}^S,\quad
\K(\Men^S)=\ind_{w,u}\K(\Men_{w,u}^S)=\bigoplus_{\a,\b}\K(\Men(\a,\b)),
\endgathered$$ where
$\Men(\a,\b)$ is as in 2.4.1. We have also $$
\K(\Nen^S_\Xen)=\prod_\a\bigoplus_\b\K(\Men(\a,\b)).$$ Therefore the
group $\K(\Nen^S_\Xen)$ can be regarded as the completion of
$\K(\Men^S)$. Note that $\Men(\a,\b)$ is a closed subscheme of
$\dot\nen(\a)\times\dot\nen(\b)$ and the later is smooth and of
finite type because $S$ is regular. So $\K(\Men^S)$,
$\K(\Nen_\Xen^S)$ are both equipped with an associative convolution
product. See Section 3.1 and the proof of the proposition below for
details.

\proclaim{2.4.9.~Proposition} The map $\star$ yields a ring
structure on $\K^{\It}(\Nen)$. If the group $S$ is regular then the
map $\rb_{\Sigma}:\K^{\It}(\Nen)\to\K^S(\Nen^S_\Xen)_{\Sigma}$ is a
ring homomorphism.
\endproclaim

\noindent{\sl Proof :} Since the group $S$ acts trivially on
$\Nen^S_\Xen$ we have
$$\K^S(\Nen^S_\Xen)_{\Sigma}=\K(\Nen^S_\Xen)\otimes\Rb^S_{\Sigma}.$$
The multiplication on the lhs is deduced by base change from the
product on $\K(\Nen^S_\Xen)$ mentioned above. It is enough to check
that we have $$\rb_{\Sigma}(x\star
y)=\rb_{\Sigma}(x)\star\rb_{\Sigma}(y),\quad\forall x,y.$$ Indeed,
setting $S=\Sigma=\Tt$, this relation and 2.4.3 imply that
$\K^{\It}(\Nen)$ is a subring of
$\K(\Nen^{S}_\Xen)\otimes\Rb^S_{S}$.

Fix $v, w\in\tilde W$. Let $u\in\tilde W$ be as in 2.3.6.
Fix $\Ec\in\Coh^{\It}(\Nen_w)$ and
$\Fc\in\Coh^{\It}(\Nen_v)$. Recall that $\Ec$, $\Fc$ denote also the
corresponding classes in $\K^{\It}(\Nen)$ and that $\Ec\star\Fc$ is
the class of an $\It$-equivariant cohomologically bounded
pseudo-coherent complex over $\Nen_u$. Let us recall the construction
of this complex. We'll regard $\Ec$ as an $\It$-equivariant coherent
$\Oc_{\dot\nen_w}$-module and $\Fc$ as an $\It$-equivariant
quasi-coherent $\Oc_{\nen\times\dot\nen_v}$-module.
Consider the diagram (2.3.4) that we reproduce below for the
comfort of the reader
$$\xymatrix{
\dot\nen_w&Y\ar[l]_-{f_2}\ar[r]^-{f_1}&\dot\nen_u\cr
&T.\ar[u]_-\nu&}$$
The map $f_2$ is flat and we have
$$\Ec\star\Fc=R(f_1)_*(\Gc),
\quad\Gc=f_2^*(\Ec)\Lotimes_{Y}\nu_*(\Fc_{(w)}).$$

We want to compute $\rb_\Sigma(\Ec\star\Fc)$. Fix an
element $x\in\tilde W$. First, we consider the induced complex
$(\Ec\star\Fc)_{(x)}$ over $(\dot\nen_u)_{(x)}$. Under induction the
maps $f_1$, $f_2$ yield flat morphisms
$$\xymatrix{
(\dot\nen_w)_{(x)}&Y_{(x)}\ar[l]_-{f_{2,(x)}}\ar[r]^-{f_{1,(x)}}
&(\dot\nen_u)_{(x)}.}$$
The induction is exact and it commutes with
tensor products. Thus we have
$$\gathered
(\Ec\star\Fc)_{(x)}=R(f_{1,(x)})_*\bigl(
f_{2,(x)}^*(\Ec_{(x)})\Lotimes_{Y_{(x)}}(\nu_*\Fc_{(w)})_{(x)}\bigr).
\endgathered$$
Fix $y,z\in\tilde W$ such that the canonical isomorphisms
$\dot\nen_\Fen=\Fen\times\dot\nen$ and
$(\nen\times\dot\nen)_\Fen=\dot\nen\times\dot\nen$
yield inclusions
$$\l:(\dot\nen_w)_{(x)}\to\Fen_x\times\dot\nen_y,\quad
\mu:(\dot\nen_u)_{(x)}\to\Fen_x\times\dot\nen_z,\quad
\nu:(\nen\times\dot\nen_v)_{(y)}\to\dot\nen_y\times\dot\nen_z.$$ We
put
$$\Gc'=\mu_*((\Ec\star\Fc)_{(x)}),\quad\Ec'=\l_*(\Ec_{(x)}),\quad\Fc'=\nu_*(\Fc_{(y)}).$$
We have
$$\Gc'=R(\pi_2)_*\bigl(\pi_3^*(\Ec')\Lotimes_{Y'}\pi_1^*(\Fc')\bigr),
\leqno(2.4.5)$$ where $Y'=\Fen_x\times\dot\nen_y\times\dot\nen_z$
and $\pi_1$, $\pi_2$, $\pi_3$ are the obvious projections
$$\xymatrix{
\dot\nen_y\times\dot\nen_z&Y'\ar[l]_-{\pi_1}\ar[r]^-{\pi_2}\ar[d]_-{\pi_3}&
\Fen_x\times\dot\nen_z\cr &\Fen_x\times\dot\nen_y.&}$$ As explained
in 2.4.2 we can regard the complex $\Gc'$, which is supported on
$\Men_{x,z}$, as a complex over $\dot\nen_x\times\Fen_z$. Let
$\Gc''$ denote the later. We have
$$\rb_{\Sigma}(\Ec\star\Fc)=\pmb\g_\Sigma(\Gc'').$$ Now we compute
$\Gc''$. Applying the base change formula to (2.4.5) we obtain the
following equality in $\K^S(\dot\nen_x\times\Fen_z)$
$$\Gc''=R(p_2)_*\bigl(p_3^*(\Ec'')\Lotimes_{Y''}
p_1^*(\Fc'')\bigr).\leqno(2.4.6)$$ Here
$Y''=\dot\nen_x\times\dot\nen_y\times\Fen_z$ and $p_1$, $p_2$, $p_3$
are the projections
$$\xymatrix{
\dot\nen_y\times\Fen_z&
Y''\ar[l]_-{p_1}\ar[r]^-{p_2}\ar[d]_-{p_3}& \dot\nen_x\times\Fen_z\cr
&\dot\nen_x\times\dot\nen_y.&}$$ Further $\Ec''$, $\Fc''$
are $S$-equivariant quasi-coherent $\Oc$-modules over
$\dot\nen_x\times\dot\nen_y$,
$\dot\nen_y\times\Fen_z$ respectively which are characterized by the
following properties
$$p_*(\Ec'')=\Ec',\quad\Fc''=q_*(\Fc'),\quad
\Ec''\text{\ is\ supported\ on\ } \Men_{x,y},$$ where $p$, $q$
are the obvious maps
$$\xymatrix{\Fen_w\times\dot\nen_u&
\dot\nen_w\times\dot\nen_u\ar[r]^-{q}\ar[l]_-p&\dot\nen_w\times\Fen_u}.$$

Now, recall that we must prove that the following formula holds in
$\K^S(\Nen^S_\Xen)_\Sigma$
$$\rb_{\Sigma}(\Ec\star\Fc)=
\rb_{\Sigma}(\Ec)\star\rb_{\Sigma}(\Fc).$$ Let $i$, $j$ be as in
2.4.2. The lhs is
$$\rb_{\Sigma}(\Ec\star\Fc)=
\pmb\g_\Sigma(\Gc'')=(i_*)^{-1} Lj^*(\Gc'').\leqno(2.4.7)$$ Let us
describe the rhs. We'll abbreviate
$$N=\dot\nen^S,\quad M=\Men^S.$$
Both are regarded
as ind-schemes of ind-finite type. Thus we have
$$\K(N)=\bigoplus_\a\K(\dot\nen(\a)),\quad
\K(M)=\bigoplus_{\a,\b}\K(\Men(\a,\b)).$$
Note that $\K(M)=\K(N^2\on M)$ for
the inclusion
$M\subset N^2$ given by
$$(x,\ben,\ben')\mapsto
\bigl((x,\ben), (x,\ben')\bigr).$$ Given $a=1,2,3$ let $q_a:N^3\to
N^2$ be the projection along the $a$-th factor. We define the
convolution product on $\K(M)$ by
$$x\star y=
R(q_2)_*\bigl(q_3^*(x)\Lotimes_{N^3} q_1^*(y)\bigr),\quad\forall
x,y\in\K(M).\leqno(2.4.8)$$ Note that $N$ is a disjoint union of
smooth schemes of finite type and that $q_1,q_3$ are flat maps.
We'll use another expression for $\star$. For this, we write
$$F=\Fen^S,\quad NF=N\times F,\quad N^2F=N\times N\times F.$$
The obvious projections below are flat morphisms
$$\xymatrix{
NF&N^2F\ar[l]_-{\rho_1}\ar[r]^-{\rho_2}\ar[d]_-{\rho_3}& NF\cr
&N^2.&}$$ We have also $\K(M)=\K(NF\on M)$
for the inclusion
$M\subset NF$ given by
$$(x,\ben,\ben')\mapsto \bigl(\ben,(x,\ben')\bigr).$$ The projection
formula yields
$$x\star y=
R(\rho_2)_*\bigl(\rho_3^*(x)\Lotimes_{N^2F}\rho_1^*(y)\bigr).
\leqno(2.4.9)$$ Note that $N^2F$ is also a disjoint union of smooth
schemes of finite type. Finally we must compute $\rb_\Sigma(\Ec)$
and $\rb_\Sigma(\Fc)$. Once again, as explained in 2.4.2, we must
first regard $\Ec'$, $\Fc'$ as complexes of $\Oc$-modules over
$\dot\nen_x\times\Fen_y$, $\dot\nen_y\times\Fen_z$ respectively, and
then we apply the map $(i_*)^{-1}\circ Lj^*$ to their class in
K-theory.

Now, by (2.4.6), (2.4.7) and (2.4.9) we are reduced to prove the
following equality
$$\gathered
(i_*)^{-1}Lj^*R(p_2)_*\bigl(p_3^*(\Ec'')\Lotimes_{Y''}
p_1^*(\Fc'')\bigr)=\cr
=R(\rho_2)_*\bigl(\rho_3^*(i'_*)^{-1}L(j')^*(\Ec'')\Lotimes_{N^2F}
\rho_1^*(i_*)^{-1}Lj^*(\Fc'')\bigr).\endgathered$$
Here $i$, $i'$, $j$ and $j'$ are the obvious inclusions in the following diagram
$$\xymatrix{
\dot\nen^S\times\Fen^S\ar[r]^{i}&
\dot\nen^S\times\Fen\ar[r]^{j}&\dot\nen\times\Fen\cr
\dot\nen^S\times\dot\nen^S\ar[r]^{i'}\ar[u]_q&
\dot\nen^S\times\dot\nen\ar[r]^{j'}\ar[u]_q&\dot\nen\times\dot\nen\ar[u]_{q}.}$$
This is an easy consequence of the base change and of the projection formula.

\qed

\vskip3mm

\subhead 2.5.~Double affine Hecke algebras
\endsubhead

\subhead 2.5.1.~Definitions\endsubhead First, let us introduce the
following notation : given any commutative ring $\Ab$ we'll write
$\Ab_t=\Ab[t,t^{-1}]$ and $\Ab_{q,t}=\Ab[q,q^{-1},t,t^{-1}]$. Recall
that $G$ is a simple, connected and simply connected linear group over $\CC$.
The double affine Hecke algebra (=DAHA) associated to $G$ is the
associative $\ZZ_{q,t}$-algebra $\Hb$ with 1 generated by the
symbols $T_w$, $X_\l$ with $w\in\tilde W$, $\l\in\tilde\Xb$ such
that the $T_w$'s satisfy the braid relations of $\tilde W$ and such that
$$\aligned
&X_\delta=q,\quad X_\mu X_\l=X_{\l+\mu},
\quad (T_{s_\a}-t)(T_{s_\a}+1)=0,
\hfill\\
&X_\l T_{s_\a}-T_{s_\a}X_{\l-r\a}=
(t-1)X_\l(1+X_{-\a}+...+X_{-\a}^{r-1}) \quad \roman{if\ }
\la\l,\alphav\ra=r\ge 0. \hfill
\endaligned$$
Here $\a$ is any simple affine root. Let $\Hb^f\subset\Hb$ be the
subring generated by $t$ and the $T_w$'s with $w\in W$. Let
$\Rb\subset\Hb$ be the subring generated by $\{X_\l\,
;\,\l\in\tilde\Xb\}$. To avoid confusions we may write $\Rb_X=\Rb$.
Finally let $\Rb_Y\subset\Hb$ be the subring generated by
$\{Y_{\check\l}\, ;\,\check\l\in\tilde\Yb\}$, where
$Y_{\check\l}=
T_{\xi_{\check\l_1}}
T_{\xi_{\check\l_2}}^{-1}$
with $\check\l=\check\l_1-\check\l_2$ and $\check\l_1$, $\check\l_2$ dominant.
The following
fundamental result has been proved by Cherednik. We'll refer to it
as the PBW theorem for $\Hb$.

\proclaim{2.5.2.~Proposition} The multiplication in $\Hb$ yields
$\ZZ_{q,t}$-isomorphisms
$$\Rb\otimes\Hb^f\otimes\Rb_Y\to\Hb,\quad
\Rb_Y\otimes\Hb^f\otimes\Rb\to\Hb.$$ The $\ZZ_{q,t}$-algebra $\Hb^f$
is isomorphic to the Iwahori-Hecke algebra (over the commutative
ring $\ZZ_{t}$) associated to the Weyl group $W$. The rings $\Rb$,
$\Rb_Y$ are the group-rings associated to the lattices $\tilde\Xb$,
$\tilde\Yb$ respectively.
\endproclaim

\subhead 2.5.3.~Remark\endsubhead The algebra $\Hb$ is the one
considered in \cite{V}. It is denoted by the symbol $\hat\Hb$ in
\cite{VV}. Note that we have
$X_{\omega_0}T_{s_0}X_{\omega_0}^{-1}=X_{\a_0}T_{s_0}^{-1}.$ Thus
$\Hb$ is a semidirect product $\CC[X_{\omega_0}^{\pm
1}]\ltimes\Hc(\tilde W,\Xb\oplus\ZZ\delta)$ with the notation in
\cite{H, sec.~5}. Changing the lattices in the definition of $\Hb$
yields different versions of the DAHA whose representation theory is
closely related to the representation theory of $\Hb$. These
different algebras are said to be {\it isogeneous}. In this paper
we'll only consider the case of $\Hb$ to simplify the exposition.
For more details the reader may consult \cite{VV, sec.~2.5}.

\vskip2mm

Let
$\Ocb(\Hb)$ be the category of all right $\CC\Hb$-modules
which are finitely generated, locally finite over $\Rb$
(i.e., for each element $m$ the $\CC$-vector space $m\Rb$
is finite dimensional), and such
that $q,t$ act by multiplication by a complex number.
It is an Abelian category. Any object has a finite
length. For any module $M$ in $\Ocb(\Hb)$ we have
$$M=\bigoplus_{h\in\tilde T} M_h,\quad
M_h=\bigcap_{\l\in\tilde\Xb}\bigcup_{r\geqslant 0} \{m\in
M;m(X_\l-\l(h))^r=0\}.$$ We'll call $M_h$ the {\it $h$-weight subspace}.
It is finite dimensional.
Next, we set $$\widehat M=\prod_h M_h.$$ The vector
space $\widehat M$ is equipped with
the product topology, the $M_h$'s being equipped with the discrete
topology. Note that $M\subset\widehat M$ is a dense subset.
The $\CC\Hb$-action on $M$ extends uniquely to a continuous $\CC\Hb$-action on
$\widehat M$.

Fix an element $h=(s,\tau)$ of $\tilde T$, i.e., we let $s\in
T_\cent$ and $\tau\in\CC^\times_\rot$. For each
$\zeta\in\CC^\times_\qua$ we can form the corresponding tuple
$(h,\zeta)\in\Tt$. Let $\Ocb_{h,\zeta}(\Hb)$ be the full subcategory
of $\Ocb(\Hb)$ consisting of the modules $M$ such that $q=\tau$,
$t=\zeta$ and $M_{h'}=0$ if $h'$ is not in the orbit of $h$
relatively to the $\tilde W$-action on $\tilde T$ in 2.1.10. Let
$\tilde W$ act on $\Tt$ so that it acts on $\tilde T$ as in 2.1.10
and it acts trivially on $\CC^\times_\qua$. We have
$$\Ocb(\Hb)=\bigoplus_{h,\zeta}\Ocb_{h,\zeta}(\Hb),$$
where $(h,\zeta)$ varies in a set of representatives of the $\tilde
W$-orbits \cite{VV, lem.~2.1.3, 2.1.6}.

\subhead 2.5.4.~Geometric construction of the DAHA\endsubhead We can
now give a geometric construction of the $\ZZ_{q,t}$-algebra $\Hb$.
First, let us introduced a few more notations. For each $\l\in\Xb$
we consider the following element of $\K^{\It}(\Nen)$
$$x_\lambda=\Oc_{\Nen_e}(\l)=\Oc_{\Nen_e}\la\l\ra.$$
Next, given a simple affine root $\a\in\tilde\Pi$ we have the good
$\It$-equivariant subscheme $\Nen'_{s_\a}\subset\Nen$ introduced in
2.4.6. For each weights $\l,\mu\in\Xb$ we define the
$\It$-equivariant coherent sheaf $\Oc_{\Nen'_{s_\a}}(\l,\mu)$ over
$\Nen$ as the direct image of the $\It$-equivariant vector bundle
$\Oc_{\Nen'_{s_\a}}(\mu)\la\l\ra$ over $\Oc_{\Nen'_{s_\a}}$, see 2.3.5$(b)$.
Assume further that
$$\l+\mu=-\a,\quad\la\l,\check\a\ra=\la\mu,\check\a\ra=-1.\leqno(2.5.1)$$
Then we consider the following element of $\K^{\It}(\Nen)$ given by
$$t_{s_\a}=-1-\Oc_{\Nen'_{s_\a}}(\l,\mu).$$

\proclaim{2.5.5.~Lemma} The element $t_{s_\a}$ is independent of the
choice of $\l$, $\mu$ as above.
\endproclaim

\noindent{\sl Proof :} It is enough to observe that if
$\la\l',\check\a\ra=0$ then the $\It$-equivariant line bundle
$\Oc_{\Nen'_{s_\a}}(\l',-\l')$ is trivial.

\qed

\vskip3mm

The assignment $\theta_\l\mapsto X_\l$ identifies $\Rb^{\tilde T}$
with the ring $\Rb=\Rb_X$, and $\Rb_t=\Rb^\Tt$ with the subring of
$\Hb$ generated by $t$ and $\Rb$. Now we can prove the main result
of this section.

\proclaim{2.5.6.~Theorem} There is an unique ring isomorphism
$\Phi:\Hb\to\K^{\It}(\Nen)$ such that $T_{s_\a}\mapsto t_{s_\a}$ and
$X_\l\mapsto x_\l$ for each $\a\in\tilde\Pi$, $\l\in\tilde\Xb$.
Under $\Phi$ and the forgetting map $\Rb^\It=\Rb^\Tt$, the canonical
(left) $\Rb^\It$-action on $\K^\It(\Nen)$ is identified with the canonical
(left) $\Rb_t$-action on $\Hb$.
\endproclaim

\noindent{\sl Proof :} First we
prove that the assignment $$T_{s_\a}\mapsto t_{s_\a},\quad
X_\l\mapsto x_\l,\quad \forall\a\in\tilde\Pi, \l\in\tilde\Xb$$
yields a $\ZZ_{q,t}$-algebra homomorphism
$$\Phi:\Hb\to\K^{\It}(\Nen).$$
We must check that the elements $t_{s_\a}$, $x_\l$ satisfy the
defining relations of $\Hb$. To do so let $S=\Sigma=\Tt$ and
consider the group homomorphism
$$\rb_S:\K^{\It}(\Nen)\to\K^S(\Nen^S_\Xen)_S.$$
Note that $\Nen_\Xen^S=\Den^S$, because $S=\Tt$.
Thus we have a $\Rb^S$-module isomorphism
$$\aligned
\K^S(\Nen^S_\Xen) &=\pro_w\ind_v\K^S\bigl((\Den_v)_{(w)}^S\bigr)\cr
&= \pro_w\ind_u\K^S(\Den_{w,u}^S)\cr &=
\prod_{w}\bigoplus_{u}\Rb^S\xb_{w,u}.
\endaligned\leqno(2.5.2)
$$
Here the symbol $\xb_{w,u}$ stands for the fundamental class of the
fixed point $(\ben_w,\ben_u)$, see 2.2.4$(b)$. The convolution product is
$\Rb^S$-linear and is given by
$$\xb_{v,w}\star\xb_{y,z}=\cases\xb_{v,z}& \roman{if}\ w=y,\cr 0& \roman{else}.\endcases$$
Let $\l,\mu$ be as in (2.5.1). Under the isomorphism above we have
$$\aligned
\Oc_{\Nen'_{s_\a}}(\l,\mu)_\Xen|_{\Den^S}
&=\Oc_{\Den_{s_\a}^S}(\l,\mu),\cr&=
\sum^\infty_w(\theta_{w\l+w\mu}\xb_{w,w}+\theta_{w\l+ws_\a\mu}\xb_{w,ws_\a}),\cr
&=\sum^\infty_w(\theta_{-w\a}\xb_{w,w}+\xb_{w,ws_\a}).
\endaligned$$ Here the symbol ${\ds\sum^\infty}$ denotes an infinite sum.
Thus 2.4.7 yields
$$\aligned
\rb_S(1+t_{s_\a})&=-\rb_S(\Oc_{\Nen'_{s_\a}}(\l,\mu)),\cr
&=-(1-\Oc_{\Den_{s_\a}^S}(\a+t,0))(1-\Oc_{\Den_{s_\a}^S}(0,-\a))^{-1}
\Oc_{\Den_{s_\a}^S}(\l,\mu),\cr &=\sum^\infty_w{1-t\theta_{w\a}\over
1-\theta_{w\a}}(\xb_{w,w}-\xb_{w,ws_\a}).
\endaligned\leqno(2.5.3)$$
By 2.4.4 we have also
$$\aligned
\rb_S(x_\l)&=\rb_S(\Oc_{\Nen_e}\la\l\ra)
=\sum^\infty_w\theta_{w\l}\xb_{w,w}.
\endaligned\leqno(2.5.4)$$
Using (2.5.3) and (2.5.4) the relations are reduced to a simple
linear algebra computation which is left to the reader.

Next we
prove that $\Phi$ is surjective. First note that
$\Rb_t=\Rb^\Tt=\Rb^\It$. We have
$$\gathered
\K^{\It}(\Nen)=\ind_w\K^{\It}(\Nen_w),\quad
\Nen_w=\bigsqcup_{v\leqslant w}\buildrel\circ\over\Nen_v,\quad
\buildrel\circ\over\Nen_v=\Nen\cap(\nen\times{\buildrel\circ\over\Fen_v}).
\endgathered$$
Further, we have $\Tb$-scheme isomorphisms
$$\buildrel\circ\over\Fen_v=\nen^\op_v,\quad\buildrel\circ\over\Nen_v=
\nen_v\times\buildrel\circ\over\Fen_v=\nen.$$ In particular
$\buildrel\circ\over\Nen_v$ is an affine space. Let $\Nen'_v$ be the
Zariski closure of $\buildrel\circ\over\Nen_v$ in $\Nen$ and let
$\gb_v=\Oc_{\Nen'_v}$, regarded as an element of $\K^\It(\Nen)$. The
direct image by the inclusion $\Nen_w\subset\Nen$ identifies the
$\Rb_t$-module $\K^{\It}(\Nen_w)$ with the direct summand
$$\bigoplus_{v\leqslant w}\Rb_t\,\gb_v\subset\K^{\It}(\Nen).$$ See
\cite{CG, sec.~7.6} for a similar argument for non-affine flags. On
the other hand the PBW theorem for $\Hb$ implies that
$\Hb=\bigoplus_{w\in\tilde W}\Rb_t T_w$ as a left $\Rb_t$-module.
Set $\Hb_w=\bigoplus_{v\leqslant w}\Rb_t T_v$. We must prove that
$\Phi$ restricts to a surjective $\Rb_t$-module homomorphism
$\Hb_w\to\K^{\It}(\Nen_w)$ for each $w$. This is proved by induction
on the length $l(w)$ of $w$. More precisely this is obvious if
$l(w)=1$ and we know that
$$l(vw)=l(v)+l(w)\ \Rightarrow\ \Hb_v\,\Hb_w=\Hb_{vw},\
\K^\It(\Nen_v)\star\K^\It(\Nen_w)\subset\K^\It(\Nen_{vw}) .$$
Therefore we are reduced to prove that under the previous assumption
we have $$\gb_v\star\gb_w\equiv a_{v,w}\,\gb_{vw},$$ with $a_{v,w}$
a unit of $\Rb_t$. Here the symbol $\equiv$ means an equality modulo
lower terms for the Bruhat order. To do that we fix $S=\Tt$ and we
consider the image of $\gb_w$ by $\rb_\Sigma$. It is, of course, to
complicated to compute the whole expression, but we only need the
terms $g^{(z)}_{y,yz}$ with $l(yz)=l(y)+l(z)$ in the sum below
$$\rb_\Sigma(\gb_x)=\sum_{y,z}^\infty g^{(x)}_{y,z}\,\xb_{y,z},$$
because the coefficient $a$ above is given by the following relation
$$g^{(w)}_{v,vw}\,g^{(v)}_{e,v}=a\,g^{(vw)}_{e,vw}.$$
The same computation as in 2.4.6 shows that
$$g^{(z)}_{y,yz}=\prod_{\a\in\tilde\Delta^\op_z}{1-t\theta_{y\a}\over
1-\theta_{-y\a}}.$$ Now, recall that
$$l(vw)=l(v)+l(w)\ \Rightarrow\
\tilde\Delta_{vw}^\op=\tilde\Delta_v^\op\sqcup
v(\tilde\Delta^\op_w).$$ Thus we have $a_{v,w}=1$.

Finally, since $\Phi$ restricts to a surjective $\Rb_t$-module
homomorphism $\Hb_w\to\K^{\It}(\Nen_w)$ for each $w$ and both sides
are free $\Rb_t$-modules of rank $l(w)$ necessarily $\Phi$ is
injective.

The last claim of the theorem follows from 2.3.9$(c)$.

\qed

\vskip3mm

\subhead 2.5.7.~Remark\endsubhead By 2.3.9 the convolution
product
$$\star:\K^\It(\Nen)\otimes\K^\It(\Nen)\to\K^\It(\Nen)$$ is
$\Rb^\It$-linear in the first variable. Recall that forgetting the
group action yields an isomorphism $\K^\It(\Nen)\to\K^\Tt(\Nen)$.
Further, since $\Nen$ has a partition into affine cell a standard
argument implies that the forgetting map gives an isomorphism
$$\Rb^S\otimes_{\Rb^\Tt}\K^\Tt(\Nen)=\K^S(\Nen)$$
for each closed subgroup $S\subset \Tt$. Thus the map $\star$
factors to a group homomorphism
$$\star:\K^S(\Nen)\otimes\K^\It(\Nen)\to\K^S(\Nen).$$
The assignment $\theta_\l\mapsto
X_\l$ identifies $\Rb^{\Tt}$ with the subring $\Rb_t\subset\Hb$ generated by
$t$ and the $X_\l$'s.
By 2.5.6 the group homomorphism above is identified, via the map $\Phi$,
with the right multiplication of $\Hb$ on $\Rb^S\otimes_{\Rb_t}\Hb.$

\vskip3mm

\head 3.~Classification of the simple admissible modules of the
Double affine Hecke algebra\endhead

\subhead 3.1.~Constructible sheaves and convolution algebras\endsubhead

The purpose of this section  is to revisit the sheaf-theoretic
analysis of convolution algebras in \cite{CG, sec.~8.6} in a more
general setting including the case of schemes locally of finite
type. It is an expanded version of \cite{V, sec.~6, app.~B}.

\subhead 3.1.1.~Convolution algebras and schemes locally of finite
type\endsubhead Let $N=\bigsqcup_{\a\in A}N(\a)$ be a disjoint union
of smooth quasi-projective connected schemes. We'll assume that the
set $A$ is countable. We'll view $N$ as an ind-scheme, by setting
$$N=\ind_{B\subset A}N(B),\quad N(B)=\bigsqcup_{\a\in B}N(\a),$$
where $B$ is any finite subset of $A$. Let $C$ be a quasi-projective
scheme (possibly singular) and $\pi:N\to C$ be an ind-proper map.
For each $\a,\b\in A$ we set
$$M(\a,\b)=N(\a)\times_CN(\b)$$
(the reduced fiber product).
It is a closed subscheme of $N(\a)\times N(\b)$. The fiber product
means indeed the reduced fiber product. Note that $N(\a)$,
$M(\a,\b)$ are complex varieties which can be equipped with their
transcendental topology. The symbol $\Hb_*(.,\CC)$ will denote the
Borel-Moore homology with complex coefficients. We'll view
$M=\bigsqcup_{\a,\b}M(\a,\b)$ as an ind-scheme in the obvious way.
We set
$$\Hb_*(M,\CC)=\bigoplus_{\a,\b}\Hb_*(M(\a,\b),\CC),\quad
\widehat\Hb_*(M,\CC)=\prod_\a\bigoplus_\b\Hb_*(M(\a,\b),\CC).$$ We'll
view $\widehat\Hb_*(M,\CC)$ as a topological $\CC$-vector space in
the following way

\vskip2mm

\item{$\bullet$}
$\bigoplus_\b\Hb_*(M(\a,\b),\CC)$ is given the discrete topology for
each $\a$,

\vskip2mm

\item{$\bullet$} $\widehat\Hb_*(M,\CC)$ is given the product
topology.

\vskip2mm

\noindent We also equip $\Hb_*(M,\CC)$ with a convolution product
$\star$ as in \cite{CG, sec.~8}. The following is immediate.

\proclaim{3.1.2.~Lemma} The multiplication on $\Hb_*(M,\CC)$ is
bicontinuous and yields the structure of a topological ring on
$\widehat\Hb_*(M,\CC)$.\endproclaim

\subhead 3.1.3.~Remark\endsubhead We may also consider the K-theory
rather than the Borel-Moore homology. Since $M$, $N$ are ind-schemes
of ind-finite type we have $$\K(N)=\bigoplus_\a \K(N(\a)),\quad
\K(M)=\bigoplus_{\a,\b} \K(M(\a,\b)).$$ We'll also set
$$\widehat\K(M)=\prod_\a\bigoplus_\b \K(M(\a,\b)).$$
Thus $\widehat\K(M)$ is again a topological ring. The multiplication
in $\K(M)$, $\widehat\K(M)$ is the convolution product associated
with the inclusion $M\subset N^2$. It is defined as in (2.4.8). By
\cite{CG, thm.~5.11.11} the bivariant Riemann-Roch map yields a
topological ring homomorphism
$$RR:\CC\widehat\K(M)\to\widehat\Hb_*(M,\CC)$$
which maps $\CC\widehat\K(M(\a,\b))$ to $\Hb_*(M(\a,\b),\CC)$ for
each $\a$, $\b$. It is invertible if all $\Hb_*(M(\a,\b),\CC)$'s are
spanned by algebraic cycles.

\subhead 3.1.4.~Admissible modules over the convolution algebra\endsubhead
Let $\Dcb(C)_{\CC\text{-}\c}^b$ be the derived
category of bounded complexes of constructible sheaves of
$\CC$-vector spaces over the quasi-projective scheme $C$. Given two complexes
$\Lc$, $\Lc'$ in
$\Dcb(C)_{\CC\text{-}\c}^b$ we'll abbreviate
$$\Ext^n(\Lc,\Lc')=\Hom(\Lc,\Lc'[n]),\quad
\Ext(\Lc,\Lc')=\bigoplus_{n\in\ZZ}\Ext^n(\Lc,\Lc'),$$
where the homomorphisms
are computed in the category $\Dcb(C)_{\CC\text{-}\c}^b$.
Now, we set
$$\Cc_\a=\CC_{N(\a)}[\dim(N(\a))],\quad\Lc_\a=\pi_*(\Cc_\a),\quad
\forall\a\in A.$$
Each
$\Lc_\a$ is a semi-simple complex by the decomposition theorem.
Assume that there is a finite set $\Xc$ of irreducible perverse
sheaves over $C$ such that
$$\Lc_\a\simeq\bigoplus_{n\in\ZZ}\bigoplus_{\Sc\in\Xc}L_{\Sc,\a,n}\otimes\Sc[n],$$
where $L_{\Sc,\a,n}$ are finite-dimensional $\CC$-vector spaces. We
set
$$L_{\Sc,\a}=\bigoplus_{n\in\ZZ}L_{\Sc,\a,n},\quad
L_{\Sc}=\bigoplus_{\a\in A}L_{\Sc,\a},\quad
L=\bigoplus_{\Sc\in\Xc}L_\Sc.$$
For each complexes $\Lc,\Lc',\Lc''$ the Yoneda product is a bilinear map
$$\Ext(\Lc,\Lc')\times\Ext(\Lc',\Lc'')\to\Ext(\Lc,\Lc'').$$
By \cite{CG, lem.~8.6.1, 8.9.1} we
have an algebra isomorphism $$\widehat
\Hb_*(M,\CC)=\prod_\a\bigoplus_\b\Ext(\Lc_\a,\Lc_\b),$$
where the rhs is given the Yoneda product. We have the following
decomposition as $\CC$-vector spaces
$\widehat\Hb_*(M,\CC)\simeq R\oplus J$, where
$$R=\bigoplus_{\Sc\in\Xc}\End
(L_{\Sc}),\quad
J=\bigoplus_{\Sc,\Tc\in\Xc}\bigoplus_{n>0}\Hom(L_\Tc,L_\Sc)\otimes
\Ext^n(\Tc,\Sc).$$ Further
$J$ is a nilpotent two-sided ideal of $\widehat\Hb_*(M,\CC)$ and the
$\CC$-algebra structure on $\widehat\Hb_*(M,\CC)/J$ is the obvious
$\CC$-algebra structure on $R$. Before to explain what is the
topology on $R$ recall the following basic fact.

\subhead 3.1.5.~Definition\endsubhead $(a)$ Let
$\Ab$ be any ring and let $M$, $N$ be $\Ab$-modules.
The {\it finite topology} on
$\Hom_\Ab(M,N)$ is the linear topology for which a basis of open
neighborhoods for 0 is given by the annihilator of $M',$ for all
finite set $M'\subset M$. This is actually the topology induced on
$\Hom_\Ab(M,N)$ from $N^M$ (a product of topological spaces where
$N$ has the discrete topology).

\vskip1mm

$(b)$ If $\Ab$ is a topological ring we'll say that a right
$\Ab$-module is {\it admissible} (or {\it smooth}) if for each
element $m$ the subset $\{x\in\Ab;mx=0\}$ is open.

\vskip3mm

Now we can formulate the following lemma.

\proclaim{3.1.6.~Lemma} The two-sided ideal $J\subset\widehat
\Hb_*(M,\CC)$ is closed. The quotient topology on $\widehat
\Hb_*(M,\CC)/J$ coincides with the finite topology on $R$.
\endproclaim

Therefore, the Jacobson density theorem implies that the set of
simple admissible right representations of $R$ is
$\{L_\Sc;\Sc\in\Xc\}$, see e.g., \cite{V, sec.~B}. This yields the
following.

\proclaim{3.1.7.~Proposition} The set of the simple admissible
right $\widehat\Hb_*(M,\CC)$-modules is canonically identified
the set $\{L_\Sc;\,\Sc\in\Xc\}$.
\endproclaim

\vskip3mm

\subhead 3.2.~Simple modules in the category $\Ocb$\endsubhead

This section reviews the classification of the simple modules in
$\Ocb(\Hb)$ from \cite{V}. The main arguments are the same as in
loc.~cit., but the use of the concentration map simplifies the
exposition. Note that $\Ocb(\Hb)$ consists of {\it right} $\CC\Hb$-modules.
This specification is indeed irrelevant because the $\ZZ_{q,t}$-algebra
$\Hb$ is isomorphic to its opposit algebra, see e.g., \cite{C, thm.~1.4.4}.

\subhead 3.2.1.~From $\Ocb(\Hb)$ to modules over the convolution algebra of
$\Men$\endsubhead In this section we apply the construction from
Section 3.1 in the following setting. Fix a regular closed subgroup
$S\subset \Tt$.
Following \cite{KL} we define the set
of the {\it topologically nilpotent elements} in
$\tilde\gen$ by
$$\Nen\ien\len=\bigcup_{\ben\in\Fen}\ben_\nil.$$
Let $N=\dot\nen^S$, $C=\Nen\ien\len^S$, and let $\pi:N\to C$ be the
obvious projection. The ind-scheme $M$ in 3.1.1 is given by
$M=\Men^S$. It is an ind-scheme of ind-finite type. We'll use the
notation from 2.4.8. Recall that
$$\K(\Men^S)=\bigoplus_{\a,\b}\K(\Men(\a,\b)),\quad
\widehat\K(\Men^S)=\prod_\a\bigoplus_\b\K(\Men(\a,\b)).$$

Now we fix an element $(h,\zeta)=(s,\tau,\zeta)$ in $\Tt$, i.e.,
we have
$h=(s,\tau)\in\tilde T$, $s\in T\times\CC^\times_\cent$,
$\tau\in\CC^\times_\rot$ and $\zeta\in\CC^\times_\qua$. Assume that
$S=\la (h,\zeta)\ra$, i.e., we assume that $S$ is the closed subgroup of $\Tt$
generated by the element $(h,\zeta)$. Let $\tilde G^h$ be the centralizer of the
element $h$ in the group $\tilde G$.

\subhead 3.2.2.~Definition\endsubhead We'll say that the pair
$(\tau,\zeta)$ is {\it regular} if $\tau$ is not a root of 1 and
$\tau^k\neq\zeta^m$ for each $m,k>0$.

\vskip2mm

For each set $X$ with a $\Tt$-action we'll abbreviate
$X^{h,\zeta}$ for the fixed points subset $X^{(h,\zeta)}$.
We have the following \cite{V, lem.~2.13}, \cite{VV, lem.~2.4.1-2}.

\proclaim{3.2.3.~Proposition} Assume that the pair $(\tau,\zeta)$ is
regular. The group $\la (h,\zeta)\ra$ is regular. The group $\tilde
G^h$ is reductive and connected.
The scheme $\Nen\ien\len^{h,\zeta}$
is of finite type and it consists of nilpotent elements of
$\tilde\gen$.
Further $\Nen\ien\len^{h,\zeta}$ contains only a finite number of $\tilde
G^h$-orbits.
\endproclaim

This proposition is essentially straightforward, except for the
connexity of the reductive group $\tilde G^h$. This is an affine
analogue of a well-known result of Steinberg which says that the
centralizer of a semi-simple element in a connected reductive group
with simply connected derived subgroup is again connected. The proof
of the connexity relies on a theorem of Kac and Peterson \cite{KP}
which says that a reductive subgroup of $\tilde G$ is always
conjugated to a subgroup of a proper L\'evi subgroup of $\tilde G$.
Since the proper L\'evi subgroups of $\tilde G$ are reductive
with simply connected
derived subgroup, because $\tilde G$ is the maximal affine Kac-Moody
group, the claim is reduced to the Steinberg theorem.

Therefore, if $(\tau,\zeta)$ is regular then the scheme
$\Nen\ien\len^{h,\zeta}$ is of finite type, the scheme
$\Men^{h,\zeta}$ is locally of finite type, the homology group
$\widehat\Hb_*(\Men^{h,\zeta},\CC)$ is a topological ring by
3.1.1, and the simple admissible
right $\widehat\Hb_*(\Men^{h,\zeta},\CC)$-modules are labeled by the set
of irreducible perverse sheaves over $\Nen\ien\len^{h,\zeta}$
which occur as a shift of a direct summand of the complex
$\pi_*(\CC_{\dot\nen^{h,\zeta}}).$

Set $\Sigma=\{(h,\zeta)\}$ and $S=\la (h,\zeta)\ra$. We'll abbreviate
$$\rb_{h,\zeta}=\rb_\Sigma,\quad\Rb_{h,\zeta}=\Rb^S_\Sigma.$$
Composing $\Phi$, $\rb_{h,\zeta}$ and the tensor product by the
character
$$\chi_{h,\zeta}:\Rb_{h,\zeta}\to\CC,\quad f\mapsto f(h,\zeta),$$ we get a $\CC$-algebra homomorphism
$$\Phi_{h,\zeta}:\CC\Hb\to\CC\widehat\K(\Men^{h,\zeta}).$$
Note that
$\widehat\K(\Men^{h,\zeta})$ is a topological ring by 3.1.3
and that the bivariant Riemann-Roch map yields a topological ring
homomorphism
$$RR:\CC\widehat\K(\Men^{h,\zeta})\to\widehat\Hb_*(\Men^{h,\zeta},\CC).$$
We'll write
$$\Psi_{h,\zeta}=RR\circ\Phi_{h,\zeta}:
\CC\Hb\to\widehat\Hb_*(\Men^{h,\zeta},\CC).$$

Throughout we'll use the following notation : for any ring homomorphism
$$\phi:\Ab\to\Bb$$ and for any (left or right) $\Bb$-module $M$
let $\phi^\bullet(M)$ be the corresponding $\Ab$-module.

\proclaim{3.2.4.~Proposition} Assume that the pair $(\tau,\zeta)$ is
regular.

(a) The map
$\Phi_{h,\zeta}:\CC\Hb\to\CC\widehat\K(\Men^{h,\zeta})$ has a
dense image.

(b) The map
$RR:\CC\widehat\K(\Men^{h,\zeta})\to\widehat\Hb_*(\Men^{h,\zeta},\CC)$
is an isomorphism.

(c) The pull-back by the composed map
$\Psi_{h,\zeta}=RR\circ\Phi_{h,\zeta}$ gives a bijection
between the set of simple right $\CC\Hb$-modules in
$\Ocb_{h,\zeta}(\Hb)$ and the set of simple admissible
right $\widehat\Hb_*(\Men^{h,\zeta},\CC)$-modules.
\endproclaim

The proof of 3.2.4 is given in 3.2.7 below. Before this we need more material.

\subhead 3.2.5.~The regular representation of $\Hb$\endsubhead
First
we define a right representation of $\widehat\K(\Men^{h,\zeta})$
on $\K(\Nen^{h,\zeta})$. We'll use the same notation as in the
previous subsection. In particular $S=\la
(h,\zeta)\ra$ is a regular closed subgroup of $\Tt$. Recall that
$\dot\nen^{h,\zeta}$ and $\Men^{h,\zeta}$ are both ind-scheme
of ind-finite type, that $\dot\nen^{h,\zeta}$ is a disjoint union of smooth
quasi-projective
varieties, and that $\Men^{h,\zeta}$ is regarded as a closed subset of
$(\dot\nen^{h,\zeta})^2$. The
convolution product on $\K(\Men^{h,\zeta})$ is given by
$$x\star y=
R(q_2)_*\bigl(q_3^*(x)\Lotimes_{(\dot\nen^{h,\zeta})^3} q_1^*(y)\bigr),
\quad\forall
x,y\in\K(\Men^{h,\zeta}),$$ where
$q_a:(\dot\nen^{h,\zeta})^3\to(\dot\nen^{h,\zeta})^2$
is the projection along the
$a$-th factor for $a=1,2,3$. The inclusion
$\Nen\subset\dot\nen$ yields an inclusion of ind-schemes
$\Nen^{h,\zeta}\subset \dot\nen^{h,\zeta}$.
For each $x\in\K(\Nen^{h,\zeta})$ and
each $y\in\K(\Men^{h,\zeta})$ we define the following element in
$\K(\Nen^{h,\zeta})$
$$x\star y=R(p_1)_*(p_2^*(x)\Lotimes_{(\dot\nen^{h,\zeta})^2} y),\leqno(3.2.1)
$$
where $p_a:(\dot\nen^{h,\zeta})^2\to\dot\nen^{h,\zeta}$
is the projection along the $a$-th factor for
$a=1,2$. It is well-known that the map (3.2.1) defines a right
representation of $\K(\Men^{h,\zeta})$ on $\K(\Nen^{h,\zeta})$,
see e.g., \cite{CG}.

\proclaim{3.2.6.~Lemma} $(a)$ The right representation of
$\K(\Men^{h,\zeta})$ on $\K(\Nen^{h,\zeta})$ extends uniquely to
an admissible right representation of
$\widehat\K(\Men^{h,\zeta})$ on $\K(\Nen^{h,\zeta})$.

$(b)$ The right $\CC\Hb$-module
$\chi_{h,\zeta}\otimes_{\Rb_t}\Hb$ belongs to $\Ocb_{h,\zeta}(\Hb)$.

$(c)$ There is
an isomorphism of right $\Hb$-modules
$\chi_{h,\zeta}\otimes_{\Rb_t}\Hb\simeq\Phi_{h,\zeta}^\bullet(\CC\K(\Nen^{h,\zeta}))$.
\endproclaim

\noindent{\sl Proof :} The first claim is obvious, because we have
$$\gathered\widehat\K(\Men^{h,\zeta})=\prod_\a\bigoplus_\b\K(\Men(\a,\b)),\quad
\K(\Nen^{h,\zeta})=\bigoplus_\a\K(\Nen(\a)),\cr
\K(\Nen(\a))\star\K(\Men(\a,\b))\subset\K(\Nen(\b)),\endgathered$$
where $\Nen(\a)=\Nen\cap\dot\nen(\a).$ Part $(b)$ is a standard computation,
see e.g., \cite{V}. Let us concentrate on part $(c)$.
Composing the map
$\chi_{h,\zeta}:\Rb_{h,\zeta}\to\CC$
with the canonical map $\Rb^\Tt\to\Rb_{h,\zeta}$ we
may regard $\chi_{h,\zeta}$ as the one-dimensional $\Rb^\Tt$-module
given by $f\mapsto f(h,\zeta)$. Recall that $\Rb_t=\Rb^\Tt$, see 2.5.7.
The vector
space $\chi_{h,\zeta}\otimes_{\Rb_t}\Hb$ has an obvious structure of
right $\Hb$-module.
The isomorphism 2.5.6 factors to a right $\Hb$-module
isomorphism
$$\chi_{h,\zeta}\otimes_{\Rb_t}\Hb\to
\chi_{h,\zeta}\otimes_{\Rb^\Tt}\K^{\Tt}(\Nen).$$
We claim that there is a right $\Hb$-module isomorphism
$$\chi_{h,\zeta}\otimes_{\Rb^\Tt}\K^{\Tt}(\Nen)\to\Phi_{h,\zeta}^\bullet(\CC\K(\Nen^{h,\zeta})).$$
To prove this, recall that composing the maps
$\rb_{h,\zeta}$ and $\chi_{h,\zeta}$
yields an algebra homomorphism
$$\K^\Tt(\Nen)\to
\CC{\K}(\Nen_\Xen^{h,\zeta})=\CC{\K}(\Men^{h,\zeta}).$$
Thus we must construct a map
$\rb:\K^\Tt(\Nen)\to
\CC{\K}(\Nen^{h,\zeta})$
which intertwines the right $\star$-product of
$\K^\Tt(\Nen)$ on itself, see 2.3.7,
with the right $\star$-product of
$\CC\K(\Men^{h,\zeta})$ on $\CC\K(\Nen^{h,\zeta})$,
see (3.2.1), relatively to the ring homomorphism
$$\chi_{h,\zeta}\circ\rb_{h,\zeta}:\K^\Tt(\Nen)\to
\CC{\K}(\Men^{h,\zeta}).$$
Further the map $\rb$ should factor to an isomorphism
$$\chi_{h,\zeta}\otimes_{\Rb^\Tt}\K^\Tt(\Nen)\to\CC\K(\Nen^{h,\zeta}).$$
Consider the following chain of inclusions
$$\xymatrix{\nen^{h,\zeta}\times\Fen^{h,\zeta}\ar[r]^i&
\nen^{h,\zeta}\times\Fen\ar[r]^j&\nen\times\Fen.}\leqno(3.2.2)$$
Since $\nen$ is pro-smooth we can consider the map $Lj^*$ in K-theory,
see 1.5.18.
Since $\Fen$ is an ind-$S$-scheme of ind-finite type we can consider the map
$i_*$ in K-theory, see 1.5.19.
Both maps are invertible, and the composed map is an isomorphism
$$(i_*)^{-1}\circ Lj^*:\K^S(\nen\times\Fen)_\Sigma
\to\K^S(\nen^{h,\zeta}\times\Fen^{h,\zeta})_\Sigma,\quad\Sigma=\{(h,\zeta)\}.$$
Now, recall that we have a good embedding
$\Nen\subset\nen\times\Fen$. Thus, we obtain also in this way an isomorphism
$\K^S(\Nen)_\Sigma \to\CC\K(\Nen^{h,\zeta}).$
Composing it with the obvious map
$\K^\Tt(\Nen)\to\K^S(\Nen)_\Sigma$
it yields a map
$$\rb:\K^\Tt(\Nen)\to\CC\K(\Nen^{h,\zeta}).$$
We must check that the map $\rb$ is compatible with the right
$\star$-product, in the above sense. This is left to the reader. The
proof is the same as the proof of 2.4.9. Compare (2.3.3), (2.4.3)
with (3.2.1), (3.2.2).

\qed

\subhead 3.2.7.~Proof of 3.2.4\endsubhead $(a)$ The map
$\Phi:\Hb\to\K^{\It}(\Nen)$ is invertible by 2.5.6. The composed map
$$\xymatrix{\K^\It(\Nen)\ar[r]^-{\rb_{h,\zeta}}&
{\K}^{S}(\Men^{h,\zeta})_{h,\zeta}=
\Rb_{h,\zeta}\otimes{\K}(\Men^{h,\zeta})
\ar[r]^-{\chi_{h,\zeta}}&
\CC{\K}(\Men^{h,\zeta})\subset\CC{\widehat\K}(\Men^{h,\zeta})}$$ has a
dense image, because the image contains $\CC\K(\Men(\a,\b))$
for each $\a$, $\b$. Composing both maps we get $\Phi_{h,\zeta}$.
Thus $\Phi_{h,\zeta}$ has a dense image.

\vskip1mm

$(b)$ It is easy to see that $\Hb_*(\Men(\a,\b),\CC)$ is spanned by
algebraic cycles for all $\a$, $\b$. Therefore we have
$\CC\widehat\K(\Men^{h,\zeta})\simeq\widehat\Hb_*(\Men,\CC)$.

\vskip1mm

$(c)$ By part $(b)$ it is enough to check that the map
$\Phi_{h,\zeta}^\bullet$ yields a bijection between the set of
simple objects in $\Ocb_{h,\zeta}(\Hb)$ and the set of simple
admissible right representations of the topological $\CC$-algebra
$\CC\widehat\K(\Men^{h,\zeta})$.
Our proof uses the following lemma, which will be checked later on.

\proclaim{3.2.8.~Lemma} (a) For each $\l\in\tilde\Xb$ the operator
of right multiplication by $\Phi_{h,\zeta}(X_\l)$ in any admissible
right $\CC\widehat\K(\Men^{h,\zeta})$-module is locally finite and
its spectrum belongs to the set $\{{}^w\l(h);w\in\tilde W\}.$

(b) If the elements $h,h'\in\tilde T$ are $\tilde W$-conjugate then
the topological rings $\CC\widehat\K(\Men^{h,\zeta})$,
$\CC\widehat\K(\Men^{h',\zeta})$ and the homomorphisms
$\Phi_{h,\zeta}$, $\Phi_{h',\zeta}$ are canonically identified.

\endproclaim

The claim 3.2.4$(c)$ is a corollary of 3.2.6 and 3.2.8.
First, let $M$ be a simple admissible right
$\CC\widehat\K(\Men^{h,\zeta})$-module. The right $\Hb$-module
$\Phi_{h,\zeta}^\bullet(M)$ belongs to $\Ocb_{h,\zeta}(\Hb)$ by 3.2.8$(a)$.
Further $\Phi_{h,\zeta}^\bullet(M)$ is a simple right $\Hb$-module.
Indeed, since $\Phi_{h,\zeta}(\CC\Hb)$ is dense in
$\CC\widehat\K(\Men^{h,\zeta})$ by 3.2.4$(a)$ and
since $M$ is admissible and simple
as a right $\CC\widehat\K(\Men^{h,\zeta})$-module, we have
$$x\star\Phi_{h,\zeta}(\CC\Hb)=x\star\CC\widehat\K(\Men^{h,\zeta})=M,\quad
\forall\, 0\neq x\in M.$$ Thus $M$ is a simple object of
$\Ocb_{h,\zeta}(\Hb)$.

Next, let $L$ be a simple object of $\Ocb_{h,\zeta}(\Hb)$. We claim
that there is a simple admissible right
$\CC\widehat\K(\Men^{h,\zeta})$-module $M$ such that
$L\simeq\Phi_{h,\zeta}^\bullet(M)$. Indeed, since $L$ belongs to
$\Ocb_{h,\zeta}(\Hb)$ there is an element $h'\in\tilde W\cdot h$
such that the $h'$-weight subspace $L_{h'}$ is non-zero. Since $L$
is simple, it is therefore a quotient of the right $\CC\Hb$-module
$\chi_{h',\zeta}\otimes_{\Rb_t}\Hb$. The later is isomorphic to
$\Phi_{h',\zeta}^\bullet(\CC\K(\Nen^{h',\zeta}))$ by 3.2.6. Let
$J$ be the kernel of the quotient map
$\Phi_{h',\zeta}^\bullet(\CC\K(\Nen^{h',\zeta}))\to L$. Hence, $J$ is a right
$\Phi_{h',\zeta}(\CC\Hb)$-submodule of $\CC\K(\Nen^{h',\zeta})$. Hence
it is also a right $\CC\widehat\K(\Men^{h',\zeta})$-module because
$\Phi_{h',\zeta}(\CC\Hb)\subset\CC\widehat\K(\Men^{h',\zeta})$ is dense
and $\CC\K(\Nen^{h',\zeta})$ is admissible.
By 3.2.8$(b)$ we can regard
$\CC\K(\Nen^{h',\zeta})$
as a right
$\CC\widehat\K(\Men^{h,\zeta})$-module and
$J$ as a right
$\CC\widehat\K(\Men^{h,\zeta})$-submodule of
$\CC\K(\Nen^{h',\zeta})$. Then the quotient
$\CC\K(\Nen^{h',\zeta})/J$ is again a right
$\CC\widehat\K(\Men^{h,\zeta})$-submodule and we have
$$L\simeq\Phi_{h,\zeta}^\bullet(\CC\K(\Nen^{h',\zeta})/J)$$
as right $\Hb$-modules. Further, since $L$ is a simple right
$\CC\Hb$-module the quotient $\CC\K(\Nen^{h',\zeta})/J$ is a
simple admissible right $\CC\widehat\K(\Men^{h,\zeta})$-module.

Finally if $M$, $M'$ are admissible right
$\CC\widehat\K(\Men^{h,\zeta})$-modules such that
$\Phi_{h,\zeta}^\bullet(M)$, $\Phi^\bullet_{h,\zeta}(M')$ are isomorphic as
right $\Hb$-modules then
$M$, $M'$ are isomorphic as right $\CC\widehat\K(\Men^{h,\zeta})$-modules,
because they are isomorphic as right
$\Phi_{h,\zeta}(\CC\Hb)$-modules and
$\Phi_{h,\zeta}(\CC\Hb)$ is a dense subring of
the topological ring $\CC\widehat\K(\Men^{h,\zeta})$.

\qed

\vskip3mm

\noindent{\sl Proof of 3.2.8 :} $(a)$ For each $\l\in\tilde\Xb$ we have
$$\Phi_{h,\zeta}(X_\l)=
x_\l=\Oc_{\Nen_e}\la\l\ra=\Oc_{\Nen_e}(\l)\in\K^\It(\Nen).$$ Note
that the set $\Men_e(\a,\b)=\Men_e\cap\Men(\a,\b)$ is empty if
$\a\neq\b$ and that it is the diagonal of $\dot\nen(\a)$ else.
Recall that $S=\la (h,\zeta)\ra$. For each $\a$ let
$\l_\a:S\to\CC^\times$ be the character of the group $S$ such that
any element $g\in S$ acts on the equivariant line bundle
$\Oc_{\Fen(\a)}(\l)$ by fiberwise multiplication by the scalar
$\l_\a(g)$. It is well-known that for each $\a$ there is an element
$w\in\tilde W$ such that $\l_\a=({}^w\l)|_S$. By 2.4.5 we have
$$\Phi_{h,\zeta}(X_\l)=\sum^\infty_\a\l_\a(h)\,\Oc_{\Men_e(\a,\a)}(\l,0)=
\sum^\infty_\a\l_\a(h)\,\Oc_{\Men_e(\a,\a)}(0,\l)\in\CC\widehat\K(\Men^{h,\zeta}).$$
Thus the operator of multiplication by $\Phi_{h,\zeta}(X_\l)$ in any
admissible $\CC\widehat\K(\Men^{h,\zeta})$-module is locally
finite and its spectrum belongs to the set $\{{}^w\l(h);w\in\tilde
W\}.$ See \cite{V, lem.~4.8} for details.

$(b)$ Since $h$ and $h'$ are $\tilde W$-conjugate they are also
$\tilde G$-conjugate. The group $\tilde G$ acts on $\Men$. This
yields an ind-scheme isomorphism
$\Men^{h,\zeta}\simeq\Men^{h',\zeta}$. The rest of the claim is
obvious.

\qed

\vskip3mm

\subhead 3.2.9.~The classification theorem\endsubhead
We can now compose 3.1.7 with 3.2.4$(c)$. We get the following
theorem \cite{V, thm.~7.6}, \cite{VV, prop.~2.5.1}
whose proof uses the connexity of the reductive
group $\tilde G^h$ in 3.2.3.
To state the theorem we need more material.
Assume that the pair $(\tau,\zeta)$ is regular.
As above, we'll write $S=\la(h,\zeta)\ra$.
Let $\Xc_{h,\zeta}$ be the set of irreducible perverse
sheaves over $\Nen\ien\len^{h,\zeta}$
which are direct summand (up to some shift)
of the complex $$\bigoplus_\a(\pi_{h,\zeta})_*\CC_{\dot\nen(\a)},\quad
\pi_{h,\zeta}:\dot\nen^{h,\zeta}=
\bigsqcup_\a\dot\nen(\a)\to\Nen\ien\len^{h,\zeta}.$$
Here the map $\pi_{h,\zeta}$ is the obvious projection.
There is a finite number of $\tilde G^h$-orbits in $\Nen\ien\len^{h,\zeta}$.
For each closed point $x\in\Nen\ien\len^{h,\zeta}$
let $A(h,\zeta,x)$ be the
group of connected components of the isotropy subgroup of $x$ in $\tilde G^h$.
The group $A(h,\zeta,x)$ acts in an obvious way on the homology space
$$H_*(\pi^{-1}_{h,\zeta}(x),\CC)=
\bigoplus_\a H_*(\pi^{-1}_{h,\zeta}(x)\cap\dot\nen(\a),\CC).$$
Let $\Irr(A(h,\zeta,x))$ be the set of irreducible representations of
the finite group $A(h,\zeta,x)$.
Each representation in $\Irr(A(h,\zeta,x))$ can be regarded
as a $\tilde G^h$-equivariant
irreducible local system over the $\tilde G^h$-orbit $O$ of $x$.
Therefore we may regard $\Xc_{h,\zeta}$ as a set of pairs
$(x,\chi)$ in $\bigsqcup_x\Irr(A(h,\zeta,x)).$

\proclaim{3.2.10.~Theorem} Assume that $(\tau,\zeta)$ is regular.

$(a)$ The set $\{\Psi_{h,\zeta}^\bullet(L_\Sc);\Sc\in\Xc_{h,\zeta}\}$
is the set of all simple objects in $\Ocb_{h,\zeta}(\Hb)$.

$(b)$ The set $\Xc_{h,\zeta}$  is identified with the
set of pairs $(x,\chi)$ such that $\chi\in\Irr(A(h,\zeta,x))$
is a Jordan-H\"older factor of the
$A(h,\zeta,x)$-module
$H_*(\pi^{-1}_{h,\zeta}(x),\CC).$

$(c)$ The simple right $\Hb$-modules
$\Psi_{h,\zeta}^\bullet(L_{x,\chi})$
and
$\Psi_{h',\zeta}^\bullet(L_{x',\chi'})$
are isomorphic iff the triplets $(h,x,\chi)$ and $(h',x',\chi')$
are $\tilde G$-conjugate.

\endproclaim

\enddocument